%% file: main.tex
\documentclass[a4paper, 11pt]{article}
\usepackage{cmap}
\usepackage[english]{babel}
\usepackage[sc,osf]{mathpazo}
\usepackage[T1]{fontenc}
\usepackage[utf8]{inputenc}
\usepackage{microtype}
\usepackage{amsfonts}       
\usepackage{amssymb}       	
\usepackage{eucal}
\usepackage{mathrsfs}
\usepackage{amsmath}        
\usepackage{amsthm}         
\usepackage{mathtools}      
\usepackage{stmaryrd}
\usepackage{thmtools}
\usepackage{graphicx}       
\usepackage[dvipsnames]{xcolor}
\usepackage{accents}
\usepackage[shortlabels]{enumitem}
\usepackage[colorlinks=true,citecolor=PineGreen,linkcolor=BrickRed]{hyperref}
\usepackage[nameinlink,capitalize]{cleveref}
\usepackage{tikz-cd}
\usetikzlibrary{shapes.geometric,patterns,arrows,decorations.pathreplacing}
\usepackage{caption}
\usepackage{subcaption}
\usepackage{subfiles}

\usepackage{csquotes}
\usepackage[backend=biber,style=alphabetic,sorting=nty,giveninits]{biblatex} 
\renewbibmacro{in:}{        
	\ifentrytype{article}{}{\printtext{\bibstring{in}\intitlepunct}}}
\DeclareFieldFormat{url}{Available at \url{#1}}
\AtEveryBibitem{%
	\savefield{note}{\myadd}%
	\clearfield{note}%
	\newbibmacro{finentry}{
		\restorefield{note}{\myadd}\printfield{note}
		\finentry}%
}

\tikzset{%
    symbol/.style={%
        draw=none,
        every to/.append style={%
            edge node={node [sloped, allow upside down, auto=false]{$#1$}}}
    }
}

\tikzset{
    rot90/.style={anchor=south, rotate=90, inner sep=.5mm}
}
\tikzset{
    rot320/.style={anchor=south, rotate=320, inner sep=.5mm}
}
\tikzset{
    rot45/.style={anchor=south, rotate=45, inner sep=.5mm}
}

\addto\extrasenglish{ 
 
 }

\crefformat{equation}{#2(#1)#3}
\Crefformat{equation}{#2(#1)#3}

\declaretheorem[name=Theorem, numberwithin=section]{theorem}
\declaretheorem[name=Theorem]{theoremA}

\declaretheorem[name=Corollary, sibling=theoremA]{corollaryA}
\declaretheorem[name=Lemma, sibling=theorem]{lemma}
\declaretheorem[name=Proposition, sibling=theorem]{proposition}
\declaretheorem[name=Corollary, sibling=theorem]{corollary}
\declaretheorem[name=Corollary, numbered=no]{corollary*}

\declaretheorem[style=definition, name=Definition, sibling=theorem]{definition}
\declaretheorem[style=definition, name=Remark, sibling=theorem]{remark}
\declaretheorem[style=definition, name=Example, sibling=theorem]{example}
\declaretheorem[style=definition, name=Construction, sibling=theorem]{construction}
\declaretheorem[style=definition, name=Assumption, sibling=theorem]{assumption}

\DeclareMathOperator{\Fun}{Fun}
\DeclareMathOperator{\Funho}{Fun^{\mathrm{ho}}}
\DeclareMathOperator{\Funred}{Fun^{\mathrm{red}}}
\DeclareMathOperator{\Funhored}{Fun^{\mathrm{hored}}}
\DeclareMathOperator{\Funlin}{Fun^{\mathrm{lin}}}
\DeclareMathOperator{\Funmred}{Fun^{\mathrm{mred}}}

\DeclareMathOperator{\Funsym}{Fun^{\mathrm{sym}}}
\DeclareMathOperator{\Funsymred}{Fun^{\mathrm{symred}}}

\DeclareMathOperator{\Funst}{Fun^\mathrm{st}}
\DeclareMathOperator{\Nat}{Nat}
\DeclareMathOperator{\Lan}{Lan}
\DeclareMathOperator{\Ran}{Ran}

\DeclareMathOperator*{\hocolim}{hocolim}
\DeclareMathOperator{\id}{id}

\DeclareMathOperator{\Ob}{Ob}

\DeclareMathOperator{\Map}{Map}

\DeclareMathOperator{\Ho}{Ho}
\DeclareMathOperator{\Sing}{Sing}

\newcommand{\wh}{\widehat}
\newcommand{\wt}{\widetilde}

\newcommand{\comp}{c}
\newcommand{\ac}{\mathsf{ac}}
\newcommand{\unit}{\mathbb 1}

\newcommand{\bbN}{\mathbb N}

\newcommand{\calE}{\mathcal E}
\newcommand{\calU}{\mathcal U}
\newcommand{\calV}{\mathcal V}
\newcommand{\calW}{\mathcal W}
\newcommand{\calN}{\mathcal N}

\newcommand{\scrW}{\mathscr W}

\newcommand{\SSS}{\mathbb S}

\newcommand{\s}{\mathbf s}

\newcommand{\Cin}{\mathcal C}
\newcommand{\Din}{\mathcal D}

\newcommand{\modN}{\mathcal N}
\newcommand{\modM}{\mathcal M}

\newcommand{\Set}{\mathbf{Set}}
\newcommand{\Cat}{\mathbf{Cat}}

\newcommand{\Top}{\mathbf{Top}}

\newcommand{\pfinCW}{\mathscr{W}}
\newcommand{\finsSet}{\s\Set^\mathrm{fin}}
\newcommand{\pfinsSet}{\s\Set^\mathrm{fin}_*}

\newcommand{\OIndex}{\mathcal{I}_S}
\newcommand{\SSSOIndex}{\underline{\mathcal{I}}_S}

\newcommand{\Sp}{\mathbf{Sp}}

\makeatletter
\newcommand{\oset}[3][0ex]{%
  \mathrel{\mathop{#3}\limits^{
    \vbox to#1{\kern-2\ex@
    \hbox{$\scriptstyle#2$}\vss}}}}
\makeatother

\newcommand{\cofarrow}{\rightarrowtail}
\newcommand{\trivcofarrow}{\oset[-.5ex]{\sim}{\rightarrowtail}}
\newcommand{\wearrow}{\oset[-.16ex]{\sim}{\longrightarrow}}
\newcommand{\wearrowback}{\oset[-.16ex]{\sim}{\longleftarrow}}
\newcommand{\fibarrow}{\twoheadrightarrow}
\newcommand{\trivfibarrow}{\oset[-.5ex]{\sim}{\twoheadrightarrow}}

\renewcommand{\subset}{\subseteq}

\mathchardef\mhyphen="2D

\title{Replacing functors with enriched ones}
\author{Thomas Blom}
\date{\today}

\begin{document}
\maketitle
\begin{abstract}
	We describe simple criteria under which a given functor is naturally equivalent to an enriched one. We do this for several bases of enrichment, namely (pointed) simplicial sets, (pointed) topological spaces and orthogonal spectra. We also describe a few corollaries, such as a new proof of a result of Lurie on Dwyer--Kan localizations.
\end{abstract}


\subfile{Sections/Introduction}

\subfile{Sections/Preliminaries}

\subfile{Sections/HomotopyFunctors}

\subfile{Sections/ReducedFunctors}

\subfile{Sections/LinearFunctors}

\printbibliography

\noindent{\sc Max Planck Institute For Mathematics, Bonn, Germany.}\\
\noindent{\emph{E-mail:} \href{mailto:blom@mpim-bonn.mpg.de}{\nolinkurl{blom@mpim-bonn.mpg.de}}}\vspace{2ex}
\end{document}

%% file: Sections/Introduction.tex
\section{Introduction}

In this paper, we will investigate the question of when a given functor is equivalent to one that respects a certain enrichment. This question is motivated by examples of the following kind: An important result in Goodwillie's calculus of functors is that any linear functor $L \colon \Sp_\mathrm{fin} \to \Sp$ is equivalent to the functor $X \mapsto X \otimes C_L$, where $C_L$ is the image of the sphere spectrum $\SSS$ under $L$ and $\otimes$ denotes the smash product of spectra. To prove this result, one notes that the functors $L$ and $- \otimes C_L$ both preserve finite colimits and shifts. Since every finite spectrum is obtained from the sphere spectrum via finite colimits and shifts, the result follows if one can construct a natural transformation $X \otimes C_L \to L(X)$ extending the identification $C_L = L(\SSS)$. Such a natural transformation can be constructed explicitly, but in a more general setting where the source or target of $L$ is different, writing down an explicit description of such a comparison map may become unwieldy. However, if one were to know that $L$ is a map of $\Sp$-enriched categories, then by noting that $- \otimes C_L$ is an $\Sp$-enriched left Kan extension in the diagram
\begin{equation*}\begin{tikzcd}[column sep = large, row sep = large]
	* \arrow[d, "\SSS"'] \arrow[r, "C_L"] & \Sp \\
	\Sp_\mathrm{fin} \arrow[ru, "- \otimes C_L"{xshift=0.3ex,yshift=-0.3ex},""{name=U,below}] \arrow[ru, "L"',""{name=D}, bend right=49] &    
	\ar[Rightarrow, from=U,to=D,"\exists !"]
\end{tikzcd}\end{equation*}
one obtains the desired comparison from the universal property of enriched left Kan extensions.  Unfortunately, if one works with $\Sp$-enriched model categories and one uses Goodwillie's constructions of the $n$-excisive and the $n$-reduced approximations \cite{Goodwillie2003CalculusIII}, then the (multi)linear functors $L$ that one obtains are generally not $\Sp$-enriched. If, however, one can show that the obtained functor $L$ is weakly equivalent to an enriched one, then this still suffices to obtain the desired natural equivalence $- \otimes C_L \simeq L$. This illustrates that if one knows that a given functor is equivalent to an enriched one, then one can often avoid explicit constructions to obtain maps and use universal properties instead, leading to cleaner and more conceptual arguments.

In this paper, we will describe criteria that ensure that a functor is equivalent to an enriched one for several bases of enrichment $\calV$: we will first consider the cases where $\calV$ is the category of simplicial sets or topological spaces, we will then consider their pointed counterparts $\calV = \s\Set_*, \Top_*$, and finally, we will consider the case where $\calV$ is the category of orthogonal spectra. In the topological case, our main result is the following.

\begin{theoremA}\label{TheoremA}
	Let $\Cin$ be a small topological category and $\modN$ a good topological model category, and assume either that $\Cin$ admits tensors by the unit interval or that $\Cin$ admits cotensors by the unit interval. Then an unenriched functor $\Cin \to \modN$ is equivalent to a $\Top$-functor if and only if it sends homotopy equivalences to weak equivalences.
\end{theoremA}

By a topological (model) category, we mean a $\Top$-enriched (model) category. For the analogous result in the context of simplicial categories, we refer the reader to \Cref{theorem:MainTheoremSimplicial1}. The definition of a \emph{good} topological model category is somewhat technical and given in \Cref{def:TopologicalGoodness} below. However, most topological model categories that one encounters in nature satisfy the assumptions of \Cref{def:TopologicalGoodness},\footnote{In fact, the author is not aware of any natural example of a topological model category that is not good in the sense of \Cref{def:TopologicalGoodness}.} meaning that \Cref{TheoremA} applies to them. Moreover, \Cref{TheoremA} can be strengthened to produce Quillen equivalences between certain model categories of unenriched and enriched functors, see \Cref{theorem:MainTheoremTopological,theorem:MainTheoremSimplicial2}.

As a byproduct of the simplicial version of this theorem, we obtain a result on the existence of simplicial (co)fibrant replacement functors in general simplicial model categories (see \Cref{cor:SimplicialCofibrantReplacment}) and an alternative proof of a result of Lurie on simplicial localizations.

\begin{corollaryA}[{\cite[Prop.\ 1.3.4.7]{Lurie2017HigherAlgebra}}]\label{TheoremB}
	Let $\Cin$ be a small simplicial category that either admits tensors by $\Delta[1]$ or cotensors by $\Delta[1]$. Then the simplicial localization of the underlying category of $\Cin$ at the simplicial homotopy equivalences is DK-equivalent to $\Cin$.
\end{corollaryA}

An interesting corollary of this result is the fact that the localization of the underlying category of the category of finite simplicial sets $\finsSet$ at the set of simplicial homotopy equivalences returns the category $\finsSet$ with its usual simplicial enrichment.\footnote{The simplicial hom-sets of $\finsSet$ generally have the wrong homotopy type; that is, the simplicial hom-set between two finite simplicial sets usually does not agree with the space of maps between the corresponding spaces. The fact that one can still recover the homotopy types of these ``wrong'' mapping spaces by localizing at the ``wrong'' maps strikes the author as remarkable.}

After treating the simplicial and topological case, we study which functors can be replaced by ones that respect an enrichment in $\s\Set_*$ or $\Top_*$. The main result is that, under mild hypotheses, a functor is equivalent to a $\Top_*$- or $\s\Set_*$-enriched one if and only if it is reduced; that is, it preserves the zero object up to weak equivalence. For precise statements, see \Cref{theorem:MainTheoremReducedSimplicial1,theorem:MainTheoremReducedTopological}. As a corollary of these results, we obtain new symmetric monoidal model categories of spectra in \Cref{example:UnenrichedLydakisModelStructure,example:UnenrichedWSpectra}: we obtain analogues of Lydakis's stable category of simplicial functors \cite{Lydakis1998SimplicialFunctors} and of the $\scrW$-spectra of \cite{MMSS2001DiagramSpectra} where one works with unenriched functors as opposed to $\s\Set_*$- or $\Top_*$-enriched functors.

Finally, we turn to the case of functors enriched in orthogonal spectra. We show that one can replace a $\Top_*$-enriched functor with an $\Sp$-enriched functor if and only if it satisfies a certain linearity condition. We call an $\Sp$-enriched (model) category a spectral (model) category.

\begin{theoremA}\label{TheoremC}
	Let $\Cin$ be a small spectral category that admits cotensors by $\SSS^{-1}$, let $\modN$ be a good spectral model category and assume that the hom-spaces of the underlying $\Top_*$-category of $\Cin$ are nondegenerately based. For any $\Top_*$-enriched functor $F \colon \Cin \to \modN$, there exists a canonical map
	\begin{equation*}\mathbb{L} \Sigma F(c) \to F(c^{\SSS^{-1}}), \end{equation*}
	and $F$ is naturally equivalent to an $\Sp$-enriched functor if and only if this map is a weak equivalence for every $c$ in $\Cin$.
\end{theoremA}

Here $\mathbb{L}\Sigma$ denotes the left derived functor of the suspension functor, and the canonical map $\mathbb{L} \Sigma F(c) \to F(c^{\SSS^{-1}})$ is constructed in \Cref{construction:Sigma1}. One can view the condition that the map $\mathbb{L} \Sigma F(c) \to F(c^{\SSS^{-1}})$ is a weak equivalence as a linearity condition: since cotensors by $\SSS^{-1}$ are a type of suspension, this map being a weak equivalence can be interpreted as saying that $F$ preserves suspensions. The definition of a good spectral model category is given in \Cref{def:SpectralGoodness} below; like in the topological case, it is a condition that is satisfied in many examples that one encounters in nature. Finally, let us also mention that we upgrade the statement of \Cref{TheoremC} to a Quillen equivalence in \Cref{theorem:MainTheoremLinearOrthogonal2}. It is worth pointing out that in the context of symmetric spectra, Pereira proved a similar result in \cite[Thm.\ 10.15]{Pereira2013SpectralOperad} for functors between categories of modules over commutative algebra objects in symmetric spectra.

We obtain an analogous result about spectral functors $\Cin \to \modN$ when $\Cin$ admits \textbf{tensors} with $\SSS^{-1}$ instead of cotensors in \cref{theorem:MainTheoremLinearOrthogonal1Tensored}, although under more restrictive conditions. As an illustration of how to combine the various main results, we also prove the following.

\begin{corollaryA}
	An unenriched functor $\Sp_{\mathrm{fin}} \to \Sp$ is naturally equivalent to a spectral functor if and only if it preserves stable equivalences and is linear.
\end{corollaryA}

In Goodwillie calculus, an important role is played by symmetric multilinear functors. Since the intended applications of this work lie in Goodwillie calculus, we also deduce analogues of our main results for (symmetric) functors of multiple variables in \Cref{ssec:MultireducedFunctors,ssec:MultilinearFunctors}.

While there are substantial differences in the details of the proofs of our main results, they all follow the same general strategy. Namely, in all cases we use (derived) enriched left Kan extensions to replace functors with enriched ones, which we construct either using bar constructions or enriched coends. Our results are then deduced through a careful analysis of these bar constructions and enriched coends.

\subsection*{Overview of the paper}
We start by discussing some preliminaries on enriched (model) category theory and the bar construction in \Cref{sec:Preliminaries}. Most results of this section are well-known, and we advise the reader to only skim this section at a first reading and refer back to it as necessary while reading the rest of the paper. 

We then prove \Cref{TheoremA} and its simplicial version, namely \Cref{theorem:MainTheoremSimplicial1}, in \Cref{sec:HomotopyFunctors}. Here we also prove \Cref{TheoremB}, our result on simplicial localizations. The reader who is solely interested in this result can skip most of this paper and restrict their attention to \Cref{ssec:SimplicialFunctors,ssec:DwyerKan} and the necessary parts of \Cref{sec:Preliminaries}.

In \Cref{sec:ReducedFunctors}, we prove our main results in the setting of $\s\Set_*$- and $\Top_*$-enriched functors. We treat the case of (symmetric) functors of several variables separately in \Cref{ssec:MultireducedFunctors}, since this requires some extra work.

Finally, we treat the case of spectrally enriched functors in \Cref{sec:LinearFunctors}. We prove \Cref{TheoremC} in \Cref{ssec:LinearOrthogonalSpectraCase} and prove an analogous result where the indexing category has tensors with $\SSS^{-1}$ in \cref{ssec:LinearOrthogonalSpectraCaseTensored}. In \cref{ssec:MultilinearFunctors},  we generalize these results to the case of (symmetric) functors of several variables.

\subsection*{Notation and terminology}

Throughout this paper, whenever an adjunction is depicted as
\begin{equation*}F : \Cin \rightleftarrows \Din : G,\end{equation*}
the left adjoint is always on the left.

Let $(\calV, \otimes, \unit)$ denote a closed symmetric monoidal category. For brevity, we will call a category enriched in $\calV$ a $\calV$-category. Similarly, $\calV$-enriched functors will be called $\calV$-functors and $\calV$-enriched natural transformations will be called $\calV$-natural transformations.

Let $\Cin$ and $\Din$ be $\calV$-categories together with a $\calV$-functor $F \colon \Cin \to \Din$, let $x$, $y$ and $z$ be objects of $\Cin$ and let $v$ be an object of $\calV$. We will use the following notation and terminology:
\begin{itemize}[noitemsep]
    \item The hom-objects of $\Cin$ are denoted by $\Cin(x,y)$.
    \item Composition in $\Cin$ is denoted by $\comp \colon \Cin(y,z) \otimes \Cin(x,y) \to \Cin(x,z)$.
    \item The tensor (if it exists) of $x$ by $v$ is denoted $v \otimes x$.
    \item The cotensor (if it exists) of $x$ by $v$ is denoted $x^v$.
    \item In case the relevant tensor in $\Din$ exists, we will write $\ac \colon \Cin(x,y) \otimes F(x) \to F(y)$ for the adjunct of $\Cin(x,y) \to \Din(F(x), F(y))$ and call it the \emph{action} of $\Cin(x,y)$ on $F(x)$.
    \item The category of $\calV$-functors $\Cin \to \Din$ and $\calV$-natural transformations between them is denoted $\Fun(\Cin, \Din)$. Note that if $\calV$ is complete and $\Cin$ is small, then $\Fun(\Cin, \Din)$ is (canonically) a $\calV$-category.
\end{itemize}
Throughout this paper, we will mainly consider the following symmetric monoidal categories as bases of enrichment:
\begin{itemize}[noitemsep]
    \item The category $\Top$ of (compactly generated weak Hausdorff) topological spaces, endowed with the Cartesian product. A $\Top$-category will be called a \emph{topological} category and a $\Top$-functor will be called \emph{continuous}.
    \item The category $\s\Set$ of simplicial sets, endowed with the Cartesian product. We will call $\s\Set$-categories and $\s\Set$-functors \emph{simplicial} categories and \emph{simplicial} functors, respectively.
    \item The category $\Top_*$ of pointed topological spaces, endowed with the smash product.
    \item The category $\s\Set_*$ of pointed simplicial sets, endowed with the smash product.
    \item The category $\Sp$ of orthogonal spectra, endowed with the Day convolution product. An $\Sp$-category will be called a \emph{spectral} category and an $\Sp$-functor will similarly be called a \emph{spectral} functor.
\end{itemize}

All adjunctions between symmetric monoidal categories considered in this paper are \textbf{strong symmetric} monoidal (meaning their left adjoint is strong symmetric monoidal), so we will simply call them \emph{monoidal adjunctions}.

Since we will frequently change the base of enrichment, our categories are often decorated by subscripts indicating which enrichment on that category is considered. For example, if $\Cin$ is a $\Top$-enriched category, then its underlying simplicial category, obtained by taking the singular complex of each hom-space, is denoted $\Cin_\Delta$, while its underlying ordinary category is denoted $\Cin_0$. When considering functors between enriched categories, we use the convention that the enrichment of the domain category indicates the enrichment of the functor. For example, if $\Cin$ and $\Din$ are topological categories, then $\Fun(\Cin,\Din)$ indicates the category of continuous functors, while $\Fun(\Cin_0, \Din)$ indicates the category of ordinary functors from $\Cin$ to $\Din$.

\subsection*{Acknowledgements}
The author wishes to thank Gregory Arone for many helpful discussions on the contents of this paper. The author is also grateful to the Max Planck Institute for Mathematics in Bonn for its hospitality, as part of this paper was written during his stay there.

%% file: Sections/Preliminaries.tex
\section{Preliminaries}\label{sec:Preliminaries}

We will first discuss some preliminaries on enriched categories, model categories, and bar constructions that will be used throughout the rest of the paper. Most of the definitions and facts discussed here are standard, although the notion of ``goodness'' discussed in \Cref{def:TopologicalGoodness,def:PointedGoodness,def:SpectralGoodness} is not. We advise the reader to only skim these preliminaries at a first reading and refer back to them as necessary while reading the rest of the paper.

\subsection{Enriched category theory}\label{ssec:EnrichedCats}

Throughout this section, let $(\calV, \otimes, \unit)$ be a closed symmetric monoidal category that is both complete and cocomplete. For the convenience of the reader, we will briefly discuss several concepts from enriched category theory that are used throughout this paper. We will assume that the reader is familiar with basic notions of enriched category theory, in particular $\calV$-categories, $\calV$-functors, $\calV$-natural transformations and $\calV$-adjunctions, the definition of $\calV$-(co)completeness, the construction of the $\calV$-category $\Fun(\Cin,\calE)$ of $\calV$-functors between two $\calV$-categories $\Cin$ and $\calE$, and the enriched Yoneda lemma.\footnote{Since we define (co)tensors and enriched (co)ends below, which are often used to define $\calV$-(co)completeness and to construct the hom-objects in $\Fun(\Cin,\calE)$, this assumption is somewhat circular. However, the purpose of this section is merely to recall the concepts from enriched category theory that play an important role throughout this paper and not to serve as a self-contained introduction.} For an introduction to these concepts, we refer the reader to \cite[\S 3]{Riehl2014Categorical} and \cite{Kelly1982BasicConcepts}.

We will first recall the definitions of (co)tensors, (co)ends, and enriched Kan extensions, and then discuss changing the base of enrichment along a symmetric monoidal adjunction.

\paragraph{(Co)tensors}

\begin{definition}
	Let $\Cin$ be a $\calV$-category and $v$ an object in $\calV$. We say that $\Cin$ \emph{admits tensors by $v$} if for any object $c$ of $\Cin$, there exists an object $v \otimes c$ in $\Cin$ together with isomorphisms
	\begin{equation*}\Cin(v \otimes c, d) \cong \calV(v, \Cin(c,d)) \end{equation*}
	that are $\calV$-natural in $d$. Dually, we say that $\Cin$ \emph{admits cotensors by $v$} if for any object $d$ of $\Cin$, there exists an object $d^v$ in $\Cin$ together with isomorphisms
	\begin{equation*}\Cin(c,d^v) \cong \calV(v, \Cin(c,d)) \end{equation*}
	that are $\calV$-natural in $c$.
\end{definition}

Throughout this paper, we will use that (co)tensors, when they exist, can be made functorial. To make this precise, let us define the product of two $\calV$-categories.

\begin{definition}\label{definition:TensorProductVCategories}
	Given $\calV$-categories $\Cin$ and $\Din$, define the $\calV$-category $\Cin \otimes \Din$ by
	\begin{enumerate}[(a)]
		\item defining the objects as $\Ob(\Cin \otimes \Din) := \Ob(\Cin) \times \Ob(\Din)$, and
		\item defining the hom-objects by $(\Cin \otimes \Din)((c,d),(c',d')) := \Cin(c,c') \otimes \Din(d,d')$.
	\end{enumerate}
	The identities and compositions are defined in the obvious way, making use of the symmetry of the monoidal structure on $\calV$.
\end{definition}

\begin{remark}
	This product can be promoted to a closed symmetric monoidal structure on the category $\calV\mhyphen\Cat$ of small $\calV$-categories. The internal hom for this symmetric monoidal structure is given by the $\calV$-category $\Fun(\Cin,\Din)$ of $\calV$-functors (see \cite[\S 2.3]{Kelly1982BasicConcepts} for details).
\end{remark}

\begin{construction}\label{construction:TensorsAreFunctorial}
	Let $\Cin$ be a $\calV$-category and suppose that $\calU$ is a class of objects in $\calV$ such that $\Cin$ admits tensors by every $u \in \calU$. Let $\calU$ abusively also denote the full sub-$\calV$-category of $\calV$ spanned by these objects. Then the tensors $u \otimes c$ can be assembled into a $\calV$-functor $\calU \otimes \Cin \to \Cin$ and the natural transformation $\Cin(u \otimes c, d) \cong \calV(u, \Cin(c,d))$ into a $\calV$-natural transformation of functors $\calU^\mathrm{op} \otimes \Cin^\mathrm{op} \otimes \Cin \to \calV$. This essentially follows by identifying $\calU \otimes \Cin$ with its image in $\Fun(\calU^\mathrm{op} \otimes \Cin^\mathrm{op}, \calV)$ under the Yoneda embedding; see \cite[\S 1.10]{Kelly1982BasicConcepts} for details. Of course, the dual statement holds for cotensors.
\end{construction}

\paragraph{(Co)ends}

In the following definition of enriched (co)ends, we will assume that $\calE$ is a $\calV$-(co)complete category. For a definition that works for any $\calV$-category $\calE$ and which is equivalent to the one given here when $\calE$ is $\calV$-(co)complete, we refer the reader to \cite[\S 3.10]{Kelly1982BasicConcepts}.

\begin{definition}
	Let $\Cin$ be a small $\calV$-category, $\calE$ a $\calV$-(co)complete $\calV$-category, and let a $\calV$-functor $F \colon \Cin^\mathrm{op} \otimes \Cin \to \calE$ be given. Then the \emph{coend} $\int^\Cin F(c,c)$ is defined as the coequalizer
	\begin{equation*} \int^\Cin F(c,c) := \mathrm{coeq}\left(\coprod_{c,c'} \Cin(c,c') \otimes F(c',c) \rightrightarrows \coprod_c F(c,c) \right)\end{equation*}
	in $\calE$, and the \emph{end} $\int_\Cin F(c,c)$ as the equalizer
	\begin{equation*} \int_\Cin F(c,c) := \mathrm{eq}\left(\prod_c F(c,c) \rightrightarrows \prod_{c,c'} F(c,c')^{\Cin(c,c')} \right).\end{equation*}
\end{definition}

\begin{example}
	Suppose $\calE$ is a $\s\Set_*$-cocomplete simplicial category and let $X_\bullet$ be a simplicial object in $\calE$. Then the \emph{geometric realization} of $X_\bullet$ in $\calE$ is the coend
	\begin{equation*}| X_\bullet | = \int^{\Delta^\mathrm{op}} \Delta^n \otimes X_n \end{equation*}
\end{example}

\begin{example}
	Given two $\calV$-functors of $\calV$-categories $F,G \colon \Cin \to \Din$ with $\Cin$ small, the \emph{natural transformations object} $\Nat(F,G)$ is defined as the end
	\begin{equation*}\Nat(F,G) = \int_\Cin \Din(F(c),G(c)). \end{equation*}
	These are the hom-objects of the $\calV$-category $\Fun(\Cin,\Din)$.
\end{example}

There exists the following ``Fubini theorem'' for enriched (co)ends (cf. (3.63) of \cite{Kelly1982BasicConcepts}).

\begin{lemma}\label{lemma:Fubini}
	Let $\Cin, \Din$ be small $\calV$-categories, $\calE$ a $\calV$-(co)complete $\calV$-category and $F \colon \Cin^\mathrm{op} \otimes \Cin \otimes \Din^\mathrm{op} \otimes \Din \to \calE$ a $\calV$-functor. Then there are canonical isomorphisms
	\begin{equation*}\int^\Cin \int^\Din F(c,c,d,d) \cong \int^{\Cin \otimes \Din} F(c,c,d,d) \cong \int^{\Din} \int^\Cin F(c,c,d,d) \end{equation*}
	and
	\begin{equation*}\int_\Cin \int_\Din F(c,c,d,d) \cong \int_{\Cin \otimes \Din} F(c,c,d,d) \cong \int_{\Din} \int_\Cin F(c,c,d,d). \end{equation*}
\end{lemma}

\paragraph{Kan extensions}

As in our definition of (co)ends, we will define enriched Kan extensions only for functors into a $\calV$-(co)complete $\calV$-category $\calE$, favouring the explicit formula that one can write down in this case over the more general definition given in \cite[\S 4.1]{Kelly1982BasicConcepts}.

\begin{definition}\label{def:EnrichedKanExtension}
	Let $\Cin$, $\Din$, and $\calE$ be $\calV$-categories, let $\alpha \colon \Cin \to \Din$ be a functor, and suppose that $\Cin$ is small and that $\calE$ is $\calV$-(co)complete. For a functor $F \colon \Cin \to \calE$, we define the \emph{(enriched) left Kan extension} $\Lan_\alpha F \colon \Din \to \calE$ of $F$ along $\alpha$ by
	\begin{equation*}\Lan_\alpha F(d) = \int^\Cin \Din(\alpha c, d) \otimes F(c) \end{equation*}
	and the \emph{(enriched) right Kan extension} $\Ran_\alpha F \colon \Din \to \calE$ of $F$ along $\alpha$ by
	\begin{equation*}\Ran_\alpha F(d) = \int_\Cin F(c)^{\Din(d,\alpha c)}. \end{equation*}
\end{definition}

By \cite[Thm.\ 4.50]{Kelly1982BasicConcepts}, Kan extension provides left and right adjoints to the functor $\Fun(\Din,\calE) \to \Fun(\Cin, \calE)$ given by precomposition with $\alpha$.

\begin{proposition}\label{prop:LeftAndRightKanExtensionAreAdjoints}
	If in the hypotheses of \Cref{def:EnrichedKanExtension} we furthermore assume that $\Din$ is small, then all (enriched) left and right Kan extensions assemble into functors $\Lan_\alpha, \Ran_\alpha \colon \Fun(\Cin,\calE) \to \Fun(\Din, \calE)$ that are left and right $\calV$-adjoint to the precomposition functor $\alpha^* \colon \Fun(\Din, \calE) \to \Fun(\Cin, \calE)$, respectively.
\end{proposition}

\paragraph{Change of base}

Recall that an adjunction $\calV \rightleftarrows \calW$ between symmetric monoidal categories is called \emph{strong symmetric monoidal} if the left adjoint is strong symmetric monoidal; note that in this case the right adjoint is automatically lax monoidal. All adjunctions between symmetric monoidal categories considered in this paper are \textbf{strong symmetric} monoidal, so we will simply call them \emph{monoidal adjunctions} from now on.

Fix a monoidal adjunction $L : \calV \rightleftarrows \calW : U$ between (co)complete closed symmetric monoidal categories. Then any $\calW$-category $\Din$ has an \emph{underlying} $\calV$-category $U\Din$ obtained by applying $U$ to each hom-object of $\Din$. Similarly, for a $\calV$-category $\Cin$, one can construct the free $\calW$-category $L\Cin$ on $\Cin$ by applying $L$ to each hom-object. These define an adjunction
\begin{equation}\label{eq:AdjunctionBaseChange}
	L : \calV\mhyphen\Cat \rightleftarrows \calW\mhyphen\Cat : U.
\end{equation}
When in the presence of a monoidal adjunction $L : \calV \rightleftarrows \calW : U$ and when given a $\calV$-category $\Cin$ and a $\calW$-category $\Din$, we will write $\Fun(\Cin, \Din)$ to denote the $\calW$-category $\Fun(L\Cin, \Din)$ (which agrees with $\Fun(\Cin, U\Din)$ as a $\calV$-category by \Cref{eq:AdjunctionBaseChange}).

If we are given a second $\calW$-category $\calE$, then the functor $\Fun(\Din, \calE) \to \Fun(U\Din, U\calE)$ sending a $\calW$-functor to its underlying $\calV$-functor agrees with the functor $\Fun(\Din, \calE) \to \Fun(LU\Din, \calE)$ given by precomposition with the counit $LU\Cin \to \Cin$. In particular, \Cref{prop:LeftAndRightKanExtensionAreAdjoints} gives the following result.

\begin{lemma}\label{lemma:BaseChangeGivenByKanExtension}
		Let $\Cin$ be a small $\calW$-category and $\calE$ a $\calW$-(co)complete $\calW$-category. Then the functor $\Fun(\Cin,\calE) \to \Fun(U\Cin, \calE)$ sending a $\calW$-functor to its underlying $\calV$-functor has both a left and right $\calW$-adjoint, which are given by left and right Kan extension along $LU\Cin \to \Cin$, respectively.
\end{lemma}

Finally, we wish to mention the following useful fact (cf.\ \cite[Thm.\ 3.7.11]{Riehl2014Categorical}).

\begin{lemma}\label{lemma:UnderlyingVCatComplete}
	Let $L : \calV \rightleftarrows \calW : U$ be a monoidal adjunction between (co)complete closed symmetric monoidal categories. Then for any $\calW$-(co)complete category $\calE$, the underlying $\calV$-category $U\calE$ is $\calV$-(co)complete. Moreover, tensors and cotensors by $v \in \Ob(\calV)$ in $U\calE$ are computed as tensors and cotensors by $Lv$ in $\calE$.
\end{lemma}

\subsection{Model categories}\label{ssec:ModelCats}

We will assume that the reader is familiar with the basic theory of model categories as explained in e.g.\ \cite[Ch.\ 1]{Hovey1999ModelCats} or \cite[Chs.\ 7-8]{Hirschhorn2003Model}.

\paragraph{Enriched model categories}
For the reader's convenience, we briefly recall a few facts about enriched model categories. We wish to mention that in all the symmetric monoidal model categories considered throughout this paper, the unit is cofibrant. Since the theory of enriched model categories simplifies if this is the case, we have decided to make this part of our definition.

\begin{definition}\label{definition:MonoidalModelCategory}
	A \emph{symmetric monoidal model category} $(\calV, \otimes, \unit)$ is a (co)complete closed symmetric monoidal category together with a model structure on $\calV$ such that
	\begin{enumerate}[(a)]
		\item\label{definition:MonoidalModelCategory:PushoutProduct} for any pair of cofibrations $v \cofarrow v'$ and $w \cofarrow w'$, the pushout-product map
		\begin{equation*}v \otimes w' \cup_{v \otimes w} v' \otimes w \to v' \otimes w' \end{equation*}
		is a cofibration which is trivial if one of $v \cofarrow v'$ and $w \cofarrow w'$ is, and
		\item the unit $\unit$ is cofibrant.
	\end{enumerate}
\end{definition}

\begin{definition}\label{definition:EnrichedModelCategory}
	Given a symmetric monoidal model category $\calV$, a \emph{$\calV$-model category} is a $\calV$-category  $\modN$ that is equipped with a model structure on its underlying category $\modN_0$ such that
	\begin{enumerate}[(a)]
		\item $\modN$ is $\calV$-(co)complete, and
		\item for any cofibration $M \cofarrow N$ in $\modN_0$ and any fibration $L \fibarrow K$ in $\modN_0$, the pullback-power map
		\begin{equation*}\modN(N,L) \to \modN(N,K) \times_{\modN(M,K)} \modN(M,L) \end{equation*}
		is a fibration in $\calV$, which is trivial if either $M \cofarrow N$ or $L \fibarrow K$ is.
	\end{enumerate}
\end{definition}

Given two symmetric monoidal model categories $\calV$ and $\calW$, we define a \emph{monoidal Quillen pair} as a (strong symmetric) monoidal adjunction $\calV \rightleftarrows \calW$ that is also a Quillen pair. Given such a monoidal Quillen pair $L : \calV \rightleftarrows \calW : U$ and a $\calW$-model category $\modN$, it follows from \Cref{lemma:UnderlyingVCatComplete} and the fact that $U$ preserves (trivial) fibrations that the underlying $\calV$-category $U\modN$ is a $\calV$-model category. The main question that this paper is concerned with is the following: given a $\calW$-category $\Cin$, which $\calV$-functors $U\Cin \to \modN$ are equivalent to a $\calW$-functor $\Cin \to \modN$?

\begin{definition}
	In a situation as above, we say that a $\calV$-functor $F \colon U\Cin \to \modN$ is \emph{equivalent to a $\calW$-functor} if there exists a $\calW$-functor $G \colon \Cin \to \modN$ and a zigzag of $\calV$-natural transformations between $F$ and $UG$ that are pointwise weak equivalences.
\end{definition}

\begin{remark}
	What we call a monoidal Quillen pair above should actually be called a strong symmetric monoidal Quillen pair. However, since all monoidal adjunctions that we consider in this paper are strong symmetric, we will drop these adjectives.
\end{remark}

\begin{example}
	The category $\s\Set$ of simplicial sets can be equipped with the Kan-Quillen model structure \cite[\S II.3]{Quillen1967HomotopicalAlgebra}. This is a symmetric monoidal model category under the Cartesian product. A $\s\Set$-model category will be called a \emph{simplicial model category}.
\end{example}

\begin{example}
	Let $\Top$ denote the category of compactly generated weak Hausdorff spaces. Then $\Top$, endowed with the Quillen model structure \cite[\S II.3]{Quillen1967HomotopicalAlgebra}, is a symmetric monoidal model category under the Cartesian product. A $\Top$-model category will be called a \emph{topological model category}. The singular complex and geometric realization form a monoidal Quillen equivalence
	\begin{equation*}|-| : \s\Set \rightleftarrows \Top : \Sing,\end{equation*}
	hence any topological model category has an underlying simplicial model category.
\end{example}

\begin{example}
	The categories $\s\Set_*$ and $\Top_*$ have model structures in which a map is a fibration, cofibration, or weak equivalence if its underlying map in $\s\Set$ or $\Top$ is. Both are symmetric monoidal model categories under the smash product, and the adjunctions
	\begin{equation*}(-)_+ : \s\Set \rightleftarrows \s\Set_* : U \quad \text{and} \quad (-)_+ : \Top \rightleftarrows \Top_* : U, \end{equation*}
	where $(-)_+$ adds a disjoint basepoint and $U$ forgets the basepoint, are monoidal Quillen adjunctions. Moreover, the singular complex and geometric realization again form a monoidal Quillen equivalence
	\begin{equation*}|-| : \s\Set_* \rightleftarrows \Top_* : \Sing.\end{equation*}
\end{example}

\paragraph{Orthogonal spectra}
In this paper, we will work with the category $\Sp$ of orthogonal spectra as our preferred symmetric monoidal model category of spectra. Recall that the category of orthogonal spectra is defined as the enriched diagram category $\Fun(\OIndex,\Top_*)$, where $\OIndex$ is the $\Top_*$-category whose objects are the natural numbers and whose hom-objects are given by $\OIndex(m,n) = O(n)_+ \wedge_{O(n-m)} S^{n-m}$. The category $\OIndex$ is symmetric monoidal, hence $\Sp$ is closed symmetric monoidal under the Day convolution product $\otimes$. A category enriched in $\Sp$ will be called a \emph{spectral} category.

We consider $\Sp$ as a model category by endowing it with the stable model structure of \cite{MMSS2001DiagramSpectra}. This is a left Bousfield localization of the projective model structure compatible with the Day convolution product. In particular, the stable model structure on $\Sp$ is a symmetric monoidal model category in the sense of \Cref{definition:MonoidalModelCategory}. The category of orthogonal spectra $\Sp$ comes with a strong symmetric monoidal left Quillen functor $\Sigma^\infty$ from $\Top_*$, hence we can speak of the underlying $\Top_*$-category of a spectral category. An $\Sp$-model category in the sense of \Cref{definition:EnrichedModelCategory} will be called a \emph{spectral} model category.

Recall that a model category is called \emph{stable} if it has a zero object and the suspension functor induces a self-equivalence on the homotopy category.

\begin{lemma}\label{lemma:SpectrallyEnrichedIsStable}
	Any spectral model category is stable.
\end{lemma}

\begin{proof}
	A spectral model category $\modN$ has a zero object since it is (co)tensored over $\Sp$ and $\Sp$ has a zero object. To see that it is stable, note that we can compute suspensions as (derived) tensors with $\SSS^1$. If we let $\SSS^{-1}$ denote the image of $1$ under the Yoneda embedding $\OIndex^\mathrm{op} \hookrightarrow \Sp$, then we have a stable equivalence $\SSS^{-1} \otimes \SSS^1 \cong \SSS^1 \otimes \SSS^{-1} \to \SSS^0$. In particular, we have natural equivalences $\SSS^1 \otimes \SSS^{-1} \otimes N \to N$ and $\SSS^{-1} \otimes \SSS^1 \otimes N \to N$ for any cofibrant $N$ in $\modN$. This shows that on the level of homotopy categories, (derived) tensors with $\SSS^{-1}$ are inverse to the suspension functor.
\end{proof}

The following fact, which follows from \cite[Prop.\ 6.3]{SchwedeShipley2003Equivalences}\footnote{Strictly speaking, one of the hypotheses of \cite[Prop.\ 6.3]{SchwedeShipley2003Equivalences} is that all objects are small with respect to the whole category, which is not the case in $\Sp$. However, all objects are small with respect to levelwise closed inclusions, which is sufficient for the proof to go through.}, will be used several times in \Cref{sec:LinearFunctors}.

\begin{lemma}\label{lemma:CofibrantReplacementOfCategories}
	Let $\Cin$ be a small spectral category. Then there exists a small spectral category $\wh \Cin$ with the same set of objects and whose hom-objects are cofibrant, together with a map $\wh \Cin \trivfibarrow \Cin$ that is the identity on objects and a trivial fibration on all hom-objects.
\end{lemma}

\paragraph{Left Bousfield localization}
While left Bousfield localizations are not necessary to prove our main results, we will make use of them throughout this paper to upgrade results such as \Cref{TheoremA} and \Cref{TheoremC} to actual Quillen equivalences of model categories. We assume that the reader is familiar with the theory of left Bousfield localizations; we refer to \cite[\S10]{Lawson2020Bousfield} for an introduction and to \cite{Hirschhorn2003Model} for more details. Let us remark that throughout this paper, a left Bousfield localization of a model category $\modN$ is simply defined as a new model structure on $\modN$ with the same class of cofibrations but a larger class of weak equivalences. We will frequently use the following result.

\begin{lemma}\label{lemma:BousfieldLocalizationLemma}
	Let $L : \modM \rightleftarrows \modN : U$ be a Quillen pair such that
	\begin{enumerate}[(a)]
		\item\label{item1:BousfieldLocalizationLemma} $U$ creates weak equivalences, and
		\item\label{item2:BousfieldLocalizationLemma} for any cofibrant object $X$ in $\modM$ that is weakly equivalent to an object in the image of $U$, the unit $X \to ULX$ is a weak equivalence.
	\end{enumerate}
	Moreover, assume that there exists a left Bousfield localization $\modM'$ of $\modM$ in which the fibrant objects are exactly the fibrant objects of $\modM$ that are weakly equivalent to an object in the image of $U$. Then \begin{equation*}L : \modM' \rightleftarrows \modN : U\end{equation*} is a Quillen equivalence.
\end{lemma}

\begin{proof}
	We first show that $L : \modM' \rightleftarrows \modN : U$ is a Quillen pair. By Corollary A.2 and Remark A.3\footnote{See \cite[Prop.\ 3.7]{nlab:bousfield_localization_of_model_categoriesRev59} for a detailed proof of the fact that a map between fibrant objects in $\modM'$ is a fibration in $\modM'$ if and only if it is in $\modM$.} of \cite{Dugger2001Replacing}, it suffices to show that $U$ preserves fibrant objects, which is true by definition. To see that this is a Quillen equivalence, it suffices to show that the derived adjunction gives an equivalence of homotopy categories. Note that $\Ho(\modM')$ is the full subcategory of $\Ho(\modM)$ on objects that are in the essential image of $U$ (cf.\ \cite[Prop.\ 10.2]{Lawson2020Bousfield}). This shows that the right derived functor $\Ho(\modN) \to \Ho(\modM')$ is essentially surjective, while item \ref{item1:BousfieldLocalizationLemma} above ensures that this functor is conservative. Moreover, item \ref{item2:BousfieldLocalizationLemma} ensures that the unit of the derived adjunction is an isomorphism, hence $\Ho(\modN) \to \Ho(\modM')$ is an equivalence.
\end{proof}

\paragraph{Goodness}
We will now discuss a notion for (pointed) topological and spectrally enriched model categories $\modN$ that we call ``goodness'', and which can be seen as ensuring that the enriched functor categories $\Fun(\Cin,\modN)$ are well-behaved even if the hom-objects of $\Cin$ are not cofibrant. While the conditions of goodness may appear to be strong, we will show that they are satisfied in many natural situations. We warn the reader that the definitions given here are not standard and that some proofs are quite technical.

Recall that, given a $\calV$-model category $\modN$ and a small $\calV$-category $\Cin$, the projective model structure on $\Fun(\Cin,\modN)$ is defined (if it exists) as the model structure in which the fibrations and the weak equivalences are the pointwise ones. We will call a functor $F \colon \Cin \to \Din$ of $\calV$-categories that is bijective on objects a \emph{weak equivalence} if $\Cin(c,d) \to \Din(F(c),F(d))$ is a weak equivalence in $\calV$ for all $c,d \in \Ob(\Cin)$.

\begin{definition}\label{def:TopologicalGoodness}
    A topological model category $\modN$ is called \emph{good} if
    \begin{enumerate}[(a)]
        \item\label{item:TopologicalGoodness1} for any small topological category $\Cin$, the projective model structure on $\Fun(\Cin,\modN)$ exists, and
        \item\label{item:TopologicalGoodness2} for any weak equivalence $\alpha \colon \Cin \to \Din$ of small topological categories that is bijective on objects, the precomposition functor $\alpha^* \colon \Fun(\Din,\modN) \to \Fun(\Cin,\modN)$ is a right Quillen equivalence.
    \end{enumerate}
\end{definition}

Let us state some direct consequences of these conditions.

\begin{proposition}\label{proposition:FunctorTensorProductInGoodTopModelCat}
	Let $\modN$ be a good topological model category, $\Cin$ a small topological category, $F \colon \Cin \to \modN$ a projectively cofibrant functor and $G \wearrow G'$ a pointwise weak equivalence between functors $\Cin^\mathrm{op} \to \Top$. Then $\int^\Cin G \otimes F \to \int^\Cin G' \otimes F$ is a weak equivalence in $\modN$. Similarly, if $G$ is any functor $\Cin^\mathrm{op} \to \Top$ and $F \wearrow F'$ is a pointwise weak equivalence between projectively cofibrant functors $\Cin \to \modN$, then $\int^\Cin G \otimes F \to \int^\Cin G \otimes F'$ is a weak equivalence.
\end{proposition}

\begin{proof}
	Given a functor $G \colon \Cin^\mathrm{op} \to \Top$, one can form the $\Top$-category $\Cin \triangleright G$ with $\Ob(\Cin \triangleright G) = \Ob(\Cin) \coprod \{*\}$ and with hom-objects given by
	\begin{equation*}
		\begin{cases}
			(\Cin \triangleright G)(c,d) = \Cin(c,d) & \quad \text{when}\quad c,d \in \Ob(\Cin) \\
			(\Cin \triangleright G)(c,*) = G(c) & \quad \text{when}\quad c \in \Ob(\Cin)\\
			(\Cin \triangleright G)(*,d) = \varnothing & \quad \text{when}\quad d \in \Ob(\Cin) \\
			(\Cin \triangleright G)(*,*) = *
		\end{cases}
	\end{equation*}
	The composition is given by composition in $\Cin$ and the action of $\Cin$ on $G$. Let $i_G \colon \Cin \hookrightarrow \Cin \triangleright G$ denote the inclusion. It follows from the definition of left Kan extension that
	\begin{equation*}
		\int^\Cin G \otimes F \cong (\Lan_{i_G} F)(*).
	\end{equation*}
	Since $\Lan_{i_G}$ is left Quillen, it preserves pointwise equivalences $F \wearrow F'$ between projectively cofibrant functors. Evaluating at $*$ yields the desired equivalence $\int^\Cin G \otimes F \wearrow \int^\Cin G \otimes F'$.
	
	On the other hand, if we are given an equivalence $G \wearrow G'$, then we obtain an equivalence $\alpha \colon \Cin \triangleright G \wearrow \Cin \triangleright G'$. The goodness of $\modN$ now implies that
	\begin{equation*}
		\int^\Cin G \otimes F \cong (\Lan_{i_G} F)(*) \wearrow (\Lan_{\alpha} \Lan_{i_G} F)(*) \cong (\Lan_{i_{G'}} F)(*) \cong \int^\Cin G' \otimes F
	\end{equation*}
	is a weak equivalence for any projectively cofibrant $F$.
\end{proof}

The special case $\Cin = *$ yields the following result on tensors in $\modN$.

\begin{proposition}\label{proposition:TensoringInGoodTopModelCat}
    Let $\modN$ be a good topological model category, $T$ a topological space, and $N$ a cofibrant object in $\modN$. Then $T \otimes - \colon \modN \to \modN$ preserves weak equivalences between cofibrant objects, while $- \otimes N \colon \Top \to \modN$ preserves all weak equivalences.
\end{proposition}

Goodness also implies the existence of projective model structures on simplicial diagram categories.

\begin{proposition}\label{proposition:GoodTopModelCatVsSimplicial}
    Let $\modN$ be a good topological model category. Then for any small simplicial category $\Cin$, the projective model structure on $\Fun(\Cin,\modN)$ exists. Moreover, for any topological category $\Din$, the forgetful functor $\Fun(\Din,\modN) \to \Fun(\Sing(\Din),\modN)$ is a right Quillen equivalence.
\end{proposition}

\begin{proof}
For the first part, note that $\Fun(\Cin,\modN) = \Fun(|\Cin|,\modN)$. For the second part, recall from \Cref{lemma:BaseChangeGivenByKanExtension} that the forgetful functor can be identified with the precomposition functor $\varepsilon^* \colon \Fun(\Din,\modN) \to \Fun(|\Sing(\Din)|,\modN)$, where $\varepsilon \colon |\Sing(\Din)| \wearrow \Din$ is the counit of the adjunction $|-| \dashv \Sing$. The result now follows from part \ref{item:TopologicalGoodness2} of \Cref{def:TopologicalGoodness}.
\end{proof}

As promised above, we will now show that many topological model categories appearing in nature are good.

\begin{proposition}\label{proposition:GoodTopologicalInNature}
    \begin{enumerate}[(i)]
        \item\label{item1:GoodTopologicalInNature} $\Top$ is a good topological model category.
        \item\label{item2:GoodTopologicalInNature} For any good topological model category $\modN$ and any small topological category $\Cin$, the projective model structure on $\Fun(\Cin,\modN)$ is good.
        \item\label{item3:GoodTopologicalInNature} For any small topological category $\Cin$ and any set of maps $S$ in $\Fun(\Cin,\Top)$, the left Bousfield localization $L_S\Fun(\Cin,\Top)$ is good.
    \end{enumerate}
\end{proposition}

\begin{proof}
Item \ref{item1:GoodTopologicalInNature} is well known. Item \ref{item2:GoodTopologicalInNature} follows since for any small topological category $\Din$, one has the equivalence of categories $\Fun(\Din, \Fun(\Cin,\modN)) \simeq \Fun(\Din \times \Cin, \modN)$.

For item \ref{item3:GoodTopologicalInNature}, we leave it to the reader to verify that the projective model structure on $\Fun(\Din, L_S\Fun(\Cin,\Top))$ exists. To verify item \ref{item:TopologicalGoodness2} of \Cref{def:TopologicalGoodness}, let $\Din \to \Din'$ be a weak equivalence that is bijective on objects. Then \begin{equation*}\Fun(\Din',\Fun(\Cin,\Top)) \to \Fun(\Din,\Fun(\Cin,\Top))\end{equation*}
is a right Quillen equivalence. Using this, one easily verifies the conditions of \cite[Lem.\ A.2.(iii)]{MMSS2001DiagramSpectra}, proving that \begin{equation*}\Fun(\Din',L_S\Fun(\Cin,\Top)) \to \Fun(\Din,L_S\Fun(\Cin,\Top))\end{equation*}
is also a right Quillen equivalence.
\end{proof}

We now define goodness for pointed topological categories. Recall that a basepoint $x_0$ of a topological space $X$ is \emph{nondegenerate} if $(X,x_0)$ is an NDR-pair in the sense of \cite[\S 6.4]{May1999Concise}. Since constructions in $\Top_*$ are generally only well-behaved with respect to nondegenerately based spaces, we will always assume that the hom-objects of our indexing categories are of this kind.

\begin{definition} \label{def:PointedGoodness}
    A pointed topological model category $\modN$ is called \emph{good} if
    \begin{enumerate}[(a)]
        \item\label{item:PointedGoodness1} for any small  $\Top_*$-category $\Cin$ whose hom-spaces are nondegenerately based, the projective model structure on $\Fun(\Cin,\modN)$ exists, and
        \item\label{item:PointedGoodness2} for any two small $\Top_*$-categories $\Cin$ and $\Din$ whose hom-objects are nondegenerately based and any weak equivalence $\alpha \colon \Cin \to \Din$ that is bijective on objects, the precomposition functor $\alpha^* \colon \Fun(\Din,\modN) \to \Fun(\Cin,\modN)$ is a right Quillen equivalence.
    \end{enumerate}
\end{definition}

\begin{remark}\label{remark:PointedGoodIsUnderlyingGood}
If $\modN$ is a good pointed topological model category, then its underlying topological model category is also good. Namely, if $\Cin$ is a small topological category, then $\Fun(\Cin,\modN) = \Fun(\Cin_+,\modN)$, where $\Cin_+$ is obtained by adding disjoint basepoints to the hom-objects of $\Cin$. The result follows since disjoint basepoints are always nondegenerate.
\end{remark}

One can formulate analogues of \Cref{proposition:FunctorTensorProductInGoodTopModelCat,proposition:TensoringInGoodTopModelCat,proposition:GoodTopModelCatVsSimplicial,proposition:GoodTopologicalInNature} in the setting of pointed topological model categories and prove them in the same way. For future reference, we state the analogue of \Cref{proposition:FunctorTensorProductInGoodTopModelCat}.

\begin{proposition}\label{proposition:FunctorTensorProductInGoodPointedModelCat}
	Let $\modN$ be a good $\Top_*$-model category, $\Cin$ a small $\Top_*$-category whose hom-spaces have nondegenerate basepoints, $F \colon \Cin \to \modN$ a projectively cofibrant functor and $G \wearrow G'$ a pointwise weak equivalence between functors $\Cin^\mathrm{op} \to \Top_*$ that are pointwise nondegenerately based. Then $\int^\Cin G \otimes F \to \int^\Cin G' \otimes F$ is a weak equivalence in $\modN$. Similarly, if $G$ is a functor $\Cin^\mathrm{op} \to \Top_*$ that is pointwise nondegenerately based and $F \wearrow F'$ is a pointwise weak equivalence between projectively cofibrant functors $\Cin \to \modN$, then $\int^\Cin G \otimes F \to \int^\Cin G \otimes F'$ is a weak equivalence.
\end{proposition}

Finally, let us define goodness for spectral model categories.

\begin{definition}\label{def:SpectralGoodness}
    A spectral model category $\modN$ is called \emph{good} if
    \begin{enumerate}[(a)]
        \item\label{item1:SpectralGoodness} for any small spectral category $\Cin$, the projective model structure on $\Fun(\Cin,\modN)$ exists, and
        \item\label{item2:SpectralGoodness} for any weak equivalence $\alpha \colon \Cin \to \Din$ of small spectral categories that is bijective on objects, the precomposition functor $\alpha^* \colon \Fun(\Din,\modN) \to \Fun(\Cin,\modN)$ is a right Quillen equivalence.
    \end{enumerate}
\end{definition}

\begin{remark}
	If $\modN$ is a good spectral model category, then its underlying $\Top_*$-model category is also good.
\end{remark}

One can mimic the proofs of \Cref{proposition:FunctorTensorProductInGoodTopModelCat,proposition:TensoringInGoodTopModelCat} to obtain analogous statements for good spectral model categories, which we state for future reference.

\begin{proposition}\label{proposition:FunctorTensorProductInGoodSpectralModelCat}
	Let $\modN$ be a good spectral model category, $\Cin$ a small spectral category, $F \colon \Cin \to \modN$ a projectively cofibrant functor and $G \wearrow G'$ a pointwise stable equivalence between functors $\Cin^\mathrm{op} \to \Sp$. Then $\int^\Cin G \otimes F \to \int^\Cin G' \otimes F$ is a weak equivalence in $\modN$. Analogously, if $G$ is any functor $\Cin^\mathrm{op} \to \Sp$ and $F \wearrow F'$ is a pointwise weak equivalence between projectively cofibrant functors $\Cin \to \modN$, then $\int^\Cin G \otimes F \to \int^\Cin G \otimes F'$ is a weak equivalence.
\end{proposition}

\begin{proposition}\label{proposition:TensoringInGoodSpectralModelCat}
	Let $\modN$ be a good spectral model category, $T$ an object in $\Sp$, and $N$ a cofibrant object in $\modN$. Then $T \otimes - \colon \modN \to \modN$ preserves weak equivalences between cofibrant objects, while $- \otimes N \colon \Sp \to \modN$ preserves all weak equivalences.
\end{proposition}

One can also prove an analogue of \Cref{proposition:GoodTopologicalInNature}, showing that many common examples of spectral model categories are good.

\begin{proposition}\label{proposition:GoodSpectralInNature}
	\begin{enumerate}[(i)]
		\item\label{item1:GoodSpectralInNature} $\Sp$ is a good spectral model category.
		\item\label{item2:GoodSpectralInNature} For any good spectral model category $\modN$ and any small $\Sp$-category $\Cin$ whose hom-objects are cofibrant, the projective model structure on $\Fun(\Cin,\modN)$ is good.
		\item\label{item3:GoodSpectralInNature} Let $\Cin$ be a small spectral category whose hom-objects are cofibrant and let $S$ be a set of maps in $\Fun(\Cin,\Sp)$. Then the enriched left Bousfield localization $L_S \Fun(\Cin,\Sp)$ of the projective model structure at $S$ is a good spectral model category.
	\end{enumerate}
\end{proposition}

\begin{proof}
	For item \ref{item1:GoodSpectralInNature}, note that properties \labelcref{item1:SpectralGoodness,item2:SpectralGoodness} of \Cref{def:SpectralGoodness} follow from \cite[Thm.\ 7.2]{SchwedeShipley2003Equivalences}.
	
	For item \ref{item2:GoodSpectralInNature}, property \ref{item1:SpectralGoodness} of \Cref{def:SpectralGoodness} follows as in \Cref{proposition:GoodTopologicalInNature}. For property \ref{item2:SpectralGoodness} of \Cref{def:SpectralGoodness}, note that since $\Cin$ has cofibrant hom-objects, one has by \cite[Prop.\ 12.3]{MMSS2001DiagramSpectra} that $\Cin \otimes \Din \wearrow \Cin \otimes \Din'$ for any weak equivalence $\Din \wearrow \Din'$ of spectral categories.
	
	Finally, note for item \ref{item3:GoodSpectralInNature} that the spectrally enriched left Bousfield localization at a set $S$ is the same as the unenriched left Bousfield localization at the set $S'$ that consists of all shifts of maps in $S$. The result now follows in the same way as in the proof of \Cref{proposition:GoodTopologicalInNature}.
\end{proof}

\subsection{Bar constructions and derived Kan extensions}\label{ssec:BarConstructions}

Given a $\calV$-model category $\modN$ and a $\calV$-functor $\alpha \colon \Cin \to \Din$ between small $\calV$-categories, derived (enriched) left Kan extension is defined as the left derived functor of the ordinary (enriched) left Kan extension functor
$\Lan_\alpha \colon \Fun(\Cin,\modN) \to \Fun(\Din, \modN)$.
One can obtain such a left derived functor by endowing these diagram categories with the projective model structures and noting that $\Lan_\alpha$ is left Quillen with respect to these. However, the projective model structures might not always exist, and the resulting description of the derived left Kan extension involves a projective cofibrant replacement functor, which is generally very inexplicit. In certain cases, one can use the bar construction to give an alternative construction of the derived (enriched) left Kan extension functor which circumvents the first problem and moreover produces an explicit formula. Throughout the rest of this section, $\modN$ will be a $\calV$-model category and $\alpha \colon \Cin \to \Din$ will be a $\calV$-functor between small $\calV$-categories. Moreover, we will assume that the symmetric monoidal model category $\calV$ comes equipped with a monoidal Quillen pair
\begin{equation*}L : \s\Set \rightleftarrows \calV : U. \end{equation*}
In particular, any $\calV$-model category has an underlying simplicial model category and one can speak of geometric realizations of simplicial objects in a $\calV$-model category.

We wish to point out that all definitions and results of this section can be dualized to obtain statements about the cobar construction and derived enriched right Kan extensions.

\begin{definition}
	Let $\calE$ be a $\calV$-cocomplete $\calV$-category, $\Cin$ a small $\calV$-category and let $F \colon \Cin \to \calE$, $G \colon \Cin^\mathrm{op} \to \calV$ be $\calV$-functors. The \emph{bar construction} $B(G,\Cin,F)$ is defined as the geometric realization of the simplicial object $B_\bullet(G,\Cin,F)$ in $\calE$ given by
	\begin{equation*}B_n(G,\Cin,F) = \coprod_{c_0,\ldots,c_n \in \Ob(\Cin)} G(c_n) \otimes \Cin(c_{n-1},c_n) \otimes \cdots \otimes \Cin(c_0,c_1) \otimes F(c_0), \end{equation*}
	with the face and degeneracy maps as in \cite[Def.\ 9.1.1]{Riehl2014Categorical}.
\end{definition}

Let us record a few lemmas on the bar construction.

\begin{lemma}\label{lemma:BarCommutesWithColimits}
	$B(-,\Cin,F) \colon \Fun(\Cin^\mathrm{op},\calV) \to \calE$ preserves colimits and tensors by objects of $\calV$.
\end{lemma}

\begin{proof}
	This follows since geometric realization commutes with colimits and tensors by objects of $\calV$.
\end{proof}

We will generally be interested in the case where we are given a $\calV$-functor $\alpha \colon \Cin \to \Din$ and $G$ is of the form $c \mapsto \Din(\alpha c, d)$, where $d$ is an object of $\Din$. By letting $d$ vary, we obtain a functor $\Din \to \calE$ given by
\begin{equation*}d \mapsto B(\Din(\alpha -, d),\Cin,F). \end{equation*}
For brevity, we will denote this functor by $B(\Din,\Cin,F)$; the map $\alpha$ should always be clear from the context. When we write $B(\Cin,\Cin,F)$, we always mean the case that $\alpha = \id_\Cin$. In this case, there is a canonical map $B(\Cin,\Cin,F) \to F$ coming from the augmentation
\begin{equation*}B_0(\Cin,\Cin,F)(c) = \coprod_d \Cin(d,c) \otimes F(d) \xrightarrow{\ac} F(c). \end{equation*}

\begin{lemma}\label{lemma:BarExtraDegeneracy}
	Let $F \colon \Cin \to \modN$ be a $\calV$-functor. Then the map $B(\Cin(-,c),\Cin,F) \to F(c)$, and hence $B(\Cin,\Cin,F) \to F$, is a natural weak equivalence.
\end{lemma}

\begin{proof}
	This can be proved by showing that $B(\Cin(-,c),\Cin,F)$ has an ``extra degeneracy'', cf.\ \cite[Ex.\ 4.5.7]{Riehl2014Categorical}.
\end{proof}

\begin{lemma}\label{lemma:BarKanExtension}
	Let $F \colon \Cin \to \modN$ be a $\calV$-functor. Then $\Lan_\alpha B(\Cin,\Cin,F) \cong B(\Din,\Cin,F)$.
\end{lemma}

\begin{proof}
	This follows from \Cref{lemma:Fubini} and the fact that $\Lan_\alpha(\Cin(c,-) \otimes E) \cong \Din(\alpha c, -) \otimes E$ for any $E$ in $\modN$.
\end{proof}

Since $B(\Cin,\Cin,F) \to F$ is a natural weak equivalence, we see that
\begin{equation*}
B(\Cin,\Cin,-) \colon \Fun(\Cin,\modN) \to \Fun(\Cin,\modN)
\end{equation*}
is a left deformation in the sense of \cite[\S 3]{Shulman2006HomotopyLimits}. If one combines this functor with a cofibrant replacement functor $Q \colon \modN \to \modN$, then one can often show that $\Lan_\alpha$ is homotopical on the image of $B(\Cin,\Cin, Q-)$ (cf.\ \cite[Thm.\ 13.7]{Shulman2006HomotopyLimits}) and hence that $\Lan_\alpha B(\Cin,\Cin,Q-) = B(\Din,\Cin,Q-)$ is a left derived functor of $\Lan_\alpha$. It is important to note that for this to work, $Q$ needs to be a $\calV$-functor and not just an ordinary functor.

\begin{lemma}\label{lemma:BarIsDerivedLeftKanExtension}
    Suppose that the hom-objects of $\Cin$ and $\Din$ are cofibrant and that $\unit \to \Cin(c,c)$ is a cofibration for any $c$ in $\Cin$. Then, for any $\calV$-model category $\modN$ with a $\calV$-enriched cofibrant replacement functor $Q$, the bar construction $B(\Din,\Cin,Q -)$ is a left derived functor of $\Lan_\alpha \colon \Fun(\Cin,\modN) \to \Fun(\Din,\modN)$.
\end{lemma}

\begin{proof}
This is shown in \cite[Thm.\ 13.7]{Shulman2006HomotopyLimits}.
\end{proof}

The following lemma is a variation on \cite[Prop.\ 23.6]{Shulman2006HomotopyLimits}.

\begin{lemma}\label{lemma:SimplicialBarIsReedyCofibrant}
	Suppose that all hom-objects of $\Cin$ are cofibrant and that $\unit \to \Cin(c,c)$ is a cofibration in $\calV$ for every $c$ in $\Cin$. Then for any pointwise cofibration $F \cofarrow F'$ in $\Fun(\Cin^\mathrm{op},\calV)$ and any pointwise cofibrant $G$ in $\Fun(\Cin,\modN)$, the map
	\begin{equation*}B_\bullet (F,\Cin,G) \to B_\bullet(F',\Cin,G)\end{equation*}
	is a Reedy cofibration in $\modN^{\Delta^\mathrm{op}}$.
\end{lemma}

\begin{proof}
	In order to prove this result, we give a slightly different presentation of the bar construction. Given a set $S$, we endow $\calV^S$, $\calV^{S \times S}$, and $\calN^S$ with the model structure in which the (co)fibrations and weak equivalences are defined pointwise. We have a product $\odot \colon \calV^{S \times S} \otimes \calV^{S \times S} \to \calV^{S \times S}$ given by
	\begin{equation*}(X \odot Y)(s,t) = \coprod_{u} X(u,s) \otimes Y(t,u),\end{equation*}
	endowing $\calV^{S \times S}$ with a (non-symmetric) monoidal structure. The monoids for this are exactly $\calV$-enriched categories whose set of objects is $S$. Given such a $\calV$-category $\Din$, we will write $\Din^{\odot}$ for the simplicial object in $\calV^{S \times S}$ given by
	\begin{equation*}(\Din^{\odot})_n := \Din^{\odot n} = \underbrace{\Din \odot \cdots \odot \Din}_{n \text{-times}}. \end{equation*}
	Observe that there are similar products $\calV^{S \times S} \otimes \modN^S \to \modN^S$, $\calV^S \otimes \calV^{S \times S} \to \calV^S$ and $\calV^S \otimes \modN^S \to \modN$, which we will also denote by $\odot$. It follows immediately from the definition that all these products satisfy the pushout-product axiom (cf.\ part \ref{definition:MonoidalModelCategory:PushoutProduct} of \Cref{definition:MonoidalModelCategory}). Taking $S = \Ob(\Cin)$, the bar construction can be rewritten as
	\begin{equation*}B_n(F,\Cin,G) = F \odot \Cin^{\odot n} \odot G.\end{equation*}
	Under this identification, the $n$-th latching map
	\begin{equation*}B_n(F,\Cin,G) \cup_{L_n B_\bullet (F,\Cin,G) } L_nB_\bullet(F',\Cin,G) \to B_n(F',\Cin,G) \end{equation*}
	is identified with the map
	\begin{equation*}F \odot \Cin^{\odot n} \odot G \cup_{F \odot L_n \Cin^{\odot} \odot G} F' \odot L_n \Cin^{\odot} \odot G \to F' \odot \Cin^{\odot n} \odot G,\end{equation*}
	which is the pushout-product of $F \to F'$ and $L_n \Cin^{\odot} \odot G \to \Cin^{\odot n} \odot G$. Since $G$ is pointwise cofibrant and $F \to F'$ is a pointwise cofibration, it follows that this map is a cofibration if $L_n \Cin^{\odot} \to \Cin^{\odot n}$ is. The latter follows exactly as in the proof of \cite[Prop.\ 23.6]{Shulman2006HomotopyLimits}.
\end{proof}

Note that for any functor $\alpha \colon \Cin \to \Din$, the precomposition functor $\alpha^* \colon \Fun(\Din,\modN) \to \Fun(\Cin,\modN)$ preserves (pointwise) weak equivalences, hence it is already right derived. In particular, if $B(\Din,\Cin,Q -)$ is a left derived functor of $\Lan_\alpha$, then
\begin{equation*} B(\Din, \Cin, Q -) : \Ho(\Fun(\Cin,\modN)) \rightleftarrows \Ho(\Fun(\Din,\modN)) : \alpha^*\end{equation*}
is an adjunction by \cite[Thm.\ 2.2.11]{Riehl2014Categorical}. A diagram chase shows that the unit and counit of this adjunction are given by
\begin{equation*} F \xleftarrow[\;\;\nu\;\;]{\sim} B(\Cin,\Cin,QF) \xrightarrow{\eta_{B(\Cin,\Cin,QF)}} \alpha^* B(\Din,\Cin,QF) \end{equation*}
and
\begin{equation*} B(\Din,\Cin,Q\alpha^*G) \xrightarrow{\Lan_\alpha(\nu)} \Lan_\alpha \alpha^* G \xrightarrow{\varepsilon_G} G ,\end{equation*}
respectively, where $\nu \colon B(\Cin,\Cin,QF) \wearrow F$ denotes the composition of the equivalence $B(\Cin,\Cin,QF) \wearrow QF$ of \Cref{lemma:BarExtraDegeneracy} and $QF \wearrow F$.

\begin{definition}\label{definition:DerivedUnitCounitBar}
    We call
    \begin{equation*}
    	\eta_{B(\Cin,\Cin,QF)} \colon B(\Cin,\Cin,QF) \to \alpha^* B(\Din,\Cin,QF) \quad \text{and} \quad \varepsilon_G \circ \Lan_\alpha(\nu) \colon B(\Din,\Cin,Q \alpha^* G) \to G
    \end{equation*}
    the \emph{derived unit} and \emph{derived counit}, respectively.
\end{definition}

%% file: Sections/HomotopyFunctors.tex
\section{Homotopy functors}\label{sec:HomotopyFunctors}

In this part, we will describe conditions under which a functor can be replaced by one that respects an enrichment in simplicial sets or topological spaces; in particular, we will prove \Cref{TheoremA}. We start by describing such conditions in the simplicial case in \Cref{ssec:SimplicialFunctors}, after which we will consider the topological case in \Cref{ssec:ContinuousFunctors}. In \Cref{ssec:DwyerKan}, we show how to obtain Proposition 1.3.4.7 of \cite{Lurie2017HigherAlgebra} as a corollary of the simplicial case.

\begin{remark}
	There are also versions of the main results of \Cref{ssec:SimplicialFunctors,ssec:ContinuousFunctors} for (symmetric) functors of several variables. We will not explicitly prove them here, but refer the reader to \Cref{ssec:MultireducedFunctors} and \Cref{remark:HomotopyMultiFunctors} instead.
\end{remark}

\subsection{Simplicial functors}\label{ssec:SimplicialFunctors}

Throughout this section $\Cin$ denotes a small simplicial category and $\modN$ a simplicial model category. The underlying category of $\Cin$ will be denoted by $\Cin_0$ whenever we need to distinguish it from $\Cin$.

\begin{definition}
\begin{enumerate}[(i)]
\item Two maps $f,g \colon C \to D$ in $\Cin$ are \emph{strictly homotopic} if there exists a $1$-simplex $H \in \Cin(C,D)_1$ such that either $d_0H = f$ and $d_1H = g$, or $d_1H = f$ and $d_0H=g$.
\item Two maps $f,g \colon C \to D$ in $\Cin$ are \emph{(simplicially) homotopic} if there exists a sequence of maps $f_0,\ldots,f_n \colon C \to D$ with $f_0 = f$ and $f_n = g$ such that for every $0 \leq i < n$, the map $f_i$ is strictly homotopic to $f_{i+1}$.
\item A map in $\Cin$ is a \emph{(simplicial) homotopy equivalence} if, up to homotopy, it has both a left and a right inverse.
\item A functor $\Cin_0 \to \modN$ is called a \emph{homotopy functor} if it sends simplicial homotopy equivalences to weak equivalences.
\end{enumerate}
\end{definition}

Note that any simplicial functor is automatically a homotopy functor since simplicial functors preserve (strict) homotopies. The goal of this section is to show that the converse of this statement holds up to natural equivalence if $\Cin$ admits either tensors or cotensors by $\Delta[1]$.

\begin{theorem}\label{theorem:MainTheoremSimplicial1}
Let $\Cin$ be a small simplicial category that admits tensors by $\Delta[1]$ or that admits cotensors by $\Delta[1]$. Then for any simplicial model category $\modN$, a functor $\Cin_0 \to \modN$ is weakly equivalent to a simplicial functor if and only if it is a homotopy functor.
\end{theorem}

If $\Cin$ admits (co)tensors by the standard simplex $\Delta[1]$, then homotopy functors can be characterized as follows.

\begin{lemma}\label{lemma:Alternative-characterization-homotopy-functor}
	In the situation of \cref{theorem:MainTheoremSimplicial1}, a functor $F \colon \Cin_0 \to \modN$ is a homotopy functor if and only if for every object $c$ in $\Cin$, the map $F(\Delta[1] \otimes c) \to F(c)$ is a weak equivalence (if $\Cin$ admits tensors by $\Delta[1]$) or for every object $c$ in $\Cin$, the map $F(c) \to F(c^{\Delta[1]})$ is a weak equivalence (if $\Cin$ admits cotensors by $\Delta[1]$).
\end{lemma}

\begin{proof}
	We treat the case where $\Cin$ admits tensors by $\Delta[1]$; the other case is dual. Since $\Delta[1] \to *$ is a simplicial homotopy equivalence, it follows that $\Delta[1] \otimes c \to c$ is a simplicial homotopy equivalence. This proves one direction. For the other direction, suppose that $F \colon \Cin_0 \to \modN$ takes $\Delta[1] \otimes c \to c$ to a weak equivalence in $\modN$ for every object $c$ in $\Cin$. To show that $F$ is a homotopy functor, it suffices to show that for any two strictly homotopic maps $f_0, f_1 \colon c \to d$, the maps $F(f_0)$ and $F(f_1)$ are identical in the homotopy category of $\modN$. If $f_0$ and $f_1$ are strictly homotopic, then without loss of generality, there exists a map $h \colon \Delta[1] \otimes c \to d$ such that $h\iota_k = f_k$ for $k=0,1$, where $\iota_k \colon c \to \Delta[1] \otimes c$ is the inclusion at the $k$-th vertex. Since $\iota_k$ is a right inverse to $\Delta[1] \otimes c \to c$, we see that $F(\iota_0)$ and $F(\iota_1)$ are both right inverses to the weak equivalence $F(\Delta[1] \otimes c) \to F(c)$. This implies that $F(\iota_0) = F(\iota_1)$ in the homotopy category of $\modN$, hence the result follows.
\end{proof}

Using that $\Delta[n]$ is a retract of $\Delta[1]^{n}$, we may assume without loss of generality that $\Cin$ admits (co)tensors by all standard simplices.

\begin{lemma}\label{lemma:Reduction-to-idempotent-completion}
	It suffices to prove \cref{theorem:MainTheoremSimplicial1} in the case where $\Cin$ admits cotensors by all the standard simplices $\Delta[n]$.
\end{lemma}

\begin{proof}
	By replacing $\Cin$ and $\modN$ with their opposite categories if necessary, it suffices to show that if \cref{theorem:MainTheoremSimplicial1} holds for any simplicial category $\Cin$ that admits tensors by all standard simplices $\Delta[n]$, then it also holds for any simplicial category $\Cin$ that only admits tensors by $\Delta[1]$. To this end, suppose that a simplicial category $\Cin$ that admits tensors by $\Delta[1]$ is given. Then $\Cin$ also admits tensors by $\Delta[1]^{n}$ for every $n$, since
	\begin{equation*}
		(\Delta[1]^n) \otimes c \cong {\underbrace{\Delta[1] \otimes \ldots \otimes \Delta[1]}_{n\text{ times}}} \otimes c.
	\end{equation*}
	Let $\wh \Cin$ denote the idempotent completion of $\Cin$. More precisely, we define $\wh\Cin$ to be the full simplicial subcategory of $\Fun(\Cin^\mathrm{op},\s\Set)$ spanned by those presheaves that are retracts of representables. Since $\Delta[n]$ is a retract of $\Delta[1]^n$, we see that the tensor $c \otimes \Delta[n]$ must exist in $\wh \Cin$ as a retract of $c \otimes \Delta[1]^n$. Now note that since $\modN$ is (co)complete, it is idempotent complete. In particular, the categories $\Fun(\Cin_0,\modN)$ and $\Fun(\Cin,\modN)$ are equivalent to $\Fun(\wh \Cin_0,\modN)$ and $\Fun(\wh \Cin, \modN)$, respectively. It follows from \cref{lemma:Alternative-characterization-homotopy-functor} that homotopy functors in $\Fun(\Cin_0,\modN)$ correspond to homotopy functors in $\Fun(\wh\Cin_0,\modN)$ under this equivalence of categories. Therefore, to deduce that \cref{theorem:MainTheoremSimplicial1} holds for the simplicial category $\Cin$, we may replace $\Cin$ with $\wh \Cin$.
\end{proof}

Note that by \Cref{lemma:BaseChangeGivenByKanExtension}, one can produce a simplicial functor from any ordinary functor $\Cin_0 \to \modN$ by left or right Kan extending along $\Cin_0 \to \Cin$, where we view $\Cin_0$ as a discrete simplicial category. We will prove \Cref{theorem:MainTheoremSimplicial1} by comparing a homotopy functor to its \emph{derived} left or right Kan extension along $\Cin_0 \to \Cin$, which by \cref{lemma:BarIsDerivedLeftKanExtension} can be computed explicitly using the (co)bar construction. Let $U$ denote the forgetful functor $\Fun(\Cin,\modN) \to \Fun(\Cin_0,\modN)$ (which we often omit from notation) and $L$ its left adjoint constructed using left Kan extension. Recall the definition of the derived (co)unit of the bar construction from \Cref{definition:DerivedUnitCounitBar}.

\begin{proposition}\label{proposition:SimplicialBarUnitCounitEquivalence}
Let $\Cin$ be a small simplicial category that admits cotensors by the standard simplex $\Delta[n]$ for every $n \geq 0$, and let $\modN$ be a simplicial model category with a cofibrant replacement functor $Q$. Then the derived unit
\begin{equation*} F \xleftarrow{\sim} B(\Cin_0,\Cin_0,QF) \xrightarrow{\eta} B(\Cin,\Cin_0,QF) \end{equation*}
is a weak equivalence for any homotopy functor $F \colon \Cin_0 \to \modN$ and the derived counit
\begin{equation*}B(\Cin,\Cin_0,QUG) \to LUG \xrightarrow{\varepsilon} G\end{equation*}
is a weak equivalence for any simplicial functor $G \colon \Cin \to \modN$.
\end{proposition}

\begin{remark}
By formally dualizing the statement of \Cref{proposition:SimplicialBarUnitCounitEquivalence}, one obtains a result about the cobar construction in the case that $\Cin$ admits tensors by the standard simplices.
\end{remark}

Before proving \Cref{proposition:SimplicialBarUnitCounitEquivalence}, let us show how to derive \Cref{theorem:MainTheoremSimplicial1} from it.

\begin{proof}[Proof of \Cref{theorem:MainTheoremSimplicial1}]
By \cref{lemma:Reduction-to-idempotent-completion}, we may assume without loss of generality that $\Cin$ admits cotensors by the standard simplices. Now suppose that $F \colon \Cin_0 \to \modN$ is a homotopy functor and write $Q$ for a cofibrant replacement functor on $\modN$. Then $B(\Cin,\Cin_0,QF)$ is a simplicial functor and $F \simeq B(\Cin_0,\Cin_0,QF) \simeq B(\Cin,\Cin_0,QF)$ by \Cref{proposition:SimplicialBarUnitCounitEquivalence}.
\end{proof}

The proof of \Cref{proposition:SimplicialBarUnitCounitEquivalence} will use the following lemmas about the bar construction and simplicial objects.

\begin{lemma}\label{lemma:BarConstructionsCommutesWithRealization}
Let $F \colon \Cin_0 \to \modN$ and $G \colon \Cin_0^\mathrm{op} \to \s\Set$ be functors, and let $G_n$ denote the functor $\Cin_0^\mathrm{op} \to \Set$ that sends an object $C$ to the set of $n$-simplices $G(C)_n$ of $G(C)$. Then $B(G,\Cin_0,F)$ is isomorphic to the geometric realization of the simplicial object $[n] \mapsto B(G_n, \Cin_0, F)$ in $\modN$, and if moreover $F$ is pointwise cofibrant, then this simplicial object is Reedy cofibrant.
\end{lemma}

\begin{proof}
The first statement follows since
\begin{equation*}B(G, \Cin_0, F) \cong B\left(\int^\Delta \Delta[n] \otimes G_n, \Cin_0, F\right) \cong \int^\Delta \Delta[n] \otimes B(G_n, \Cin_0, F) \end{equation*}
holds by \Cref{lemma:Fubini}.
For the second statement, note that the latching object $L_n B(G_\bullet, \Cin_0, F)$ agrees with $B(L_n G_\bullet, \Cin_0, F)$ by \Cref{lemma:BarCommutesWithColimits}. The map $B(L_n G_\bullet, \Cin_0, F) \to B(G_n, \Cin_0, F)$ is the inclusion of a summand into a coproduct of cofibrant objects, hence it is a cofibration.
\end{proof}

For the following lemmas, recall that a simplicial object $X_\bullet$ in a model category is called \emph{homotopically constant} if for any $f \colon [n] \to [m]$, the map $f^* \colon X_m \to X_n$ is a weak equivalence.

\begin{lemma}\label{lemma:HomotopicallyConstantObject}
Let $F \colon \Cin \to \modN$ be a homotopy functor and suppose that $\Cin$ admits cotensors by the standard simplices $\Delta[n]$. Then for any object $c$ in $\Cin$, the simplicial object $[n] \mapsto F(c^{\Delta[n]})$ is homotopically constant.
\end{lemma}

\begin{proof}
By the two-out-of-three property and the fact that $F$ is a homotopy functor, it suffices to show that for every $n \geq 0$, the map $c^{\Delta[0]} \to c^{\Delta[n]}$ is a simplicial homotopy equivalence. Now note that the map $\Delta[n] \to \Delta[0]$ is a simplicial deformation retract by \cite[Beispiel I.5.4]{Lamotke1968Semisimpliziale}. Since simplicial cotensors assemble into a simplicial functor by \Cref{construction:TensorsAreFunctorial}, we see that $c^{\Delta[0]} \to c^{\Delta[n]}$ is indeed a simplicial homotopy equivalence.
\end{proof}

\begin{lemma}\label{lemma:HomotopicallyConstantSimplicialReedyObject}
Let $N_\bullet$ be a homotopically constant Reedy cofibrant simplicial object in $\modN$. Then the canonical map $N_0 \to |N_\bullet|$ is a weak equivalence.
\end{lemma}

\begin{proof}
Let $cN_0$ denote the constant simplicial object on $N_0$. Then $cN_0 \to N_\bullet$ is a levelwise weak equivalence between Reedy cofibrant objects, hence $N_0 = |cN_0| \to |N_\bullet|$ is a weak equivalence.
\end{proof}

\begin{proof}[Proof of \Cref{proposition:SimplicialBarUnitCounitEquivalence}]
We first show that the derived unit is a weak equivalence for any homotopy functor $F \colon \Cin_0 \to \modN$. For this, it suffices to show that for any pointwise cofibrant $F$, the map $B(\Cin_0,\Cin_0,F) \to B(\Cin,\Cin_0,F)$ is a weak equivalence. Applying \Cref{lemma:BarConstructionsCommutesWithRealization} to $G = \Cin(-,c) \cong \{\Cin_0(-,c^{\Delta[n]})\}_{n \geq 0}$ yields that $B(\Cin,\Cin_0,F)$ is naturally isomorphic to the geometric realization of the Reedy cofibrant simplicial object
\begin{equation*} [n] \mapsto B(\Cin_0(-,c^{\Delta[n]}),\Cin_0,F). \end{equation*}
By \Cref{lemma:HomotopicallyConstantSimplicialReedyObject}, the map $B(\Cin_0,\Cin_0,F) \to B(\Cin,\Cin_0,F)$ is a weak equivalence if this simplicial object is homotopically constant. Since $B(\Cin_0(-,c^{\Delta[\bullet]}),\Cin_0,F) \simeq F(c^{\Delta[\bullet]})$ by \Cref{lemma:BarExtraDegeneracy}, it follows from \Cref{lemma:HomotopicallyConstantObject} that this is indeed the case.

To see that the derived counit is a weak equivalence, it suffices to show that the counit of
\begin{equation*}B(\Cin,\Cin_0,Q-) : \Ho(\Fun(\Cin_0,\modN)) \rightleftarrows \Ho(\Fun(\Cin,\modN)) : U \end{equation*}
is a natural isomorphism. This follows from the fact that $U$ creates weak equivalences, that any simplicial functor is a homotopy functor, and that the derived unit is a weak equivalence on homotopy functors.
\end{proof}

\begin{remark}
One can also use the ideas of \cite[\S 6]{RezkShwedeShipley2001Simplicial} to show that homotopy functors $\Cin \to \modN$ can be replaced by simplicial functors if $\Cin$ admits tensors or cotensors by standard simplices. To see this, let $R$ denote a simplicial Reedy fibrant replacement functor in $\modN^{\Delta^\mathrm{op}}$ and $Q$ a simplicial Reedy cofibrant replacement functor in $\modN^{\Delta}$, assuming these exist. A homotopy functor $F \colon \Cin \to \modN$ is then equivalent to the totalization of $[n] \mapsto RF(\Delta[n] \otimes -)$ or to the realization of $[n] \mapsto QF((-)^{\Delta[n]})$ if the relevant (co)tensors exist, and these functors can be shown to be simplicial functors. The advantage of this construction is that it also works for simplicial indexing categories $\Cin$ that are not small. However, to ensure the existence of these \textbf{simplicial} Reedy (co)fibrant replacement functors on $\modN^{\Delta^\mathrm{op}}$ and $\modN^{\Delta}$, one needs to put extra conditions on $\modN$ that we do not require in \Cref{theorem:MainTheoremSimplicial1}.
\end{remark}

We will now show that \Cref{theorem:MainTheoremSimplicial1} can be upgraded to a Quillen equivalence if the relevant model structures exist. Define (if it exists) the \emph{homotopy model structure} $\Funho(\Cin_0, \modN)$ on $\Fun(\Cin_0, \modN)$ to be the left Bousfield localization of the projective model structure in which the fibrant objects are the pointwise fibrant homotopy functors.

\begin{theorem}\label{theorem:MainTheoremSimplicial2}
Let $\Cin$ be a small simplicial category and $\modN$ a simplicial model category. If the projective model structures on $\Fun(\Cin, \modN)$ and $\Fun(\Cin_0, \modN)$ exist, then
\begin{equation*}L : \Fun(\Cin_0, \modN) \rightleftarrows \Fun(\Cin, \modN) : U  \end{equation*}
is a Quillen pair. If moreover $\Cin$ admits tensors or cotensors by $\Delta[1]$ and if the homotopy model structure $\Funho(\Cin_0,\modN)$ exists, then this Quillen pair becomes a Quillen equivalence between $\Funho(\Cin_0, \modN)$ and the projective model structure on $\Fun(\Cin, \modN)$.
\end{theorem}

\begin{remark}\label{remark:ExistenceHomotopyModelStructure}
The projective model structures on $\Fun(\Cin, \modN)$ and $\Fun(\Cin_0, \modN)$ always exist if $\modN$ is cofibrantly generated. Moreover, the homotopy model structure $\Funho(\Cin_0, \modN)$ exists whenever $\modN$ is left proper and combinatorial or cellular: In these cases, there always exists a set $\mathcal{T}$ of cofibrant objects in $\modN$ such that $N \to M$ is a weak equivalence if and only if $\Map(T,N) \to \Map(T,M)$ is a weak equivalence for every $T \in \mathcal{T}$. Here $\Map$ denotes the derived mapping space of $\modN$. The homotopy model structure is then obtained as the left Bousfield localization of $\Fun(\Cin_0,\modN)$ with respect to the set of maps
\begin{equation*}\{f^* \colon \Cin_0(d,-) \otimes T \to \Cin_0(c,-) \otimes T \mid f \colon c \to d \text{ is a homotopy equivalence in } \Cin \text{ and } T \in \mathcal{T} \}.\end{equation*}
\end{remark}

\begin{proof}[Proof of \Cref{theorem:MainTheoremSimplicial2}]
By the same proof as \cref{lemma:Reduction-to-idempotent-completion}, we may assume without loss of generality that $\Cin$ admits tensors by all standard simplices or cotensors by all standard simplices.

It is clear that the forgetful functor $U$ preserves pointwise (trivial) fibrations, hence that $L \dashv U$ is a Quillen pair with respect to the projective model structures. For the second statement, let us first assume that $\Cin$ admits cotensors by the standard simplices. Then \Cref{proposition:SimplicialBarUnitCounitEquivalence} shows that the derived unit of $L \dashv U$ is a weak equivalence for any homotopy functor $\Cin_0 \to \modN$, hence the result follows from \Cref{lemma:BousfieldLocalizationLemma}. If, on the other hand, $\Cin$ admits tensors by the standard simplices, then the proof of \Cref{lemma:BousfieldLocalizationLemma} only shows that $L \dashv U$ is a Quillen pair. To see that it is a Quillen equivalence, note that since $U$ preserves weak equivalences, its left and right derived functors agree. In particular, it suffices to show the left derived functor
\begin{equation*}\Ho(\Fun(\Cin,\modN)) \to \Ho(\Funho(\Cin_0,\modN))\end{equation*}
of $U$ is an equivalence of categories. This follows since by the dual of \Cref{lemma:BarIsDerivedLeftKanExtension}, the cobar construction can be used to derive the right adjoint of $U$, and by the dual of \Cref{proposition:SimplicialBarUnitCounitEquivalence}, the derived unit and counit of the resulting adjunction are weak equivalences.
\end{proof}

\begin{example}\label{example:FiniteSimplicialSets}
Let $\s\Set^{\mathrm{fin}}$ denote the category of finite simplicial sets, i.e.\ simplicial sets with finitely many non-degenerate simplices. Then the homotopy model structure on $\Fun((\s\Set^\mathrm{fin})_0, \s\Set)$ exists by \Cref{remark:ExistenceHomotopyModelStructure}, hence the category of simplicial functors $\Fun(\s\Set^\mathrm{fin}, \s\Set)$ equipped with the projective model structure is Quillen equivalent to the homotopy model structure $\Funho((\s\Set^\mathrm{fin})_0, \s\Set)$.\footnote{In \cite{Lydakis1998SimplicialFunctors}, Lydakis considers functors between $\finsSet$ and $\s\Set$ that preserve weak homotopy equivalences and calls them homotopy functors. Note that this terminology conflicts with ours, since in our terminology a homotopy functor only sends \textbf{simplicial} homotopy equivalences to weak homotopy equivalences.} There are of course many variations on this, such as a Quillen equivalence between $\Fun(\pfinsSet,\s\Set_*)$\footnote{Here we view $\pfinsSet$ as $\s\Set$-enriched and not as $\s\Set_*$-enriched.} and $\Funho((\pfinsSet)_0,\s\Set_*)$ to which we will return in \Cref{example:UnenrichedLydakisModelStructure}.
\end{example}

\begin{example}\label{example:FreelyAddingTensors}
Let $\Cin$ be any small simplicial category. One can ``freely'' add tensors by $\Delta[1]$ in the following way: First, identify $\Cin$ with its image under the enriched Yoneda embedding $\Cin \to \s\Set^{\Cin^\mathrm{op}}$; let us write $y_c$ for the presheaf represented by an object $c$ of $\Cin$. The category $\s\Set^{\Cin^\mathrm{op}}$ admits tensors by $\Delta[1]$, defined pointwise. If $\wt \Cin$ denotes the full subcategory of $\s\Set^{\Cin^\mathrm{op}}$ spanned by all objects of the form $\Delta[1]^n \otimes y_c$, then $\wt \Cin$ is a simplicial category that admits tensors $\Delta[1]$. \Cref{theorem:MainTheoremSimplicial2} yields a right Quillen equivalence $\Fun(\wt \Cin,\s\Set) \to \Funho(\wt \Cin_0,\s\Set)$, while the restriction $\Fun(\wt \Cin, \s\Set) \to \Fun(\Cin,\s\Set)$ is a right Quillen equivalence since $\Cin \to \wt \Cin$ is a DK-equivalence (cf.\ \Cref{definition:DKEquivalence} below).
\end{example}

In a cofibrantly generated simplicial model category, one can use the enriched small object argument to obtain simplicial fibrant and cofibrant replacement functors (cf.\ \cite[Thm.\ 13.2.1]{Riehl2014Categorical}). \Cref{theorem:MainTheoremSimplicial1} allows us to obtain a similar result for simplicial model categories $\modN$ that are not necessarily cofibrantly generated, but at the cost of restricting to a small full subcategory of $\modN$. Given a full subcategory $\modN' \subset \modN$, we define a \emph{cofibrant replacement functor on $\modN'$} as a functor $F \colon \modN' \to \modN$ landing in the full subcategory of cofibrant objects together with a natural weak equivalence $F \wearrow i$, where $i$ denotes the inclusion $\modN' \hookrightarrow \modN$. A \emph{fibrant replacement functor on $\modN'$} is defined dually.

\begin{corollary}\label{cor:SimplicialCofibrantReplacment}
Let $\modN$ be a simplicial model category. For any small full subcategory $\modN'$ of $\modN$, there exist simplicial fibrant and cofibrant replacement functors on $\modN'$.
\end{corollary}

\begin{proof}
We will only treat the case of the cofibrant replacement functors, since the statement about fibrant replacement functors is obtained by applying this case to $\modN^\mathrm{op}$. We may assume without loss of generality that $\modN'$ is closed under cotensors by the standard simplices $\Delta[n]$, since otherwise we can enlarge $\modN'$ to make this the case.

Let $Q$ be an (unenriched) cofibrant replacement functor on $\modN$. Applying $B(\modN',\modN'_0, Q-)$ to the inclusion $i \colon \modN' \hookrightarrow \modN$, we see by \Cref{proposition:SimplicialBarUnitCounitEquivalence} that there is a natural equivalence
\begin{equation*} B(\modN',\modN'_0,Qi) \wearrow i, \end{equation*}
while $B(\modN',\modN'_0,Qi)$ is pointwise cofibrant by \Cref{lemma:SimplicialBarIsReedyCofibrant}.
\end{proof}

\subsection{Continuous functors}\label{ssec:ContinuousFunctors}

We will now prove the topological counterparts of \Cref{theorem:MainTheoremSimplicial1,theorem:MainTheoremSimplicial2}. Throughout this section, $\Cin$ denotes a small topological category, $\modN$ a topological model category, and $\Cin_0$ the underlying category of $\Cin$. Note that just like a simplicial category, any topological category $\Cin$ comes with an obvious notion of \emph{homotopy} defined via paths in the hom-spaces, and we will call a functor $\Cin_0 \to \modN$ a \emph{homotopy functor} if it takes homotopy equivalences to weak equivalences. Clearly any continuous functor is a homotopy functor.

Recall the definition of goodness from \Cref{def:TopologicalGoodness}. As before, we define (if it exists) the \emph{homotopy model structure} $\Funho(\Cin_0,\modN)$ as the left Bousfield localization of the projective model structure on $\Fun(\Cin_0,\modN)$ in which the fibrant objects are the pointwise fibrant homotopy functors.

\begin{theorem}\label{theorem:MainTheoremTopological}
Let $\Cin$ be a small topological category that admits either tensors or cotensors by the unit interval, and let $\modN$ be a good topological model category. Then a functor $\Cin_0 \to \modN$ is weakly equivalent to a continuous functor if and only if it is a homotopy functor. Moreover, the forgetful functor
\begin{equation*}U \colon \Fun(\Cin, \modN) \to \Funho(\Cin_0, \modN)\end{equation*}
from the projective model structure to the homotopy model structure is a Quillen equivalence whenever the homotopy model structure exists.
\end{theorem}

\begin{proof}
Let $\Cin_\Delta$ denote the simplicial category $\Sing(\Cin)$. If $\Cin$ has (co)tensors by the unit interval $I$, then since $|\Delta[1]| \cong I$, it follows that $\Cin_\Delta$ admits (co)tensors by $\Delta[1]$.

By the goodness of $\modN$, we see that the projective model structures exist on $\Fun(\Cin_0,\modN)$, $\Fun(\Cin_\Delta,\modN)$, and $\Fun(\Cin,\modN)$. Moreover, by left Kan extending along the maps $C_0 \to |\Cin_\Delta| \to \Cin$ and their composite, we obtain the Quillen pairs
\begin{equation*}\begin{tikzcd}
\Fun(\Cin_\Delta,\modN) \ar[r, shift left=1] \ar[d, shift left = 1] & \Fun(\Cin,\modN), \ar[l, shift left=1] \ar[dl, shift left = 2, xshift=6] \\
 \Fun(\Cin_0,\modN) \ar[u, shift left=1] \ar[ur, shift right = 0, xshift=6] &
\end{tikzcd}\end{equation*}
where the left adjoints are drawn on the left or on top. The right adjoints in this diagram are the forgetful functors. By \Cref{proposition:GoodTopModelCatVsSimplicial}, the top horizontal adjunction is a Quillen equivalence. Moreover, it is clear that a functor $\Cin \to \modN$ is a homotopy functor if and only if the underlying simplicial functor $\Cin_\Delta \to \modN$ is a homotopy functor, so the result now follows by applying \Cref{theorem:MainTheoremSimplicial1,theorem:MainTheoremSimplicial2} to $\Fun(\Cin_\Delta,\modN)$.
\end{proof}

Note that many common topological model categories $\modN$ are good, cellular, and left proper; for example, all left Bousfield localizations of categories of diagrams in $\Top$ (cf.\ \Cref{proposition:GoodTopologicalInNature}). In these cases, the existence of the homotopy model structure $\Funho(\Cin_0,\modN)$ follows from the same argument as in \Cref{remark:ExistenceHomotopyModelStructure}, hence \Cref{theorem:MainTheoremTopological} provides us with a Quillen equivalence $\Funho(\Cin_0, \modN) \to \Fun(\Cin, \modN)$ whenever $\Cin$ admits tensors or cotensors by the unit interval.

\begin{example}\label{example:FiniteCWComplexes}
Let $\pfinCW$ denote the category of pointed finite CW-complexes, viewed as a $\Top$-enriched category, and write $\pfinCW_0$ for the underlying ordinary category. Then the homotopy model structure on $\Fun(\pfinCW_0, \Top_*)$ exists and is Quillen equivalent to the projective model structure on $\Fun(\pfinCW, \Top_*)$. We will see in \Cref{example:UnenrichedWSpectra} below that if one left Bousfield localizes the homotopy model structure further so that the fibrant objects become the linear functors, then one obtains a symmetric monoidal model category of spectra.
\end{example}

\begin{example}
Let $\Cin$ be any small $\Top$-enriched category. Analogously to \Cref{example:FreelyAddingTensors}, one can construct a category $\wt \Cin$ by freely adding tensors with the unit interval $I$ to $\Cin$. One then obtains right Quillen equivalences $\Fun(\Cin,\Top) \wearrowback \Fun(\wt \Cin, \Top) \wearrow \Funho(\wt \Cin_0,\Top)$.
\end{example}

\subsection{Digression: Dwyer--Kan localization}\label{ssec:DwyerKan}

Given a category $\Cin$ together with a subcategory $\calW \subset \Cin$, Dwyer and Kan in the series of papers \cite{DwyerKan1980SimplicialLocalizationsCategories,DwyerKan1980CalculatingSimplicialLocalizations,DwyerKan1980FunctionComplexesHomotopical,DwyerKan1987EquivalencesHomotopyDiagrams} developed methods to formally turn the maps in $\calW$ into homotopy equivalences, producing a simplicial category $L(\Cin,\calW)$ called its \emph{simplicial localization}. While their constructions are very concrete and have good formal properties, it is generally hard to compute these simplicial localizations explicitly in examples. In this section, we give an alternative proof of \cite[Prop.\ 1.3.4.7]{Lurie2017HigherAlgebra} which gives an explicit description of this simplicial localization in certain situations.

We will call a simplicial category $\Cin$ equipped with an ordinary subcategory $\calW$ a \emph{relative category}. A \emph{functor of relative categories} $(\Cin,\calW) \to (\Din,\calV)$ is a functor $F \colon \Cin \to \Din$ such that $F(\calW) \subset \calV$.

Let us recall the definition of a DK-equivalence.

\begin{definition}\label{definition:DKEquivalence}
	A functor $F \colon \Cin \to \Din$ of simplicial categories is called a \emph{Dwyer--Kan equivalence} or \emph{DK-equivalence} for short if
	\begin{enumerate}[(a),noitemsep]
		\item for any two objects $c,d \in \Cin$, the map $\Cin(c,d) \to \Din(Fc, Fd)$ is a weak homotopy equivalence, and
		\item for any $d \in \Din$, there exists a $c \in \Cin$ such that there is a simplicial homotopy equivalence $Fc \simeq d$.
	\end{enumerate}
\end{definition}

\begin{theorem}[{\cite[Prop.\ 1.3.4.7]{Lurie2017HigherAlgebra}}]\label{theorem:DwyerKanLocalizations}
Let $\Cin$ be a small simplicial category that either admits tensors by $\Delta[1]$ or cotensors by $\Delta[1]$, and let $\mathcal{W} \subset \Cin_0$ denote the subcategory of simplicial homotopy equivalences. Then the simplicial localization $L(\Cin_0, \mathcal{W})$ of $\Cin_0$ with respect to $\mathcal{W}$ is DK-equivalent to $\Cin$.
\end{theorem}

The idea of the proof is to identify both $L(\Cin_0,\calW)$ and $\Cin$ with the full subcategory of $\Funho(\Cin_0,\s\Set)^\mathrm{op}$ spanned by the representable functors through a ``derived Yoneda embedding''. To make this precise, we use Theorem 2.2 of \cite{DwyerKan1987EquivalencesHomotopyDiagrams}.

\begin{remark}
	One can get rid of the smallness assumption on $\Cin$ in \cref{theorem:DwyerKanLocalizations} by enlarging the universe.
\end{remark}

\begin{proof}[Proof of \cref{theorem:DwyerKanLocalizations}]
If $\Cin$ is a small simplicial category and $\calW$ is an ordinary subcategory, then by an argument similar to that of \Cref{remark:ExistenceHomotopyModelStructure}, the category $\Fun(\Cin,\s\Set)$ can be endowed with a model structure in which the cofibrations are the projective ones and in which the fibrant objects are precisely the pointwise fibrant functors that send the maps of $\calW$ to homotopy equivalences. Let us denote this model category by $\Fun^\calW(\Cin,\s\Set)$. It is proved in \cite[Thm.\ 2.2]{DwyerKan1987EquivalencesHomotopyDiagrams} that if $(\Cin,\calW) \to (\Din,\calV)$ is a functor of relative categories such that the induced Quillen pair
\begin{equation*}\Fun^\calW(\Cin,\s\Set) \rightleftarrows \Fun^\calV(\Din,\s\Set) \end{equation*}
is a Quillen equivalence, then the induced functor $L(\Cin,\calW) \to L(\Din,\calV)$ on simplicial localizations is a DK-equivalence ``up to retracts'' (cf.\ \cite[1.3.(iii)]{DwyerKan1987EquivalencesHomotopyDiagrams}).\footnote{The result of Dwyer--Kan is phrased differently since they did not have the general machinery of left Bousfield localizations in combinatorial model categories available yet, but it comes down to the same thing.} In particular, if $\Cin \to \Din$ is an isomorphism on objects, then this implies that $L(\Cin,\calW) \to L(\Din,\calV)$ is a DK-equivalence.

Now assume that $\Cin$ is a small simplicial category that either admits tensors or cotensors by $\Delta[1]$, and let $\calW$ be the full subcategory of simplicial homotopy equivalences. Then we have functors of relative categories
\begin{equation*}(\Cin,\varnothing) \to (\Cin,\calW) \leftarrow (\Cin_0,\calW) \end{equation*}
which induce right Quillen functors
\begin{equation*}\Fun(\Cin,\s\Set) \leftarrow \Fun^\calW(\Cin,\s\Set) \to \Fun^\calW(\Cin_0,\s\Set) \end{equation*}
The left-hand functor is a Quillen equivalence since it is the identity, while the right-hand functor is a Quillen equivalence by \Cref{theorem:MainTheoremSimplicial2}. We deduce from \cite[Thm.\ 2.2]{DwyerKan1987EquivalencesHomotopyDiagrams} that there are DK-equivalences $\Cin \simeq L(\Cin,\varnothing) \simeq L(\Cin,\calW) \simeq L(\Cin_0,\calW)$.
\end{proof}

\begin{remark}
	Strictly speaking, the statements of \cref{theorem:DwyerKanLocalizations} and \cite[Prop.\ 1.3.4.7]{Lurie2017HigherAlgebra} are different, so let us indicate why they are equivalent. It follows from \cite[Prop.\ 1.2.1]{Hinich2016DwyerKanLocalization} that a map $\Cin_0 \to \Cin$ exhibits $\Cin$ as the Dwyer--Kan localization of $\Cin_0$ at $\calW \subset \Cin_0$ precisely if the composite $(N(\Cin_0),\calW) \to N(\Cin)^\natural \to N(\Cin')^\natural$ is a weak equivalence of marked simplicial sets, where $\Cin \trivcofarrow \Cin'$ denotes a fibrant replacement of $\Cin$. Applying this to the case where $\calW$ consists of all homotopy equivalences, we obtain that \cref{theorem:DwyerKanLocalizations} is equivalent to saying that $(N(\Cin_0),\mathrm{hoeq}) \to N(\Cin ')^\natural$ is a weak equivalence of marked simplicial sets. By taking $\calW$ to be the homotopy equivalences in \cite[Prop.\ 1.3.4.7]{Lurie2017HigherAlgebra}, we see that this proposition implies \cref{theorem:DwyerKanLocalizations}.
	
	Conversely, if $\Cin$ and $\calW$ are as in \cite[Prop.\ 1.3.4.7]{Lurie2017HigherAlgebra}, then $\calW$ contains the canonical maps $\Delta[1] \otimes c \to c$. By the argument of \cref{lemma:Alternative-characterization-homotopy-functor}, any functor out of $\Cin_0$ sending the maps in $\calW$ to isomorphisms must send all homotopy equivalences to isomorphisms. In particular, we may assume without loss of generality that $\calW$ contains the subcategory $\mathrm{hoeq}$ of all homotopy equivalences. Since $(N(\Cin_0),\mathrm{hoeq}) \to N(\Cin')^\natural$ is a weak equivalence of marked simplicial sets by \cref{theorem:DwyerKanLocalizations}, it follows that $(N(\Cin_0),\calW) \to (N(\Cin'),\calW)$ must also be a weak equivalence of marked simplicial sets.
\end{remark}

\begin{example}
One can apply \Cref{theorem:DwyerKanLocalizations} to all examples considered in \Cref{ssec:ContinuousFunctors,ssec:SimplicialFunctors}. For example, the simplicial localization of the category $\wt \Cin_0$ constructed in \Cref{example:FreelyAddingTensors} with respect to the simplicial homotopy equivalences recovers the original category $\Cin$. Similarly, when applying \Cref{theorem:DwyerKanLocalizations} to the category of finite simplicial sets, one obtains the (to the author) surprising result that if one localizes $(\finsSet)_0$ at the simplicial homotopy equivalences, then one obtains the simplicial category $\finsSet$.
\end{example}

\begin{example}
Let $\mathbf{SC}^*$ denote the ($\Top_*$-enriched) category of separable $C^*$-algebras. Since this category has cotensors by the unit interval, its underlying simplicial category admits cotensors by the standard simplices as explained in the proof of \Cref{theorem:MainTheoremReducedTopological}. In particular, the simplicial localization $L(\mathbf{SC}^*_0,\mathrm{hoeq})$ of the underlying category $\mathbf{SC}^*_0$ at the homotopy equivalences is DK-equivalent to $\mathbf{SC}^*$, recovering Proposition 3.17 of \cite{BarneaJoachimMahanta2017ModelStructureProjective}.
\end{example}

%% file: Sections/ReducedFunctors.tex
\section{Reduced functors}\label{sec:ReducedFunctors}

We will now describe conditions under which one can replace a simplicial or continuous functor with one that is $\s\Set_*$- or $\Top_*$-enriched. In \Cref{ssec:PointedSimplicialFunctors}, we describe the simplicial case, from which we deduce the topological case in \Cref{ssec:PointedContinuousFunctors}. We will then extend these results to (symmetric) multifunctors in \Cref{ssec:MultireducedFunctors}, which are a type of functor that naturally occurs in Goodwillie calculus.

\subsection{Pointed simplicial functors}\label{ssec:PointedSimplicialFunctors}

Throughout this section, $\Cin_*$ denotes a small $\s\Set_*$-enriched category and $\modN$ a $\s\Set_*$-enriched model category. We will write $\Cin$ for the underlying $\s\Set$-enriched category of $\Cin_*$ obtained by forgetting the basepoints of all the hom-objects.

Note that in a $\s\Set_*$-enriched category $\Cin_*$, any initial object $X$ is automatically terminal and vice-versa, and that this is the case if and only if the basepoint of $\Cin_*(X,X)$ is the identity map $\id_X$. We call an object that is both initial and terminal a \emph{zero object}. If $\Cin$ is a simplicial category that has a zero object $*$,\footnote{This statement should be interpreted in the enriched sense and not just on the level of the underlying category, meaning that for every object $X$ of $\Cin$, there should be isomorphisms $\Cin(*,X) \cong * \cong \Cin(X,*)$ of simplicial sets.} then there is a unique way of upgrading $\Cin$ to a $\s\Set_*$-enriched category $\Cin_*$ by defining the basepoint of $\Cin(X,Y)$ to be the image of
\begin{equation}\label{eq:ZeroObjectDeterminesBasepoints}
* \cong \Cin(*,X) \wedge \Cin(X,Y) \to \Cin(X,Y). 
\end{equation}
In particular, a $\s\Set_*$-enriched model category is the same as a simplicial model category with a zero object. For this reason, we will call a $\s\Set_*$-enriched category with a zero object a \emph{pointed simplicial category}, and a $\s\Set_*$-enriched model category a \emph{pointed simplicial model category}. We will always denote zero objects by $*$.

\begin{definition}\label{definition:PointedReducedSimplicial} Let $\Cin$ be a pointed simplicial category, $\modN$ a pointed simplicial model category and $F \colon \Cin \to \modN$ a (not necessarily enriched) functor. Then
\begin{enumerate}[(i)]
    \item the functor $F$ is called \emph{pointed} if it preserves zero objects (i.e.\ $F(*) \cong *$), and
    \item the functor $F$ is called \emph{reduced} if the map $F(*) \to *$ is a weak equivalence in $\modN$.
\end{enumerate}
\end{definition}

\begin{remark}\label{remm:PointedFunctorIsPointedEnriched}
A $\s\Set_*$-enriched functor is by definition the same as a simplicial functor that preserves the basepoints of all the hom-objects. Since the zero object determines the basepoints of the hom-objects through the map \Cref{eq:ZeroObjectDeterminesBasepoints}, a simplicial functor between pointed simplicial categories is a $\s\Set_*$-enriched functor if and only if it is pointed.
\end{remark}

Clearly any pointed functor is a reduced functor. We will show that the converse holds up to natural equivalence; by the previous remark this means that any reduced functor $F \colon \Cin \to \modN$ can be replaced by a $\s\Set_*$-functor.

\begin{theorem}\label{theorem:MainTheoremReducedSimplicial1}
Let $\Cin$ be a small pointed simplicial category and let $\modN$ be a pointed simplicial model category. A simplicial functor $\Cin \to \modN$ is weakly equivalent to a pointed simplicial functor if and only if it is reduced.
\end{theorem}

The proof strategy is the same as for \Cref{theorem:MainTheoremSimplicial1}. Let $\Cin_+$ denote the $\s\Set_*$-category obtained by adjoining a disjoint basepoint to all hom-objects of $\Cin$, and let $L$ denote the left adjoint of the forgetful functor
\begin{equation*}U \colon \Fun(\Cin_*,\modN) \to \Fun(\Cin,\modN) .\end{equation*}
Note that $L$ is given by left Kan extension along $\Cin_+ \to \Cin_*$.

\begin{proposition}\label{proposition:PointedDerivedUnitCounitEquivalence}
Let $\Cin$ be a small pointed simplicial category and let $\modN$ be a pointed simplicial model category. Then
\begin{equation*}B(\Cin_+,\Cin_+,F) \to B(\Cin_*,\Cin_+, F) \end{equation*}
is a weak equivalence for any pointwise cofibrant reduced simplicial functor $F \colon \Cin \to \modN$, and the composite
\begin{equation*}B(\Cin_*,\Cin_+,UG) \to LUG \to G \end{equation*}
is a weak equivalence for any pointwise cofibrant pointed simplicial functor $G \colon \Cin_* \to \modN$.
\end{proposition}

\begin{remark}
We would like to call the two maps of this proposition the derived (co)units. However, since there might not exist a simplicial cofibrant replacement functor $Q$ defined on the whole category $\modN$, this would be slightly incorrect terminology.
\end{remark}

\begin{proof}[Proof of \Cref{proposition:PointedDerivedUnitCounitEquivalence}]
For the first equivalence, note that the square
\begin{equation*}\begin{tikzcd}
*_+ \cong \Cin_+(-,*) \ar[r,tail] \ar[d] & \Cin_+(-,d) \ar[d] \\
* \ar[r] & \Cin_*(-,d) \arrow[ul,phantom, very near start, "\ulcorner"]
\end{tikzcd}\end{equation*}
is a pushout. By \Cref{lemma:BarCommutesWithColimits} and \Cref{lemma:SimplicialBarIsReedyCofibrant}, the square
\begin{equation*}
\begin{tikzcd}
B(*_+,\Cin_+,F) \cong B(\Cin_+(*,-),\Cin_+,F) \ar[r,tail] \ar[d] & B(\Cin_+,\Cin_+,F) \ar[d] \\
* \cong B(*,\Cin_+,F) \ar[r] & B(\Cin_*,\Cin_+,F) \arrow[ul,phantom, very near start, "\ulcorner"]
\end{tikzcd}
\end{equation*}
is a pushout square in which the top horizontal arrow is a pointwise cofibration. Since the leftmost vertical map is a weak equivalence between pointwise cofibrant objects, we conclude that the rightmost vertical map is a weak equivalence, which was what we needed to show.

For the second map, consider the commutative triangle
\begin{equation*}\begin{tikzcd}
{UB(\Cin_*,\Cin_+,UG)} \arrow[r]           & UG. \\
{B(\Cin_+,\Cin_+,UG)} \arrow[u] \arrow[ru] &   
\end{tikzcd}\end{equation*}
The vertical map is a weak equivalence by the above, while the diagonal map is a weak equivalence by \Cref{lemma:BarExtraDegeneracy}, hence the horizontal map is a weak equivalence. Since the forgetful functor $U$ creates weak equivalences, the desired result follows.
\end{proof}

\begin{proof}[Proof of \Cref{theorem:MainTheoremReducedSimplicial1}]
Clearly any pointed functor is reduced. For the converse, let $F \colon \Cin \to \modN$ be a reduced simplicial functor. By \Cref{cor:SimplicialCofibrantReplacment}, we may assume without loss of generality that $F$ is pointwise cofibrant, so by \Cref{proposition:PointedDerivedUnitCounitEquivalence}, $B(\Cin_*,\Cin_+,F)$ is a pointed simplicial functor weakly equivalent to $F$.
\end{proof}

Again, this result can be upgraded to a Quillen equivalence if the required model structures exist. The \emph{reduced model structure} $\Funred(\Cin,\modN)$ is defined, if it exists, as the left Bousfield localization of the projective model structure on $\Fun(\Cin,\modN)$ in which the fibrant objects are the reduced pointwise fibrant functors.

\begin{theorem}\label{theorem:MainTheoremReducedSimplicial2}
Let $\Cin$ be a small pointed simplicial category and $\modN$ a simplicial model category. If the projective model structures on $\Fun(\Cin_*,\modN)$ and $\Fun(\Cin,\modN)$ exist, then
\begin{equation*}L : \Fun(\Cin,\modN) \rightleftarrows \Fun(\Cin_*,\modN) : U \end{equation*}
is a Quillen pair. Furthermore, if the reduced model structure $\Funred(\Cin,\modN)$ exists, then this adjunction becomes a Quillen equivalence between $\Funred(\Cin,\modN)$ and the projective model structure on $\Fun(\Cin_*,\modN)$.
\end{theorem}

\begin{proof}
    If the projective model structure on $\Fun(\Cin,\modN)$ exists, then we can use its cofibrant replacement functor $Q$ to replace any functor $F \colon \Cin \to \modN$ with the pointwise cofibrant functor $QF$. In particular, $B(\Cin_*,\Cin_+,Q-)$ is a left derived functor of $L \colon \Fun(\Cin,\modN) \to \Fun(\Cin_*,\modN)$. The rest of the proof is analogous to that of \Cref{theorem:MainTheoremSimplicial2}, using \Cref{proposition:PointedDerivedUnitCounitEquivalence} instead of \Cref{proposition:SimplicialBarUnitCounitEquivalence}.
\end{proof}

\begin{example}\label{example:FinitePointedSimplicialSets}
	Let $(\pfinsSet)_0$ denote the underlying category of $\pfinsSet$. Extending \Cref{example:FiniteSimplicialSets}, we obtain a Quillen equivalence
	\begin{equation}\label{SomeQuillenEquivalence}
		\Funhored((\pfinsSet)_0,\s\Set_*) \rightleftarrows \Fun(\pfinsSet,\s\Set_*). 
	\end{equation}
	Here $\Funhored((\pfinsSet)_0,\s\Set_*)$ denotes the left Bousfield localization of the projective model structure in which the fibrant objects are the pointwise fibrant homotopy functors that are also reduced, and $\Fun(\pfinsSet,\s\Set_*)$ denotes the category of $\s\Set_*$-enriched functors with the projective model structure.
\end{example}

\begin{example}\label{example:UnenrichedLydakisModelStructure}
	For brevity, let us write $\Cin_0$ for $(\pfinsSet)_0$. In the previous example, one can localize both sides of the Quillen equivalence further such that the fibrant objects become the pointwise fibrant functors that send weak homotopy equivalences to weak equivalences, that are reduced, and that are linear in the sense that for every finite pointed simplicial set $X$,
	\begin{equation*}\begin{tikzcd}
	F(X) \ar[r] \ar[d] & F(CX) \ar[d] \\
	F(*) \ar[r] & F(\Sigma X)
	\end{tikzcd}\end{equation*}
	is a homotopy pullback square. Here $CX$ denotes the reduced cone on $X$ and $\Sigma X$ the suspension. On the left-hand side of \Cref{SomeQuillenEquivalence}, this localization $\Funlin(\Cin_0,\s\Set_*)$ is obtained by localizing at the following two sets of maps:
	\begin{align*}
	&\{\Cin_0(Y,-)_+ \to \Cin_0(X,-)_+ \mid X \to Y\text{ is a weak homotopy equivalence}\}, \\
	&\{\Cin_0(*,-)_+ \cup^h_{\Cin_0(\Sigma X,-)_+} \Cin_0(CX,-)_+ \to \Cin_0(X,-)_+ \mid X \in \pfinsSet\}.
	\end{align*}
	Moreover, the Day convolution product on $\Fun(\Cin_0,\s\Set_*)$ endows this localization with the structure of a symmetric monoidal model category. If one localizes the right-hand side of \Cref{SomeQuillenEquivalence} at the analogous sets of maps, then one obtains Lydakis's stable model category of simplicial functors \cite{Lydakis1998SimplicialFunctors} together with a monoidal Quillen equivalence between $\Funlin(\Cin_0,\s\Set_*)$ and Lydakis's model category.
\end{example}

Finally, let us mention the following pointed analogue of \Cref{cor:SimplicialCofibrantReplacment}.

\begin{corollary}\label{cor:PointedCofibrantReplacment}
	Let $\modN$ be a pointed simplicial model category. For any small full subcategory $\modN'$ of $\modN$, there exist pointed simplicial fibrant and cofibrant replacement functors on $\modN'$.
\end{corollary}

\begin{proof}
	Assume without loss of generality that $\modN'$ contains the zero object of $\modN$. By \Cref{cor:SimplicialCofibrantReplacment}, there exists a simplicial cofibrant replacement functor $Q \colon \modN' \to \modN$. It follows as in the proof of \Cref{cor:SimplicialCofibrantReplacment} that  $B(\Cin_*, \Cin_+,Q(-))$ is a pointed simplicial cofibrant replacement functor $\modN' \to \modN$.
	
	A pointed simplicial fibrant replacement functor $\modN' \to \modN$ is constructed by applying the previous construction to $\modN^\mathrm{op}$.
\end{proof}

\subsection{Pointed continuous functors}\label{ssec:PointedContinuousFunctors}

In this section, we will deduce the topological analogue of the main result of \Cref{ssec:PointedSimplicialFunctors}. The proof strategy is similar to that of \Cref{theorem:MainTheoremTopological}, but note that it will be necessary to assume that the hom-spaces of the indexing category have nondegenerate basepoints.

We say that a topological category $\Cin$ is \emph{pointed} if it has a zero object. Similar to \Cref{ssec:PointedSimplicialFunctors}, there is a unique way of upgrading a pointed topological category $\Cin$ to a $\Top_*$-enriched category $\Cin_*$. Since a $\Top_*$-enriched model category has a zero object by definition, we will call such a model category a \emph{pointed topological model category}.

Pointed and reduced functors are defined as in \Cref{definition:PointedReducedSimplicial}; that is, $F \colon \Cin \to \modN$ is reduced if $F(*) \simeq *$ and pointed if $F(*) \cong *$. As in \Cref{remm:PointedFunctorIsPointedEnriched}, a pointed continuous functor $F \colon \Cin \to \modN$ is the same as a $\Top_*$-enriched functor $\Cin_* \to \modN$. The definition of the \emph{reduced model structure} $\Funred(\Cin,\modN)$ is the same as in the simplicial case.

\begin{theorem}\label{theorem:MainTheoremReducedTopological}
Let $\Cin$ be a small pointed topological category whose hom-spaces have nondegenerate basepoints, and let $\modN$ be a good pointed topological model category. Then a continuous functor $\Cin \to \modN$ is weakly equivalent to a pointed functor if and only if it is reduced. Moreover, the forgetful functor
\begin{equation*}U \colon \Fun(\Cin_*, \modN) \to \Funred(\Cin, \modN) \end{equation*}
from the projective model structure to the reduced model structure is a Quillen equivalence whenever the reduced model structure exists.
\end{theorem}

\begin{proof}
This is analogous to the proof of \Cref{theorem:MainTheoremTopological}.
\end{proof}

\begin{example}\label{example:UnenrichedWSpectra}
	Recall the category $\pfinCW$ and the Quillen equivalence 
	\begin{equation*}\Funho(\pfinCW_0,\Top_*) \leftrightarrows \Fun(\pfinCW_\mathbf{T},\Top_*)\end{equation*}
	from example \Cref{example:FiniteCWComplexes}. Here we view $\pfinCW$ as a $\Top_*$-category and $\pfinCW_\mathbf{T}$ denotes its underlying topological category. As in \Cref{example:FinitePointedSimplicialSets}, we obtain a Quillen equivalence
	\begin{equation*}\Funhored(\pfinCW_0,\Top_*) \leftrightarrows \Fun(\pfinCW,\Top_*),\end{equation*}
	and we can localize both sides of this Quillen equivalence as in \Cref{example:UnenrichedLydakisModelStructure} such that the fibrant objects become the linear reduced homotopy functors. The localization on the right-hand side then becomes the absolute stable model category of $\scrW$-spectra defined in \cite[\S 17]{MMSS2001DiagramSpectra}. Let us denote the localization on the left-hand side by $\Funlin(\pfinCW_0,\Top_*)$. As in \Cref{example:UnenrichedLydakisModelStructure}, one can verify that the Day convolution makes $\Funlin(\pfinCW_0,\Top_*)$ into a symmetric monoidal model category and that the Quillen equivalence with the absolute stable model category of $\scrW$-spectra is monoidal.
\end{example}

\subsection{(Symmetric) multireduced functors}\label{ssec:MultireducedFunctors}

In this section, we extend the previous results to symmetric multifunctors. Such functors occur naturally when studying the layers of the Taylor tower in Goodwillie calculus, which is the main motivation for discussing them here. We will first treat the case of non-symmetric multifunctors and then lift the obtained results to the case of symmetric multifunctors.

\begin{definition}\label{definition:Multifunctor}
	Let $\calV$ be a closed symmetric monoidal category. A \emph{multifunctor} is a $\calV$-functor of the form
	\begin{equation*}\Cin^{\otimes n} \to \Din,\end{equation*}
	where $\Cin$ and $\Din$ are $\calV$-categories and $\Cin^{\otimes n}$ denotes
	\begin{equation*}\Cin^{\otimes n} = \underbrace{\Cin \otimes \cdots \otimes \Cin}_{n \text{-times}}\end{equation*}
	(cf.\ \Cref{definition:TensorProductVCategories}).
\end{definition}

Throughout this section, we will work with topological categories. The results also hold in the simplicial case, in fact, even under slightly more general assumptions. We will come back to this in \Cref{remark:MultireducedSimplicialVersion} below. To avoid an overload of parentheses, we shall write $\Cin_*^{\wedge n}$ for $(\Cin_*)^{\wedge n}$ and $\Cin_+^{\wedge n}$ for $(\Cin_+)^{\wedge n}$ throughout this section.

Note that if $\Cin_*$ is a $\Top_*$-category, then the underlying topological category of $\Cin_*^{\wedge n}$ is not $\Cin^{\times n}$. This means that we can't apply the previous results directly and that we need a different definition of ``reduced'' in the case of a multifunctor.

\begin{definition}\label{definition:MultipointedMultireduced}
	Let $\Cin_*$ be a pointed topological category, $\modN$ a pointed topological model category and $F \colon \Cin^{\times n} \to \modN$ a continuous multifunctor. Then
	\begin{enumerate}[(i)]
		\item the multifunctor $F$ is called \emph{(multi)pointed} if $F(c_1,\ldots,c_n) \cong *$ whenever $c_i = *$ for some $1 \leq i \leq n$, and
		\item the multifunctor $F$ is called \emph{(multi)reduced} if $F(c_1,\ldots,c_n) \simeq *$ whenever $c_i = *$ for some $1 \leq i \leq n$.
	\end{enumerate}
\end{definition}

\begin{remark}
	As in \Cref{remm:PointedFunctorIsPointedEnriched}, a $\Top_*$-functor $\Cin_*^{\wedge n} \to \modN$ is the same as a multipointed topological functor $F \colon \Cin^{\times n} \to \modN$.
\end{remark}

The \emph{multireduced model structure} $\Funmred(\Cin^{\times n},\modN)$ is defined, if it exists, as the left Bousfield localization of the projective model structure on $\Fun(\Cin^{\times n},\modN)$ in which the fibrant objects are the multireduced pointwise fibrant functors.

\begin{theorem}\label{theorem:MainTheoremMultireduced}
	Let $\Cin$ be a small pointed topological category whose hom-spaces have nondegenerate basepoints, and let $\modN$ be a good pointed topological category. Then a continuous multifunctor $\Cin^{\times n} \to \modN$ is weakly equivalent to a multipointed functor if and only if it is multireduced. Moreover, the forgetful functor
	\begin{equation*}U \colon \Fun(\Cin_*^{\wedge n}, \modN) \to \Funmred(\Cin^{\times n}, \modN) \end{equation*}
	from the projective model structure to the multireduced model structure is a Quillen equivalence whenever the multireduced model structure exists.
\end{theorem}

\begin{proof}
	The forgetful functor can be viewed as the restriction along
	\begin{equation}\label{map:MultiContinuousToMultiPointedIndex}
	(\Cin^{\times n})_+ \cong \Cin_+^{\wedge n} \to \Cin_*^{\wedge n}.
	\end{equation}
	In particular, left Kan extension along this map provides a left adjoint to the forgetful functor.
	
	Now let $F \colon \Cin^{\times n}  \to \modN$ be multireduced and assume without loss of generality that $F$ is cofibrant in the projective model structure on $\Fun(\Cin^{\times n},\modN)$. We will show that the left Kan extension of $F$ along \Cref{map:MultiContinuousToMultiPointedIndex} is weakly equivalent to $F$. In order to see this, note that this map factors as
	\begin{equation*}\Cin_+^{\wedge n} \to \cdots \to \Cin_+^{\wedge (n-k+1)} \wedge \Cin_*^{\wedge (k-1)} \to \Cin_+^{\wedge (n-k)} \wedge \Cin_*^{\wedge k} \to \cdots \to \Cin_*^{\wedge n}.  \end{equation*}
	In particular, we can prove the claim for each of these maps separately, so it suffices to show that the left Kan extension of a reduced functor
	\begin{equation*}G \colon \Cin_+ \to \Fun(\Cin_+^{\wedge (n-k)} \wedge \Cin_*^{\wedge (k-1)},\modN)\end{equation*}
	along $\Cin_+ \to \Cin_*$ is weakly equivalent to $G$. This was proved in \Cref{theorem:MainTheoremReducedTopological}.
	
	The claim about the Quillen equivalence follows directly from \Cref{lemma:BousfieldLocalizationLemma}.
\end{proof}

We now extend this result to symmetric multifunctors.

\begin{definition}\label{definition:SymmetricMultifunctor}
	Let $\calV$ be a closed symmetric monoidal category and let $\Cin$ and $\Din$ be $\calV$-categories. A \emph{symmetric multifunctor} is a multifunctor $F \colon \Cin^{\otimes n} \to \Din$ equipped with $\calV$-natural isomorphisms
	\begin{equation*}\sigma_F \colon F(c_1,\ldots,c_n) \cong F(c_{\sigma(1)},\ldots,c_{\sigma(n)}) \end{equation*}
	for every $\sigma \in \Sigma_n$, satisfying $(\tau\sigma)_F = \sigma_F \tau_F$. A \emph{symmetric $\calV$-natural transformation} is a $\calV$-natural transformation between symmetric functors that is compatible with the symmetry, and the category of symmetric multifunctors and such natural transformations is denoted $\Funsym(\Cin^{\otimes n},\Din)$.
\end{definition}

\begin{remark}
	In the definition above, it follows automatically that for the identity permutation $e \in \Sigma_n$, one has $e_F = \id$.
\end{remark}

In the proof of \Cref{theorem:MainTheoremSymmetricReduced} below, we will use a different description of the category $\Funsym(\Cin^{\otimes n},\modN)$, namely as the homotopy fixed points of $\Fun(\Cin^{\otimes n},\modN)$. Note that  $\Sigma_n$ acts from the left on $\Cin^{\otimes n}$ by permuting the factors, giving a right $\Sigma_n$-action on $\Fun(\Cin^{\otimes n},\modN)$. Define $E\Sigma_n$ to be the groupoid whose objects are the elements of $\Sigma_n$ and which has exactly one arrow between any pair of objects. $\Sigma_n$ acts on this groupoid by multiplication from the right. The category $\Funsym(\Cin^{\otimes n}, \Din)$ is easily seen to be equivalent to the category of $\Sigma_n$-equivariant functors
\begin{equation}\label{map:HomotopyFixedPointsIsSymmetricFunctors}
E\Sigma_n \to \Fun(\Cin^{\otimes n}, \Din)
\end{equation}
and $\Sigma_n$-equivariant $\calV$-natural transformations between such functors.

In \cite[Lem.\ 3.6]{BiedermannRondigs2014CalculusII}, yet another description of $\Funsym(\Cin^{\otimes n},\Din)$ is given. Namely, it is shown that $\Funsym(\Cin^{\otimes n},\Din)$ is equivalent to $\Fun(\Sigma_n \wr \Cin^{\otimes n},\Din)$ for a certain $\calV$-category $\Sigma_n \wr \Cin^{\otimes n}$, namely the \emph{wreath product category} (cf.\ \cite[Def.\ 3.3]{BiedermannRondigs2014CalculusII}). In particular, we see that  $\Funsym(\Cin^{\otimes n},\Din)$ is an ordinary functor category, so the usual methods for constructing projective model structures apply to it. In the theorem below, the \emph{symmetric multireduced model structure} $\Funsymred(\Cin^{\times n},\modN)$ is defined, if it exists, as the left Bousfield localization of the projective model structure on $\Funsym(\Cin^{\times n},\modN)$ in which the fibrant objects are the symmetric multireduced pointwise fibrant functors.

\begin{theorem}\label{theorem:MainTheoremSymmetricReduced}
	Let $\Cin$ be a small pointed topological category whose hom-spaces have nondegenerate basepoints, and let $\modN$ be a good pointed topological category. Then a continuous symmetric multifunctor $\Cin^{\times n} \to \modN$ is weakly equivalent (as symmetric multifunctors) to a symmetric multipointed functor if and only if it is multireduced. Moreover, the forgetful functor
	\begin{equation*}U \colon \Funsym(\Cin_*^{\wedge n}, \modN) \to \Funsymred(\Cin^{\times n}, \modN) \end{equation*}
	from the projective model structure to the symmetric multireduced model structure is a Quillen equivalence whenever the symmetric multireduced model structure exists.
\end{theorem}

\begin{proof}
	The main idea of this proof is that one can work underlying in $\Fun(\Cin^{\times n},\modN)$ and $\Fun(\Cin_*^{\wedge n},\modN)$, where we have already proved the result. Since $\Cin_+^{\wedge n} \to \Cin_*^{\wedge n}$ is a $\Sigma_n$-equivariant functor, the adjunction
	\begin{equation}\label{adjunction:LeftKanExtensionIsEquivariant}
	L : \Fun(\Cin^{\times n},\modN) \rightleftarrows \Fun(\Cin_*^{\wedge n},\modN) : U
	\end{equation}
	obtained by restriction and left Kan extension along this functor is $\Sigma_n$-equivariant.
	By identifying symmetric functors with $\Sigma_n$-equivariant functors out of $E \Sigma_n$ as in \Cref{map:HomotopyFixedPointsIsSymmetricFunctors}, we obtain an adjunction
	\begin{equation}\label{adjunction:LeftKanExtensionSymmetricVersion}
		L : \Funsym(\Cin^{\times n},\modN) \rightleftarrows \Funsym(\Cin_*^{\wedge n},\modN) : U
	\end{equation}
	where the functors are computed underlying as in \Cref{adjunction:LeftKanExtensionIsEquivariant} (hence the abuse of notation).
	Furthermore, note that the forgetful functor
	\begin{equation}\label{map:ForgetfulFunctorSymmetric}
	\Funsym(\Cin^{\times n},\modN) \to \Fun(\Cin^{\times n},\modN)
	\end{equation}
	admits a right adjoint given by $F \mapsto \prod_{\sigma \in \Sigma_n} \sigma^* F$. Since this right adjoint preserves pointwise fibrations, we see that the forgetful functor \Cref{map:ForgetfulFunctorSymmetric} preserves projectively cofibrant functors. This implies that the derived adjunction of \Cref{adjunction:LeftKanExtensionSymmetricVersion} can be computed underlying as the derived adjunction of \Cref{adjunction:LeftKanExtensionIsEquivariant}.
	By the proof of \Cref{theorem:MainTheoremMultireduced}, the derived unit of this adjunction is a weak equivalence on any multireduced functor.
	
	The claim about the Quillen equivalence again follows directly from \Cref{lemma:BousfieldLocalizationLemma}.
\end{proof}

\begin{remark}\label{remark:MultireducedSimplicialVersion}
	If the relevant projective model structures exist, then the above proofs also go through for simplicial functors. However, if these projective model structures do not exist, then it is still true that any (symmetric) multireduced simplicial functor $F$ is equivalent to a (symmetric) multipointed one. Namely, if one uses \Cref{cor:PointedCofibrantReplacment} to replace $F$ with a pointwise cofibrant functor, then the proofs of \Cref{theorem:MainTheoremMultireduced,theorem:MainTheoremSymmetricReduced} still go through if one uses the bar construction instead of left Kan extension everywhere.
\end{remark}

\begin{remark}\label{remark:HomotopyMultiFunctors}
	One can also prove similar results for homotopy functors. Namely, if $\Cin$ is a topological or simplicial category, then a multifunctor $(\Cin_0)^{\times n} = (\Cin^{\times n})_0 \to \modN$ is a homotopy functor if and only if it sends homotopy equivalences to weak equivalences in each variable separately. In particular, the analogue of \Cref{theorem:MainTheoremMultireduced} follows directly from \Cref{theorem:MainTheoremSimplicial1,theorem:MainTheoremSimplicial2,theorem:MainTheoremTopological}. This can then be upgraded to a statement about symmetric multifunctors by the same proof as \Cref{theorem:MainTheoremSymmetricReduced}.
\end{remark}

%% file: Sections/LinearFunctors.tex
\section{Linear functors}\label{sec:LinearFunctors}

In the final part of this paper, we will study in which cases a given functor is equivalent to one that respects an enrichment in the category $\Sp$ of orthogonal spectra. In \Cref{ssec:LinearOrthogonalSpectraCase}, under the assumption that the indexing category has cotensors by $\SSS^{-1}$, we describe a condition called \emph{strict linearity} (cf.\ \Cref{definition:LinearOrthogonal}) and show that a functor is equivalent to an $\Sp$-functor if and only if it satisfies this condition. In particular, we prove \Cref{TheoremC}. In \Cref{ssec:LinearOrthogonalSpectraCaseTensored}, we instead study the case where the indexing category admits tensors by $\SSS^{-1}$. Unfortunately, one cannot simply dualize the results of \Cref{ssec:LinearOrthogonalSpectraCaseTensored}, so this case requires a separate treatment with slightly more restrictive assumptions. We conclude by proving analogous results for (symmetric) functors of multiple variables in \Cref{ssec:MultilinearFunctors}.

\subsection{Linear functors in the presence of cotensors with \texorpdfstring{$\SSS^{-1}$}{S\textasciicircum -1}}\label{ssec:LinearOrthogonalSpectraCase}

Recall the category of orthogonal spectra $\Sp$ described in \Cref{ssec:ModelCats}. We will write $\SSS^{-n}$ for the image of $n \in \OIndex$ under the Yoneda embedding $\OIndex^\mathrm{op} \hookrightarrow \Sp$.\footnote{To avoid awkward notation, we will write $\SSS$ or $\SSS^{0}$ and not $\SSS^{-0}$ for the image of $0 \in \Sigma_S$ under the Yoneda embedding.} Given a spectral category $\Cin$, we will write $\Cin_*$ for its underlying $\Top_*$-category.

\begin{construction}\label{construction:Sigma1}
	Let $\Cin$ be a spectral category that admits cotensors by $\SSS^{-1}$ and let $\modN$ be a spectral model category. Since $\SSS^{-k} \otimes \SSS^{-1} = \SSS^{-k-1}$, we see inductively that $\Cin$ admits cotensors by $\SSS^{-n}$ for every $n \geq 0$. By \Cref{construction:TensorsAreFunctorial}, these cotensors assemble into a functor $\SSSOIndex \otimes \Cin \to \Cin$, where $\SSSOIndex$ is the full spectral subcategory of $\Sp^\mathrm{op}$ on the objects $\{\SSS^{-n}\}_{n \geq 0}$. By the enriched Yoneda lemma, the underlying $\Top_*$-category of $\SSSOIndex$ is $\OIndex$. In particular, restricting the cotensor functor to the underlying $\Top_*$-categories yields a functor $\OIndex \wedge \Cin_* \to \Cin_*$. Composing this functor with a $\Top_*$-functor $F \colon \Cin_* \to \modN$ yields maps $\OIndex(0,n) \to \modN(F(c),F(c^{\SSS^{-n}}))$ and hence, by adjunction, maps
	\begin{equation*}
		\sigma_n \colon \Sigma^n F(c) = \OIndex(0,n) \otimes F(c) \to F(c^{\SSS^{-n}})
	\end{equation*}
	that are natural in $c$.
\end{construction}

\begin{definition}\label{definition:LinearOrthogonal}
	 Let $\Cin$ be a spectral category that admits cotensors by $\SSS^{-1}$, let $\modN$ be a spectral model category, and let $F \colon \Cin_* \to \modN$ be a $\Top_*$-functor. Then $F$ is \emph{strictly linear} if for every $c$ in $\Cin$, the composite
	 \begin{equation*}
	 	(\mathbb{L} \Sigma) F(c) \to \Sigma F(c) \xrightarrow{\sigma_1} F(c^{\SSS^{-1}})
	 \end{equation*}
 	is a weak equivalence in $\modN$.
\end{definition}

By $\mathbb{L} \Sigma$, we mean the left derived functor of the suspension functor $\Sigma = S^1 \otimes - \colon \modN \to \modN$.

\begin{remark}\label{remark:Sigma1EquivalenceImpliesSigmanEquivalence}
	Using the associativity of the cotensor, one can show that $\sigma_{n+1} = \sigma_1 \circ \Sigma \sigma_n$. In particular, if $F \colon \Cin_* \to \modN$ is strictly linear, then it follows inductively that $(\mathbb{L} \Sigma^n)F(c) \to F(c^{\SSS^{-n}})$ is a weak equivalence for every $n \geq 0$.
\end{remark}

\begin{example}
	The spectral category $\SSSOIndex$ from \Cref{construction:Sigma1} is dual to the full subcategory of $\Sp$ spanned by $\{\SSS^{-n}\}$, hence it admits cotensors by $\SSS^{-1}$. A functor $X = \{X_n\}_{n \in\bbN} \colon \OIndex \to \modN$ is strictly linear if and only if for every $n \geq 0$, the map $\mathbb L \Sigma X_n \to X_{n+1}$ is a weak equivalence. Since $\modN$ is stable by \Cref{lemma:SpectrallyEnrichedIsStable}, this is equivalent to $X_n \to \mathbb R \Omega X_{n+1}$ being a weak equivalence. That is, strictly linear functors $\OIndex \to \modN$ are the same as orthogonal $\Omega$-spectrum objects in $\modN$.
\end{example}

The main result of this section is the following.

\begin{theorem}\label{theorem:MainTheoremLinearOrthogonal1}
	Let $\Cin$ be a small spectral category that admits cotensors by $\SSS^{-1}$, let $\modN$ be a good spectral model category, and assume that the hom-spaces of the underlying $\Top_*$-category $\Cin_*$ of $\Cin$ are nondegenerately based. Then a $\Top_*$-functor $\Cin_* \to \modN$ is equivalent to a spectral functor if and only if it is strictly linear.
\end{theorem}

Let us start by showing that spectral functors are strictly linear. We will use the following lemmas.

\begin{lemma}\label{lemma:SuspensionsProjectivelyCofibrantFunctors}
	Let $\Cin$ be a small spectral category and $\modN$ a good spectral model category. Denote the left derived functor of $\Sigma \colon \Fun(\Cin,\modN) \to \Fun(\Cin,\modN)$ by $\mathbb{L} \Sigma_P$ and the left derived functor of $\Sigma \colon \modN \to \modN$ by $\mathbb{L}\Sigma_\modN$. For any $F \colon \Cin \to \modN$ and any $c$ in $\Cin$, one has $(\mathbb{L} \Sigma_P)F(c) \simeq (\mathbb{L} \Sigma_\modN)F(c)$.
\end{lemma}

\begin{proof}
	The subtlety of this statement lies in the fact that projectively cofibrant functors $F \colon \Cin \to \modN$ need not be pointwise cofibrant. To see that nonetheless $(\mathbb{L} \Sigma_P)F(c) \simeq (\mathbb{L} \Sigma_\modN)F(c)$, note that by \Cref{lemma:CofibrantReplacementOfCategories} there exists a weak equivalence of spectral categories $\wh \Cin \wearrow \Cin$ where $\wh\Cin$ has cofibrant hom-objects. By the goodness of $\modN$, restriction along this weak equivalence induces a Quillen equivalence
	\begin{equation*}
		L : \Fun(\wh\Cin,\modN) \rightleftarrows \Fun(\Cin,\modN) : U.
	\end{equation*}
	In particular, one can write $F \simeq LG$ for a projectively cofibrant functor $G \colon \wh\Cin \to \modN$. Since $LG$ is projectively cofibrant, we have $(\mathbb{L} \Sigma_P)F(c) \simeq \Sigma L G(c)$. Because $\wh\Cin$ has cofibrant hom-objects, the functor $G$ is pointwise cofibrant, hence $G(c)$ is a cofibrant replacement of $F(c)$. Since $\Sigma G \wearrow UL(\Sigma G) \cong \Sigma ULG$, we have $(\mathbb{L} \Sigma_\modN) F(c) \simeq \Sigma G(c) \simeq \Sigma LG(c)$.
\end{proof}

\begin{lemma}\label{lemma:TwoCanonicalMapsAreEquivalences}
	The canonical maps $\Sigma \SSS^{-1} \to (\SSS^{-1})^{\SSS^{-1}}$ and $\SSS^{0} \to (\SSS^{-1})^{\SSS^{-1}}$ are stable equivalences.
\end{lemma}

\begin{proof}
	Note that cotensors by $\SSS^{-1}$ in $\Sp$ are the same as shifts. In particular, the map $\Sigma \SSS^{-1} \to (\SSS^{-1})^{\SSS^{-1}}$, constructed as in \Cref{construction:Sigma1}, is a stable equivalence by \cite[Prop.\ 3.1.25.(ii)]{Schwede2018GlobalHomotopyTheory} (taking $G$ to be the trivial group, $X = \SSS^{-1}$, and $V = \mathbb{R}$ in that proposition).
	
	For the second map, note that the functor $(\mhyphen)^{\SSS^{-1}} = \mathrm{sh}$ preserves stable equivalences by \cite[Prop.\ 3.1.25.(ii)]{Schwede2018GlobalHomotopyTheory} and that it is a right Quillen equivalence. In particular, we see that $\SSS^{0} \to (\SSS^{-1})^{\SSS^{-1}}$ is a stable equivalence if and only if its adjunct $\id \colon \SSS^{-1} \to \SSS^{-1}$ is a stable equivalence, which is the case since it is the identity map.
\end{proof}

\begin{proposition}\label{proposition:SpectralFunctorIsStrictlyLinear}
	Let $\Cin$ and $\modN$ be as in \Cref{theorem:MainTheoremLinearOrthogonal1} and let $F \colon \Cin \to \modN$ be a spectral functor. Then $F$ is strictly linear.
\end{proposition}

\begin{proof}
	By the goodness of $\modN$ and the fact that the strict linearity of functors $\Cin \to \modN$ is a homotopy invariant property, we may assume without loss of generality that $F$ is projectively cofibrant. By \Cref{lemma:SuspensionsProjectivelyCofibrantFunctors}, $F$ is strictly linear if and only if for every $c$ in $\Cin$, the map $\sigma_1 \colon \Sigma F(c) \to F(c^{\SSS^{-1}})$ is a weak equivalence in $\modN$. We will denote the object $\SSS^{-n}$ of the spectral category $\SSSOIndex$ defined in \Cref{construction:Sigma1} by $n$. By construction, $\sigma_1$ agrees with the map
	\begin{equation*}
		\SSS^1 \otimes F(c) \cong \SSSOIndex(0,1) \otimes F(c) \to F(c^{\SSS^{-1}}).
	\end{equation*}
	Since $\SSSOIndex(1,0) \cong \SSS^{-1}$, $\SSSOIndex(0,0) \cong \SSS^0$, and $\SSSOIndex(1,1) \cong (\SSS^{-1})^{\SSS^{-1}}$, we obtain a commutative diagram
	\begin{equation}\label{diag:SpectralImpliesLinear}\begin{tikzcd}[column sep = large]
	\SSS^1 \otimes \SSS^{-1} \otimes \SSS^1 \otimes F(c) \ar[r, "\id \otimes \id \otimes \ac"] \ar[d, "\id \otimes \comp \otimes \id"] & \SSS^1 \otimes \SSS^{-1} \otimes F(c^{\SSS^{-1}}) \ar[r, "\comp \otimes \id"] \ar[d, "\id \otimes \ac"] & (\SSS^{-1})^{\SSS^{-1}} \otimes F(c^{\SSS^{-1}}) \ar[d, "\ac"] \\
	\SSS^1 \otimes \SSS^0 \otimes F(c) \ar[r, "\id \otimes \ac"] & \SSS^1 \otimes F(c) \ar[r,"\sigma_1"] & F(c^{\SSS^{-1}}).
	\end{tikzcd}\end{equation}
	The leftmost vertical map is a weak equivalence since $F$ is projectively cofibrant and $c \colon \SSS^{-1} \otimes \SSS^{1} \to \SSS^{0}$ is a stable equivalence between cofibrant orthogonal spectra. Moreover, the bottom left horizontal map is an isomorphism. To see that the composition of the top right horizontal map and the rightmost vertical map is a weak equivalence, let $\wh{\SSSOIndex} \trivfibarrow \SSSOIndex$ be a cofibrant replacement in the sense of \Cref{lemma:CofibrantReplacementOfCategories}. In particular, the hom-object $\wh{\SSSOIndex}(1,1)$ is cofibrant. We then obtain the commutative diagram
	\begin{equation*}\begin{tikzcd}
		\wh{\SSSOIndex}(0,1) \otimes \wh{\SSSOIndex}(1,0) \otimes F(c^{\SSS^{-1}}) \ar[r,"c \otimes \id"] \ar[d] & \wh{\SSSOIndex}(1,1) \otimes F(c^{\SSS^{-1}}) \ar[r,"\ac"] \ar[d] & F(c^{\SSS^{-1}}) \ar[d,"\id"] \\
		\SSS^{1} \otimes \SSS^{-1} \otimes F(c^{\SSS^{-1}}) \ar[r,"c \otimes \id"] & (\SSS^{-1})^{\SSS^{-1}} \otimes F(c^{\SSS^{-1}}) \ar[r,"\ac"] & F(c^{\SSS^{-1}})
	\end{tikzcd}\end{equation*}
	The leftmost vertical map is a weak equivalence since $\wh{\SSSOIndex}(0,1) \otimes \wh{\SSSOIndex}(1,0) \to \SSS^{1} \otimes \SSS^{-1}$ is a stable equivalence between cofibrant orthogonal spectra. Moreover, $\wh{\SSSOIndex}(0,1) \otimes \wh{\SSSOIndex}(1,0) \to \wh{\SSSOIndex}(1,1)$ is a stable equivalence between cofibrant orthogonal spectra since $\SSS^{1} \otimes \SSS^{-1} \to (\SSS^{-1})^{\SSS^{-1}}$ is a stable equivalence by \Cref{lemma:TwoCanonicalMapsAreEquivalences}. This shows that the top left horizontal map is a weak equivalence. Finally, to see that the top right horizontal map is a weak equivalence, consider the diagram
	\begin{equation*}\begin{tikzcd}
			\SSS^{0} \otimes F(c^{\SSS^{-1}}) \ar[r] \ar[dr,"\cong"'] & \wh{\SSSOIndex}(1,1) \otimes F(c^{\SSS^{-1}}) \ar[d] \\
			& F(c^{\SSS^{-1}})
	\end{tikzcd}\end{equation*}
	The unit map $\SSS^{0} \to \wh{\SSSOIndex}(1,1)$ is a weak equivalence between cofibrant objects since $\SSS^{0} \to (\SSS^{-1})^{\SSS^{-1}}$ is a weak equivalence by \Cref{lemma:TwoCanonicalMapsAreEquivalences}, hence the map $\wh{\SSSOIndex}(1,1) \otimes F(c^{\SSS^{-1}}) \to F(c^{\SSS^{-1}})$ is a weak equivalence by the two-out-of-three property. In particular, the composition
	\begin{equation*}
		\SSS^{1} \otimes \SSS^{-1} \otimes F(c^{\SSS^{-1}}) \xrightarrow{c \otimes \id} (\SSS^{-1})^{\SSS^{-1}} \otimes F(c^{\SSS^{-1}}) \xrightarrow{\ac} F(c^{\SSS^{-1}})
	\end{equation*}
	is a weak equivalence, hence an application of the two-out-of-six property to diagram \Cref{diag:SpectralImpliesLinear} shows that $\sigma_1 \colon \Sigma F(c) \to F(c^{\SSS^{-1}})$ is a weak equivalence.
\end{proof}

Like in the proofs of \Cref{theorem:MainTheoremSimplicial1,theorem:MainTheoremReducedSimplicial1}, we will show that any strictly linear functor is equivalent to a spectral one by carefully analyzing derived left Kan extensions. In our analysis, we will make use of (orthogonal) spectrum objects in $\modN$. Recall that an \emph{(orthogonal) spectrum object} in $\modN$ is a $\Top_*$-functor $\OIndex \to \modN$. The category $\Sp(\modN) := \Fun(\OIndex, \modN)$ of spectrum objects in $\modN$ comes with the usual adjunction
\begin{equation*}\Sigma^\infty : \modN \rightleftarrows \Sp(\modN) : \Omega^\infty,\end{equation*}
where $\Sigma^\infty$ sends an object $N$ of $\modN$ to the spectrum object $\{\Sigma^n N\}_{n \in \bbN}$ and where $\Omega^\infty$ evaluates a spectrum object at the $0$-th level. Moreover, since $\modN$ is a spectral model category, it comes with a second adjunction
\begin{equation*}|-|_\SSS : \Sp(\modN) \rightleftarrows \modN : \Sing_\SSS,\end{equation*}
where the right adjoint $\Sing_\SSS$ sends an object $N$ of $\modN$ to $\{N^{\SSS^{-n}}\}_{n \in \bbN}$ while the left adjoint is given by $|\{M_n\}_{n \in \bbN}|_\SSS = \int^{\OIndex} \SSS^{-n} \otimes M_n$.

\begin{remark}\label{remark:StrictLinearityAsLevelwiseEquivalence}
	If $\Cin$ is a spectral category that admits cotensors by $\SSS^{-1}$, then given a $\Top_*$-functor $F \colon \Cin_* \to \modN$, one can form the spectrum object $\{F(c^{\SSS^{-n}})\}_{n \in \bbN}$. Since $F(c) \cong F(c^{\SSS^{0}}) = \Omega^\infty \{F(c^{\SSS^{-n}})\}_{n \in \bbN}$, we obtain a map
	\begin{equation}\label{equation:LevelwiseSigman}
	\{\Sigma^n F(c)\}_{n \in \bbN} \to \{F(c^{\SSS^{-n}})\}_{n \in \bbN}
	\end{equation}
	by adjunction. It follows by construction that this is levelwise the map $\sigma_n$ of \Cref{construction:Sigma1}. In particular, if $F$ is projectively cofibrant, then by \Cref{remark:Sigma1EquivalenceImpliesSigmanEquivalence} and \Cref{lemma:SuspensionsProjectivelyCofibrantFunctors}, the functor $F$ is strictly linear if and only if the map \Cref{equation:LevelwiseSigman} is a level equivalence for every $c$ in $\Cin$.
\end{remark}

We will endow $\Sp(\modN)$ with the projective model structure and call the weak equivalences of the projective model structure \emph{level equivalences}. The two adjunctions described above are Quillen pairs with respect to the projective model structure on $\Sp(\modN)$ since their right adjoints clearly preserve (trivial) fibrations.

\begin{lemma}\label{lemma:SpectralRealization1}
	There is a natural isomorphism $|\Sigma^\infty N|_\SSS \cong N$.
\end{lemma}

\begin{proof}
	This follows since $N \cong N^{\SSS^0} = \Omega^\infty \Sing_\SSS(N)$.
\end{proof}

\begin{lemma}\label{lemma:SpectralRealization2}
	Let $\Cin_*$ be a small $\Top_*$-category and let $F \colon \Cin_* \to \modN$ and $G \colon \Cin_*^\mathrm{op} \to \Sp$ be $\Top_*$-functors. Then there is a natural isomorphism 
	\begin{equation*}
	\int^{\Cin_*} G \otimes F \cong \left|\left\{\int^{\Cin_*} G_n \otimes F\right\}_{n \in \bbN}\right|_\SSS
	\end{equation*}
	where $G_n$ denotes the functor $\Cin_* \to \Top_*$ that sends $c$ to the $n$-th space of $G(c)$.
\end{lemma}

\begin{proof}
	Note that if $X = \{X_n\}_{n \in \bbN}$ is an orthogonal spectrum, then $\int^{\OIndex} \SSS^{-n} \otimes X_n \cong X$ by the enriched coYoneda lemma (see e.g.\ (3.72) of \cite{Kelly1982BasicConcepts}). In particular, we see that
	\begin{equation*}\int^{\Cin_*} G \otimes F \cong \int^{\Cin_*} \int^{\OIndex} \SSS^{-n} \otimes G_n \otimes F \cong \int^{\OIndex} \SSS^{-n} \otimes \int^{\Cin_*} G_n \otimes F \end{equation*}
	by \Cref{lemma:Fubini}.
\end{proof}

\begin{lemma}\label{lemma:SpectralRealization3}
	Let $\modN$ be a good spectral model category and let $\Cin_*$ be a small $\Top_*$-category. Then for any projectively cofibrant $F \colon \Cin_* \to \modN$ and any pointwise cofibrant $G \colon \Cin_* \to \Sp$, the spectrum object $\{\int^{\Cin_*} G_n \otimes F\}_{n \in \bbN}$ is projectively cofibrant in $\Sp(\modN)$.

\end{lemma}

\begin{proof}
	We need to show that if $\{L_n\}_{n \in \bbN} \to \{K_n\}_{n \in \bbN}$ is levelwise a trivial fibration, then any map $\{\int^{\Cin_*} G_n \otimes F\}_{n \in \bbN} \to \{K_n\}_{n \in \bbN}$ lifts to a map to $\{L_n\}_{n \in \bbN}$. By adjunction, this is equivalent to constructing a lift in
	\begin{equation*}
	\begin{tikzcd}
	 & \int_{\OIndex} L_n^{G_n} \ar[d] \\
	 F \ar[r] \ar[ur,dashed] & \int_{\OIndex} K_n^{G_n},
	\end{tikzcd}
	\end{equation*}
	hence the result follows if we can show that $\int_{\OIndex} L_n^{G_n} \to \int_{\OIndex} K_n^{G_n}$ is pointwise a trivial fibration. In particular, it suffices to show that for any cofibrant orthogonal spectrum $X = \{X_n\}_{n \in \bbN}$, the functor
	\begin{equation*}
		\Sp(\modN) \to \modN; \quad \{N_n\}_{n \in \bbN} \mapsto \int_{\OIndex} N_n^{X_n}
	\end{equation*}
	is right Quillen. This follows since its left adjoint is the functor $\modN \to \Sp(\modN)$ defined by $N \mapsto \{X_n \otimes N\}_{n \in \bbN}$, which is easily seen to be left Quillen.
\end{proof}

Finally, let us mention the following criterion for determining whether $F \colon \Cin_* \to \modN$ is strictly linear.

\begin{lemma}\label{lemma:CriterionStrictLinearity}
	Let $\Cin$ and $\modN$ be as in \Cref{theorem:MainTheoremLinearOrthogonal1}, and let $p \colon \wh \Cin \trivfibarrow \Cin$ be a cofibrant replacement of $\Cin$ in the sense of \Cref{lemma:CofibrantReplacementOfCategories}. Then $F \colon \Cin_* \to \modN$ is strictly linear if and only if there exists a projectively cofibrant $G \colon \wh\Cin_* \to \modN$ such that $F \simeq \Lan_p G$ and
	\begin{equation}\label{equation:CriterionStrictLinearity}
		\{\Sigma^n G \}_{n \in \bbN} \to \left\{\int^{\wh \Cin_*} \wh \Cin(c,-)_n \otimes G(c)\right\}_{n \in \bbN}
	\end{equation}
	is pointwise a level equivalence in $\Sp(\modN)$.
\end{lemma}

\begin{proof}
	Note that $p \colon \wh\Cin_* \to \Cin_*$ is a weak equivalence between $\Top_*$-categories whose hom-spaces are nondegenerately based, hence the goodness of $\modN$ implies that $F \simeq \Lan_p G$ for some projectively cofibrant $G \colon \wh\Cin_* \to \modN$. Since strict linearity is a homotopy invariant property, it suffices to show that $\Lan_p G$ is strictly linear if and only if \Cref{equation:CriterionStrictLinearity} is a level equivalence. To see that this holds, note that
	\begin{equation*}
	\int^{\wh \Cin_*} \Cin(c,d)_n \otimes G(c) \cong \int^{\Cin_*} \Cin_*(c,d^{\SSS^{-n}}) \otimes \Lan_p G(c) \cong \Lan_p G(d^{\SSS^{-n}})
	\end{equation*}
	naturally in $d$. In particular, we obtain the commutative diagram
	\begin{equation*}
	\begin{tikzcd}[column sep = small]
		\{\Sigma^n G(d)\}_{n \in \bbN} \ar[rr] \ar[d,"\sim" rot90] & & \{\int^{\wh \Cin_*} \wh\Cin(c,d)_n \otimes G(c)\}_{n \in \bbN} \ar[d,"\sim" rot90] \\
		\{\Sigma^n \Lan_p G(d)\}_{n \in \bbN} \ar[r] & \{ \Lan_p G(d^{\SSS^{-n}}) \}_{n \in \bbN} \ar[r,"\cong"] & \{\int^{\wh \Cin_*} \Cin(c,d)_n \otimes G(c)\}_{n \in \bbN}
	\end{tikzcd}
	\end{equation*}
	Using that $\wh \Cin \trivfibarrow \Cin$ is a level equivalence on hom-objects, we see that the left-hand vertical arrow is a level equivalence since $\Sigma^n G$ is projectively cofibrant. On the other hand, the right-hand vertical arrow is a level equivalence by \Cref{proposition:FunctorTensorProductInGoodPointedModelCat}. It follows from \Cref{remark:StrictLinearityAsLevelwiseEquivalence} that $\Lan_p G$ is strictly linear if and only if \Cref{equation:CriterionStrictLinearity} is a level equivalence.
\end{proof}

Combining these lemmas, we prove the other direction of \Cref{theorem:MainTheoremLinearOrthogonal1}. Write
\begin{equation*}
L : \Fun(\Cin_*,\modN) \rightleftarrows \Fun(\Cin,\modN) : U
\end{equation*}
for the Quillen pair induced by restriction and left Kan extension along $\Sigma^\infty \Cin_* \to \Cin$, where $\Sigma^\infty \Cin_*$ denotes the spectral category obtained by applying $\Sigma^\infty$ to the hom-spaces of $\Cin_*$.

\begin{proposition}\label{proposition:StrictlyLinearImpliesUnitWE}
	Let $\Cin$ be a small spectral category that admits cotensors by $\SSS^{-1}$, let $\modN$ be a good spectral model category, and assume that the hom-spaces of the underlying $\Top_*$-category of $\Cin$ are nondegenerately based. Then for any strictly linear and projectively cofibrant $\Top_*$-functor $F \colon \Cin_* \to \modN$, the unit $F \to ULF$ is a pointwise weak equivalence.
\end{proposition}

\begin{proof}
	Let $\wh \Cin \trivfibarrow \Cin$ be a cofibrant replacement in the sense of \Cref{lemma:CofibrantReplacementOfCategories}. Restriction and left Kan extension give us a diagram of Quillen pairs
	\begin{equation*}
	\begin{tikzcd}
	{\Fun(\wh\Cin_*,\modN)} \arrow[d, shift right] \arrow[r, shift left,"L'"]          & {\Fun(\wh\Cin,\modN)} \arrow[l, shift left,"U'"] \arrow[d, shift right] \\
	{\Fun(\Cin_*,\modN)} \arrow[u, shift right] \arrow[r, shift left,"L"] & {\Fun(\Cin,\modN),} \arrow[l, shift left,"U"] \arrow[u, shift right]
	\end{tikzcd}
	\end{equation*}
	The vertical adjunctions are Quillen equivalences by the goodness of $\modN$. In particular, by \Cref{lemma:CriterionStrictLinearity}, the proposition follows if we can show that for any projectively cofibrant $G \colon \wh \Cin_* \to \modN$ such that the map \Cref{equation:CriterionStrictLinearity} is a level equivalence, the unit $G \to U' L' G = \int^{\wh \Cin_*} \wh \Cin(c,-) \otimes G(c)$ is a pointwise equivalence. By \Cref{lemma:SpectralRealization1,lemma:SpectralRealization2}, this map is obtained by applying $|-|_\SSS$ to
	\begin{equation*}
	\{\Sigma^n G(d) \}_{n \in \bbN} \to \left\{\int^{\wh \Cin_*} \wh \Cin(c,d)_n \otimes G(c)\right\}_{n \in \bbN}
	\end{equation*}
	This map is a level equivalence by assumption. Since $G$ is projectively cofibrant and $\Cin$ has cofibrant hom-objects, we see that $G(d)$ and hence $\Sigma^\infty G(d)$ are cofibrant. Moreover, $\{\int^{\wh \Cin_*} \wh \Cin(c,d)_n \otimes G(c) \}_{n \in \bbN}$ is cofibrant by \Cref{lemma:SpectralRealization3}, so we conclude that the unit $G \to U'L'G$ is a pointwise equivalence.
\end{proof}

\begin{proof}[Proof of \Cref{theorem:MainTheoremLinearOrthogonal1}]
This follows from \Cref{proposition:SpectralFunctorIsStrictlyLinear} and \Cref{proposition:StrictlyLinearImpliesUnitWE}.
\end{proof}

The \emph{linear model structure} $\Funlin(\Cin_*,\modN)$ is defined, if it exists, as the left Bousfield localization of the projective model structure on $\Fun(\Cin_*,\modN)$ in which the fibrant objects are the strictly linear pointwise fibrant functors. As in the previous cases, we can upgrade \Cref{theorem:MainTheoremLinearOrthogonal1} to a Quillen equivalence.

\begin{theorem}\label{theorem:MainTheoremLinearOrthogonal2}
Let $\Cin$ be a small spectral category that admits cotensors by $\SSS^{-1}$, let $\modN$ be a good spectral model category and assume that the hom-spaces of the underlying $\Top_*$-category of $\Cin$ are nondegenerately based. Then the forgetful functor
\begin{equation*}
	U \colon \Fun(\Cin, \modN) \to \Funlin(\Cin_{*}, \modN)
\end{equation*}
from the projective model structure to the linear model structure is a right Quillen equivalence whenever the linear model structure exists.
\end{theorem}

\begin{proof}
Write $L$ for the left adjoint of the forgetful functor $U$. It follows from the proof of \Cref{theorem:MainTheoremLinearOrthogonal1} that the unit $F \to ULF$ is a weak equivalence for any projectively cofibrant strictly linear functor $F \colon \Cin_* \to \modN$. By \Cref{lemma:BousfieldLocalizationLemma}, $L \dashv U$ is a Quillen equivalence.
\end{proof}

\subsection{Linear functors in the presence of tensors with \texorpdfstring{$\SSS^{-1}$}{S\textasciicircum -1}}\label{ssec:LinearOrthogonalSpectraCaseTensored}

We will now consider the case where the indexing category $\Cin$ has tensors by $\SSS^{-1}$ instead of cotensors. Note that we can't simply dualize \Cref{theorem:MainTheoremLinearOrthogonal1}, since the notion of a good spectral model category is not self-dual. However, we will show that in certain special cases, the proof of \Cref{theorem:MainTheoremLinearOrthogonal1} can be modified in such a way that it works when $\Cin$ has tensors instead of cotensors by $\SSS^{-1}$.

\begin{assumption}\label{assumption:The-tensored-linear-case}
	Let $\Cin$ denote a small spectral category such that the following holds:
	\begin{itemize}[noitemsep]
		\item $\Cin$ admits tensors by $\SSS^{-1}$ and $\SSS^1$.
		\item The hom-spaces of the underlying $\Top_*$-category $\Cin_*$ are nondegenerately based.
		\item For any two objects $c$ and $d$ in $\Cin$, the map
		\[\Cin(c,\SSS^1 \otimes \SSS^{-1} \otimes d) \to \Cin(c, d)\]
		induced by the canonical map $\SSS^1 \otimes \SSS^{-1} \to \SSS^{0}$ is a stable equivalence in $\Sp$.
	\end{itemize}
	Moreover, let $\Din$ be a small spectral category with cofibrant hom-objects and let $\modN$ denote the category $\Fun(\Din,\Sp)$, endowed with any model structure that is a left Bousfield localization of the projective model structure. In particular, $\modN$ is a good spectral model category by \Cref{proposition:GoodSpectralInNature}.
\end{assumption}

\begin{remark}\label{remark:ExistenceTestObjectsN}
	These assumptions imply that $\modN = \Fun(\Din,\Sp)$ is left proper and cellular, hence in particular that it admits left Bousfield localizations at any set of maps by \cite[Thm.\ 4.1.1]{Hirschhorn2003Model}. Moreover, it follows that there exists a set $\mathcal{T} \subset \Ob(\modN)$ of objects such that a map $M \to N$ in $\modN$ is a weak equivalence if and only if for every $T \in \mathcal{T}$, the map of derived mapping spaces $\Map(T,M) \to \Map(T,N)$ is a weak equivalence. An explicit example of such a set is given by
	\[\mathcal{T} = \{\Din(d,-) \otimes \SSS^n \mid d \in \Ob(\Din), n \in \mathbb{Z}\}.\]
\end{remark}

\begin{example}
	Let $\Sp_\mathrm{fin}$ denote the full subcategory of $\Sp$ spanned by those orthogonal spectra $X$ for which $\varnothing \to X$ is a finite composition of pushouts of generating cofibrations. One can show inductively that $\Sp_\mathrm{fin}$ has tensors by $\SSS^{1}$ and $\SSS^{-1}$ and that the hom-spaces of the underlying $\Top_*$-category of $\Sp_\mathrm{fin}$ are nondegenerately based. The fact that $\Sp_\mathrm{fin}(X,\SSS^{-1} \otimes \SSS^{1} \otimes Y) \to \Sp_\mathrm{fin}(X,Y)$ is a stable equivalence follows from the lemma below, so we conclude that $\Cin = \Sp_\mathrm{fin}$ satisfies \cref{assumption:The-tensored-linear-case}.
\end{example}

\begin{lemma}\label{lemma:Mapping-spaces-from-finite-spectrum-derived}
	Let $X$ be a finite orthogonal spectrum and $Y \to Z$ a stable equivalence of orthogonal spectra. Then $\Sp(X,Y) \to \Sp(X,Z)$ is a stable equivalence of orthogonal spectra.
\end{lemma}

\begin{proof}
	Since $X$ is cofibrant, this result holds whenever $Y$ and $Z$ are fibrant in $\Sp$. In particular, it suffices to show that for any orthogonal spectrum $Y$, there exists a fibrant replacement $Y \wearrow Y'$ such that $\Sp(X,Y) \to \Sp(X,Y')$ is a stable equivalence. In orthogonal spectra, an explicit fibrant replacement is given by $Y \wearrow Y' = \hocolim_n \Omega^n Y^{\SSS^{-n}}$. Since $X$ is a finite spectrum, we see that
	$$\Sp(X,\hocolim_n \Omega^n Y^{\SSS^{-n}}) \simeq \hocolim_n \Sp(X,\Omega^n Y^{\SSS^{-n}}) \cong \hocolim_n \Omega^n \Sp(X,Y)^{\SSS^{-n}}.$$
	Since $\Sp(X,Y) \to \hocolim_n \Omega^n \Sp(X,Y)^{\SSS^{-n}}$ is a stable equivalence, the result follows.
\end{proof}

The definition of strict linearity is also slightly more involved than in \Cref{ssec:LinearOrthogonalSpectraCase}.

\begin{definition}\label{definition:StableAndStrictlyLinearFunctors}
	Let $\Cin$ and $\modN = \Fun(\Din,\Sp)$ be as in \Cref{assumption:The-tensored-linear-case}. A $\Top_*$-functor $F \colon \Cin_* \to \modN$ is called \emph{stable} if for every $c \in \Cin$, the map
	\[F(\SSS^1 \otimes \SSS^{-1} \otimes c) \to F(c)\]
	induced by the canonical map $\SSS^1 \otimes \SSS^{-1} \to \SSS^0$ is a weak equivalence in $\modN$. The functor $F$ is called \emph{strictly linear} if moreover for every $c$ in $\Cin$, the composite
	\[\mathbb{L} \Sigma F(\SSS^{-1} \otimes c) \to \Sigma F(\SSS^{-1} \otimes c) \to F(c)\]
	is a weak equivalence in $\modN$. This map is constructed analogously to \Cref{construction:Sigma1}.
\end{definition}

The main result is now as follows.

\begin{theorem}\label{theorem:MainTheoremLinearOrthogonal1Tensored}
	Let $\Cin$ and $\modN = \Fun(\Din,\Sp)$ be as in \Cref{assumption:The-tensored-linear-case} and let $F \colon \Cin_* \to \modN$ be a $\Top_*$-functor. Then $F$ is equivalent to a spectral functor if and only if it is strictly linear.
\end{theorem}

\begin{remark}
	The author suspects that this theorem should hold under more general assumptions on $\modN$ than those of \Cref{assumption:The-tensored-linear-case}. However, many $\Sp$-enriched model categories considered in practice are of the kind described in \Cref{assumption:The-tensored-linear-case}.
\end{remark}

The proof of \Cref{theorem:MainTheoremLinearOrthogonal1Tensored} in the ``only if'' direction uses the following lemma.

\begin{lemma}\label{lemma:RepresentablesStrictlyLinearAndStable}
	Let $\Cin$ be as in \Cref{assumption:The-tensored-linear-case}. For any $c$ in $\Ob(\Cin)$, the functor $\Cin(c,-) \colon \Cin \to \Sp$ is strictly linear.
\end{lemma}

\begin{proof}
	For ease of notation, we write $F$ for $\Cin(c,-)$. The fact that $F$ is stable is part of \Cref{assumption:The-tensored-linear-case}. For strict linearity, note that $\mathbb{L}\Sigma F(\SSS^{-1} \otimes d) \simeq \Sigma F(\SSS^{-1} \otimes d)$. Since $F$ is an $\Sp$-functor, note that for any spectrum $X$ such that $\Cin$ admits tensors by $X$, there is a canonical natural map $X \otimes F(d) \to F(X \otimes d)$. In particular, we obtain a commutative diagram
	\begin{equation*}
		\begin{tikzcd}
			\SSS^{-1} \otimes \SSS^1 \otimes \SSS^{-1} \otimes F(d) \ar[r] \ar[d,"\id \otimes f \otimes \id"] & \SSS^{-1} \otimes \SSS^1 \otimes F(\SSS^{-1} \otimes d) \ar[r, "f \otimes \id"] \ar[d] & F(\SSS^{-1} \otimes d) \\
			\SSS^{-1} \otimes F(d) & \SSS^{-1} \otimes F(\SSS^1 \otimes \SSS^{-1} \otimes d) \ar[r] \ar[l,"\id \otimes F(f \otimes \id)","\sim"'] & F(\SSS^{-1} \otimes \SSS^1 \otimes \SSS^{-1} \otimes d) \ar[u,"F(f \otimes \id)"',"\sim" rot90]
		\end{tikzcd}
	\end{equation*}
	where $f$ denotes the map $\SSS^1 \otimes \SSS^{-1} \to \SSS^0$. The stability of $F$ implies that the rightmost vertical map is a weak equivalence as well as the bottom left horizontal map. Since $f$ is a stable equivalence between cofibrant objects, we also see that the top right horizontal map and the leftmost vertical map are stable equivalences. By the two-out-of-six property, all maps are stable equivalences. Now note the map $\Sigma F(\SSS^{-1} \otimes c) \to F(c)$, which we need to show is a stable equivalence, agrees with the composite
	\[\SSS^1 \otimes F(\SSS^{-1} \otimes c) \to F(\SSS^1 \otimes \SSS^{-1} \otimes c) \to F(c).\]
	After tensoring with $\SSS^{-1}$, these maps become the middle vertical and bottom left horizontal map in the diagram above, hence weak equivalences. We conclude that $\Sigma F(\SSS^{-1} \otimes c) \to F(c)$ is a weak equivalence.
\end{proof}

\begin{proposition}\label{lemma:SpectralFunctorIsStrictlyLinearTensoredCase}
	Let $\Cin$ and $\modN$ be as in \Cref{assumption:The-tensored-linear-case} and let $F \colon \Cin \to \modN$ be a spectral functor. Then $F$ is strictly linear.
\end{proposition}

\begin{proof}
	Without loss of generality, assume that $F$ is projectively cofibrant. By the coYoneda lemma \cite[(3.71)]{Kelly1982BasicConcepts}, we may write
	\[F(d) \cong \int^\Cin \Cin(c,d) \otimes F(c),\]
	while by \Cref{lemma:SuspensionsProjectivelyCofibrantFunctors} we have $\mathbb{L} \Sigma F(d) \simeq \Sigma F(d)$, so it suffices to show that
	\[\int^\Cin \Cin(c,\SSS^{1} \otimes \SSS^{-1} \otimes d) \otimes F(c) \to \int^\Cin \Cin(c,d) \otimes F(c)\]
	and
	\[\Sigma \int^\Cin \Cin(c,\SSS^{-1} \otimes d) \otimes F(c) \cong \int^\Cin \Sigma\Cin(c,\SSS^{-1} \otimes d) \otimes F(c) \to \int^\Cin \Cin(c,d) \otimes F(c)\]
	are weak equivalences. This follows from \Cref{proposition:FunctorTensorProductInGoodSpectralModelCat} and \Cref{lemma:RepresentablesStrictlyLinearAndStable}. 
\end{proof}

For the converse, note that the only place in the proof of \Cref{proposition:StrictlyLinearImpliesUnitWE} where cotensors by $\SSS^{-1}$ are used is in the proof of \Cref{lemma:CriterionStrictLinearity}: there they are used to show that
\[\Sigma^n F \to \int^{\Cin_*} \Cin(c,-)_n \otimes F(c)\]
is a weak equivalence when $F$ is strictly linear and cofibrant. In the current setting this need not be true, so we will use a slight variation on this result instead. Let $\Funst(\Cin_*,\modN)$ denote the left Bousfield localization of the projective model structure on $\Fun(\Cin_*,\modN)$ in which the fibrant objects are precisely the projectively fibrant functors that are stable in the sense of \Cref{definition:StableAndStrictlyLinearFunctors}. This model structure exists since it is the left Bousfield localization of the projective model structure at the maps
\[\{\Cin_*(c,-) \otimes T \to \Cin_*(\SSS^{1} \otimes \SSS^{-1} \otimes c, -) \otimes T \mid c \in \Ob(\Cin), T \in \mathcal{T}\},\]
where $\mathcal{T}$ is as in \Cref{remark:ExistenceTestObjectsN}.

\begin{lemma}\label{lemma:MapUnitSpectralTensoredCaseUnrealizedEquivalence}
	Let $F \colon \Cin_* \to \modN$ be projectively cofibrant and strictly linear. Then
	\begin{equation}\label{equation:MapUnitSpectralTensoredCaseUnrealizedLevel}
		\Sigma^n F \to \int^{\Cin_*} \Cin(c,-)_n \otimes F(c)
	\end{equation}
	is a weak equivalence in $\Funst(\Cin_*,\modN)$ for every $n \geq 0$.
\end{lemma}

\begin{proof}
	For $n = 0$ this map is an isomorphism, so suppose that $n \geq 1$. We will write $q$ and $w$ for the functors
	\[\SSS^{-n} \otimes - \colon \Cin \to \Cin \quad \text{and} \quad \SSS^n \otimes - \colon \Cin \to \Cin,\]
	respectively. Precomposition then gives right Quillen functors
	\[q^*, w^* \colon \Fun(\Cin_*,\modN) \to \Fun(\Cin_*,\modN). \]
	It is clear that $q^*$ and $w^*$ take stable functors to stable functors and strictly linear functors to strictly linear ones. In particular, they preserve the fibrant objects of $\Fun^\mathrm{st}(\Cin_*,\modN)$. By Corollary A.2 and Remark A.3 of \cite{Dugger2001Replacing}, they are right Quillen with respect to the stable model structure. Since the map $\SSS^1 \otimes \SSS^{-1} \to \SSS^0$ induces a natural weak equivalence $(q^*w^*F)(c) \cong (w^*q^*F)(c) \wearrow F(c)$ for any stable $F \colon \Cin_* \to \modN$, we see that on the homotopy category of $\Fun^\mathrm{st}(\Cin_*,\modN)$, the right derived functors of $q^*$ and $w^*$ are inverse to each other. In particular, $q^*$ and $w^*$ are right Quillen equivalences.
	
	Note that since $\Cin(c,-)_n \cong \Cin_*(\SSS^{-n} \otimes c, -)$, the target of the map \Cref{equation:MapUnitSpectralTensoredCaseUnrealizedLevel} is simply $q_!F$. Moreover, this map arises from a natural transformation $\phi \colon \Sigma^n \Rightarrow q_!$. Using the counit $\varepsilon$ of $q_! \dashv q^*$ we can form the composite
	\begin{equation}\label{equation:CommutativeTriangleNaturalTransformations}
		\Sigma^n q^* \xRightarrow{\phi_{q^*}} q_! q^* \xRightarrow{\varepsilon} \id
	\end{equation}
	which, for a functor $G$, is simply the map $\Sigma G(\SSS^{-1} \otimes c) \to G(c)$ from \Cref{definition:StableAndStrictlyLinearFunctors} iterated $n$ times.
	
	Now note that $F$ is equivalent to a functor of the form $q^*G$ with $G$ strictly linear, for example $G = w^*F$. In particular, it suffices to show that
	\begin{equation}\label{equation:FinalMapSpectralProofTensoredCase}
		\Sigma^n \wh{q^*G} \to q_! \wh{q^*G}
	\end{equation}
	is a weak equivalence in $\Fun^\mathrm{st}(\Cin_*,\modN)$, where $\wh{q^*G}$ denotes a projectively cofibrant replacement of $q^*G$.
	The natural transformations \Cref{equation:CommutativeTriangleNaturalTransformations} give us a commutative triangle
	\begin{equation*}
		\begin{tikzcd}
			\Sigma^n \wh{q^*G} \ar[r] \ar[d] & q_! \wh{q^*G} \ar[dl] \\
			G &
		\end{tikzcd}
	\end{equation*}
	The vertical map is the one appearing in the definition of strict linearity (iterated $n$ times), so since $\wh{q^*G}$ is projectively cofibrant and strictly linear, it is a pointwise weak equivalence. Moreover, the diagonal map is the derived counit of $q_! \dashv q^*$, so it is a weak equivalence in $\Fun^\mathrm{st}(\Cin_*,\modN)$. We conclude that \Cref{equation:FinalMapSpectralProofTensoredCase} is a weak equivalence in $\Fun^\mathrm{st}(\Cin_*,\modN)$.
\end{proof}

Now let $p \colon \wh \Cin \trivfibarrow \Cin$ be a cofibrant replacement in the sense of \Cref{lemma:CofibrantReplacementOfCategories} and let
\[p_! : \Fun(\wh \Cin_*,\modN) \rightleftarrows \Fun(\Cin_*,\modN) : p^*\]
be the adjunction given by restriction and left Kan extension along $p$. Note that this is a Quillen equivalence by the goodness of $\modN$ and the fact that $\wh \Cin_* \to \Cin_*$ is a weak equivalence of $\Top_*$-categories. We call a $\Top_*$-functor $F \colon \wh \Cin_* \to \modN$ \emph{stable} if $\mathbb{L}p_! F \colon \Cin_* \to \modN$ is stable. Let $\Funst(\wh\Cin_*,\modN)$ again denote the left Bousfield localization of the projective model structure in which the fibrant objects are the stable functors that are projectively fibrant. To see that it exists, note that since $\wh \Cin_*(c,d) \to \Cin_*(c,d)$ is a trivial fibration for any $c,d \in \Ob(\Cin)$, we can choose for every $c \in \Ob(\Cin)$ a map $\wh f_c \in \wh \Cin_*(\SSS^1 \otimes \SSS^{-1} \otimes c, c)$ such that $p(\wh f_c)$ is the canonical map $\SSS^1 \otimes \SSS^{-1} \otimes c \to c$ in $\Cin_*$. Then $\Funst(\wh\Cin_*,\modN)$ is obtained as the left Bousfield localization of the projective model structure at the set of maps
\[\left\{\wh\Cin_*(c,-) \otimes T \xrightarrow{\wh f_c^*} \wh\Cin_*(\SSS^{1} \otimes \SSS^{-1} \otimes c, -) \otimes T \,\middle\vert\, c \in \Ob(\Cin), T \in \mathcal{T}\right\}.\]
Note that $\Funst(\wh \Cin_*,\modN) \rightleftarrows \Funst(\Cin_*,\modN)$ is a Quillen equivalence by \cite[Thm.\ 3.3.20(1)]{Hirschhorn2003Model}.

\begin{lemma}\label{lemma:DerivedUnitSpectralTensoredCase}
	Let $F \colon \wh \Cin_* \to \modN$ be projectively cofibrant and suppose that $\Lan_p F \colon \Cin_* \to \modN$ is strictly linear. Then
	\[F \to \int^{\wh \Cin_*} \wh \Cin(c,-) \otimes F(c) \]
	is a pointwise weak equivalence.
\end{lemma}

\begin{proof}
	Consider the commutative square
	\[\begin{tikzcd}
		\Sigma^n F \ar[rr] \ar[d] & & \int^{\wh \Cin_*} \wh\Cin(c,-)_n \otimes F(c) \ar[d] \\
		\Sigma^n \Lan_p F \ar[r] & \int^{\Cin_*} \Cin(c,-)_n \otimes \Lan_p F \ar[r,"\cong"] & \int^{\wh \Cin_*} \Cin(c,-)_n \otimes F(c)
	\end{tikzcd}\]
	Since $\wh\Cin(c,d) \trivfibarrow \Cin(c,d)$ is a level trivial fibration for any $c,d \in \Ob(\Cin)$, we deduce from \Cref{proposition:FunctorTensorProductInGoodPointedModelCat} that the rightmost vertical map is a pointwise equivalence. Moreover, the leftmost vertical map is a pointwise equivalence since $\Sigma^n F$ is projectively cofibrant, while the bottom horizontal map is a weak equivalence in $\Funst(\Cin_*,\modN)$ by \Cref{lemma:MapUnitSpectralTensoredCaseUnrealizedEquivalence}. In particular,
	\[\Sigma^\infty F = \{\Sigma^n F\}_{n \in \bbN} \to \left\{\int^{\wh \Cin_*} \wh\Cin(c,-)_n \otimes F(c)\right\}_{n \in \bbN}\]
	is a weak equivalence in the projective model structure on $\Sp(\Funst(\wh\Cin_*,\modN))$. Let us write $G$ for the target of this map. We see that $\Sigma^\infty F$ is projectively cofibrant since $F$ is. On the other hand, $G$ is pointwise cofibrant by \Cref{lemma:SpectralRealization3}, in the sense that $G(d)$ is cofibrant in $\Sp(\modN)$ for any $d \in \Ob(\wh \Cin_*)$. Let us factor $\Sigma^\infty F \to G$ as a cofibration $\Sigma^\infty F \cofarrow \wh G$ followed by a pointwise trivial fibration $\wh G \trivfibarrow G$. Then $\Sigma^\infty F \cofarrow \wh G$ is a trivial cofibration in the projective model structure on \(\Sp(\Funst(\wh\Cin_*,\modN))\), hence $|\Sigma^\infty F|_\SSS \to |G|_\SSS$ is a trivial cofibration in $\Funst(\wh\Cin_*,\modN)$. Moreover, $\wh G$ and $G$ are both pointwise cofibrant, hence $|\wh G|_\SSS \to |G|_\SSS$ is a pointwise weak equivalence. We deduce from \Cref{lemma:SpectralRealization1,lemma:SpectralRealization2} that
	\begin{equation}\label{equation:SomeMap}
		F \cong |\Sigma^\infty F|_\SSS \to |G|_\SSS \cong \int^{\wh \Cin_*} \wh \Cin(c,-) \otimes F(c)
	\end{equation}
	is a weak equivalence in $\Funst(\wh\Cin_*,\modN)$. It follows by \Cref{lemma:SpectralFunctorIsStrictlyLinearTensoredCase} that the functor $\int^{\wh \Cin_*} \wh \Cin(c,-) \otimes F(c)$ is strictly linear, so in particular stable. By \cite[Thm.\ 3.2.13(1)]{Hirschhorn2003Model}, the map \Cref{equation:SomeMap} is a pointwise weak equivalence.
\end{proof}

As in the proof of \Cref{proposition:StrictlyLinearImpliesUnitWE}, we deduce from this that any strictly linear functor is naturally equivalent to a spectral functor.

\begin{proof}[Proof of \Cref{theorem:MainTheoremLinearOrthogonal1Tensored}]
	It was proved in \Cref{lemma:SpectralFunctorIsStrictlyLinearTensoredCase} that any spectral functor is strictly linear. The other direction follows exactly as in \Cref{proposition:StrictlyLinearImpliesUnitWE}: in the square of Quillen pairs
	\begin{equation*}
		\begin{tikzcd}
			{\Fun(\wh\Cin_*,\modN)} \arrow[d, shift right] \arrow[r, shift left]          & {\Fun(\wh\Cin,\modN)} \arrow[l, shift left] \arrow[d, shift right] \\
			{\Fun(\Cin_*,\modN)} \arrow[u, shift right] \arrow[r, shift left] & {\Fun(\Cin,\modN),} \arrow[l, shift left] \arrow[u, shift right]
		\end{tikzcd}
	\end{equation*}
	the vertical adjunctions are Quillen equivalences by the goodness of $\modN$, so it follows from \Cref{lemma:DerivedUnitSpectralTensoredCase} that for any $F \colon \Cin_* \to \modN$ that is projectively cofibrant and strictly linear, the unit
	\[F \to \int^{\Cin_{*}} \Cin(c,-) \otimes F(c)\]
	is a pointwise weak equivalence.
\end{proof}

As always, this result can be upgraded to a Quillen equivalence. Let $\Funlin(\Cin_*,\modN)$ denote the left Bousfield localization of the projective model structure in which the fibrant objects are the projectively fibrant functors that are strictly linear. Note that this model structure exists since it can be constructed by localizing $\Funst(\Cin_*,\modN)$ at the set of maps
\[\{\Sigma \Cin_*(c,-) \otimes T \to \Cin_*(\SSS^{-1} \otimes c, -) \otimes T \mid c \in \Ob(\Cin), T \in \mathcal{T}\}.\]
Namely, a projectively fibrant stable functor $F$ is local with respect to this set of maps if and only if
$F(\SSS^{-1} \otimes c) \to \Omega F(c)$
is a weak equivalence for every $c \in \Ob(\Cin)$, which by adjunction is equivalent to $F$ being strictly linear.

\begin{theorem}\label{theorem:MainTheoremLinearOrthogonal2Tensored}
	The adjunction
	\[L : \Funlin(\Cin_*,\modN) \rightleftarrows \Fun(\Cin,\modN) : U\]
	given by restriction and left Kan extension along $\Sigma^\infty \Cin_* \to \Cin$ is a Quillen equivalence.
\end{theorem}

\begin{proof}
	It follows from the proof of \Cref{theorem:MainTheoremLinearOrthogonal1Tensored} that the unit $F \to ULF$ is a weak equivalence for any projectively cofibrant strictly linear functor $F \colon \Cin_* \to \modN$. The result now follows from \Cref{lemma:BousfieldLocalizationLemma}.
\end{proof}

By way of illustration, we will now combine our main results to describe which functors $\Sp_\mathrm{fin} \to \Sp$ are naturally equivalent to spectral functors. We will call an unenriched functor $F \colon \Sp_{\mathrm{fin}} \to \modN$ \emph{linear} if it is reduced and if for any finite orthogonal spectrum $X$, the pushout square
\begin{equation}\label{equation:SuspensionPushout}
	\begin{tikzcd}
		X \ar[r] \ar[d] & CX \ar[d] \\
		* \ar[r] & \Sigma X,
	\end{tikzcd}
\end{equation}
where $CX$ is the reduced cone on $X$, is taken to a homotopy pullback square by $F$.

\begin{proposition}
	An unenriched functor $\Sp_{\mathrm{fin}} \to \Sp$ is naturally equivalent to a spectral functor if and only if it preserves stable equivalences and is linear.
\end{proposition}

\begin{proof}
	We will show below that a $\Top_*$-functor $\Sp_\mathrm{fin} \to \Sp$ that preserves stable equivalences is linear in the sense above if and only if it is strictly linear in the sense of \cref{definition:StableAndStrictlyLinearFunctors}. First, let us show how to deduce the proposition from this. By \cite[(3.71)]{Kelly1982BasicConcepts}, any spectral functor $F \colon \Sp_{\mathrm{fin}} \to \Sp$ can be written as an enriched coend $\int \Sp_\mathrm{fin}(X,-) \otimes F(X)$, hence it follows from \cref{lemma:Mapping-spaces-from-finite-spectrum-derived} and \cref{proposition:FunctorTensorProductInGoodSpectralModelCat} that $F$ preserves stable equivalences. Since $F$ is strictly linear by \cref{proposition:SpectralFunctorIsStrictlyLinear}, it is also linear in the sense above.
	
	Conversely, suppose that $F$ is linear and preserves stable equivalences. Then $F$ is a reduced homotopy functor, so by \cref{theorem:MainTheoremTopological} and \cref{theorem:MainTheoremReducedTopological}, we may assume that $F$ is $\Top_*$-enriched. In particular, $F$ is strictly linear and hence equivalent to a spectral functor by \cref{theorem:MainTheoremLinearOrthogonal1Tensored}.
	
	What remains is to show the claim about linearity and strict linearity made at the start of the proof. If $F \colon \Sp_\mathrm{fin} \to \Sp$ is a $\Top_*$-functor, then it is linear precisely if $F(X) \to \mathbb{R} \Omega F(\Sigma X)$ is a weak equivalence for every finite spectrum $X$, which is the case precisely if the assembly map $\mathbb{L} \Sigma F(X) \to F(\Sigma X)$ is strictly linear. We can factor $\Sigma F(\SSS^{-1} \otimes X) \to F(X)$ as
	\[\begin{tikzcd}
		\Sigma F(\SSS^{-1} \otimes X) \ar[r] \ar[d] & F(X) \\
		F(\Sigma \SSS^{-1} \otimes X) \cong F(\SSS^1 \otimes \SSS^{-1} \otimes X) \ar[ur]
	\end{tikzcd}\]
	where the vertical map is the assembly map and the diagonal map is the one appearing in \cref{definition:StableAndStrictlyLinearFunctors}. In particular, if $F$ preserves stable equivalences, then $F$ is strictly linear if and only if for any finite spectrum $X$, the assembly map $\mathbb{L} \Sigma F(\SSS^{-1} \otimes X) \to F(\Sigma \SSS^{-1} \otimes X)$ is a weak equivalence. Since any finite orthogonal spectrum is stably equivalent to one of the form $\SSS^{-1} \otimes X$, it follows that a $\Top_*$-functor that preserves stable equivalences is linear if and only if it is strictly linear.
\end{proof}

\subsection{(Symmetric) multilinear functors}\label{ssec:MultilinearFunctors}

We will now extend the results of \Cref{ssec:LinearOrthogonalSpectraCase,ssec:LinearOrthogonalSpectraCaseTensored} to (symmetric) multifunctors. Recall the definition of a (symmetric) multifunctor from \Cref{definition:Multifunctor,definition:SymmetricMultifunctor}.

Let $\Cin$ denote a spectral category and $\modN$ a spectral model category. Given a multifunctor $F \colon  \Cin_*^{\wedge n} \to \modN$ and $1 \leq i \leq n$, any $(n-1)$-tuple $\vec{c} = (c_1,\ldots,c_{i-1},c_{i+1},\ldots,c_n)$ of objects in $\Cin$ defines a functor
\[F^i_{\vec{c}} \colon \Cin_* \to \modN; \quad c \mapsto F(c_1,\ldots,c_{i-1},c,c_{i+1},\ldots,c_n).\]

\begin{definition}
	Let $\Cin$ be a spectral category, $\modN$ a spectral model category, and $F \colon \Cin_*^{\wedge n} \to \modN$ a (possibly symmetric) multifunctor. If $\Cin$ admits cotensors by $\SSS^{-1}$, then we say that $F$ is \emph{strictly multilinear} if it is strictly linear in every variable separately; that is, for every $1 \leq i \leq n$ and every $(n-1)$-tuple $\vec{c}$ of objects in $\Cin$, the functor $F^i_{\vec c}$ is strictly linear in the sense of \Cref{definition:LinearOrthogonal}. If $\Cin$ instead admits tensors by $\SSS^{-1}$ and $\SSS^1$, then we will call $F$ \emph{stable} or \emph{strictly multilinear} if it is stable or strictly linear in each variable separately in the sense of \Cref{definition:StableAndStrictlyLinearFunctors}.
\end{definition}

\begin{assumption}\label{assumption:The-multilinear-case}
	Throughout this section, $\Cin$ denotes a small spectral category and $\modN$ a spectral model category such that one of the following is satisfied:
	\begin{enumerate}[(1)]
		\item\label{item:assumption1} $\Cin$ admits cotensors by $\SSS^{-1}$, the underlying $\Top_*$-category of $\Cin$ has nondegenerately based hom-spaces and $\modN$ is a good spectral model category.
		\item\label{item:assumption2} $\Cin$ and $\modN$ satisfy \Cref{assumption:The-tensored-linear-case}.
	\end{enumerate}
\end{assumption}

We wish to show that strictly multilinear functors $\Cin_*^{\wedge n} \to \modN$ are equivalent to spectral functors $\Cin^{\otimes n} \to \modN$. However, one has to be careful since the tensor product of two hom-objects of $\Cin$ might not be the derived tensor product, hence $\Cin^{\otimes n}$ might have the wrong homotopy type.\footnote{To the author's knowledge, the question of whether the tensor product of orthogonal spectra is fully homotopical or not is open (cf.\ \cite[Question 1.5]{SagaveSchwede2021Convolution}). If it were fully homotopical, then this section could be significantly simplified.} For that reason, we need to work with a cofibrant replacement $\wh \Cin$ of $\Cin$ and consider spectral functors $\wh \Cin^{\otimes n} \to \modN$ instead.

Throughout the rest of this section, we fix a cofibrant replacement $q \colon \wh \Cin \trivfibarrow \Cin$ in the sense of \Cref{lemma:CofibrantReplacementOfCategories}. On underlying $\Top_*$-categories, we obtain a weak equivalence $p \colon \wh \Cin_* \to \Cin_*$ of $\Top_*$-categories whose hom-spaces are nondegenerately based, so by the goodness of $\modN$, we obtain a Quillen equivalence
\begin{equation}\label{equation:MultivariateRestrictionKanExtenstion}
	p_! : \Fun(\wh\Cin_*^{\wedge n}, \modN) \rightleftarrows \Fun(\Cin_*^{\wedge n}, \modN) : p^*.
\end{equation}
Let us write
\begin{equation}\label{equation:MultivariateSpectralization}
	L : \Fun(\wh\Cin_*^{\wedge n}, \modN) \rightleftarrows \Fun(\wh \Cin^{\otimes n}, \modN) : U
\end{equation}
for the adjunction obtained by restriction and left Kan extension along $(\Sigma^\infty \wh \Cin_*)^{\wedge n} \to \wh \Cin^{\otimes n}$.
The main result of this section is now the following.

\begin{theorem}\label{theorem:MainTheoremMultilinear1}
	Let $\Cin$ and $\modN$ be as in \Cref{assumption:The-multilinear-case} and let a multifunctor $F \colon \Cin_*^{\wedge n} \to \modN$ be given. Then $F$ is strictly multilinear if and only if there exists a spectral functor $G \colon \wh \Cin^{\otimes n} \to \modN$ such that $p^*F \simeq UG$.
\end{theorem}

We will first do the ``if'' direction.

\begin{proposition}\label{lemma:SpectralMultifunctorIsStrictlyLinear}
	Let $\Cin$ and $\modN$ be as in \Cref{assumption:The-multilinear-case} and let a multifunctor $F \colon \Cin_*^{\wedge n} \to \modN$ be given such that $p^*F \simeq UG$ for some spectral multifunctor $G \colon \wh \Cin^{\otimes n} \to \modN$. Then $F$ is strictly multilinear.
\end{proposition}

\begin{proof}
	Given $1 \leq i \leq n$ and an $(n-1)$-tuple $\vec c = (c_1,\ldots,c_{i-1},c_{i+1},\ldots,c_n)$ of objects in $\Ob(\Cin)$, let $\iota^i_{\vec c}$ denote the functor
	\[\Cin \to \Cin^{\otimes n}; \quad c \mapsto (c_1,\ldots,c_{i-1},c,c_{i+1},\ldots,c_n).\]
	The analogous functors for $\Cin_*$, $\wh \Cin$ and $\wh \Cin_*$ will also be denoted $\iota^i_{\vec c}$. By definition, $F$ is strictly multilinear if and only if all the functors $(\iota^i_{\vec c})^* F$ are strictly linear. Write $q$ and $p$ for the functors $\wh \Cin \trivfibarrow \Cin$ and $\wh \Cin_* \to \Cin_*$, respectively. Since $p^*F \simeq UG$, we see that $p^*(\iota^i_{\vec c})^* F \simeq U (\iota^i_{\vec c})^* G$. Since
	\[q_! : \Fun(\wh\Cin, \modN) \rightleftarrows \Fun(\Cin, \modN) : q^*\]
	is a Quillen equivalence by the goodness of $\modN$, this implies that $(\iota^i_{\vec c})^*F \simeq UG'$ for some $G' \colon \Cin \to \modN$. By \Cref{proposition:SpectralFunctorIsStrictlyLinear,lemma:SpectralFunctorIsStrictlyLinearTensoredCase}, we conclude that $(\iota^i_{\vec c})^* F$ is strictly linear.
\end{proof}

For the other direction, we will use the following. Given an $n$-tuple $\vec m = (m_1,\ldots,m_n)$ of natural numbers, we shall write $\Sigma^{\vec m}$ for the suspension $\Sigma^{m_1} \cdots \Sigma^{m_n} = \Sigma^{m_1 + \cdots + m_n}$. Let $\Funst(\Cin_*^{\wedge n}, \modN)$ denote the left Bousfield localization of the projective model structure on $\Fun(\Cin_*^{\wedge n},\modN)$ in which the fibrant objects are the projectively fibrant functors that are stable in each variable separately.

\begin{lemma}
	If item \ref{item:assumption1} of \Cref{assumption:The-multilinear-case} holds and $F \colon \Cin_*^{\wedge n} \to \modN$ is strictly multilinear and projectively cofibrant, then the map
	\begin{equation}\label{equation:MapUnitSpectralMultifunctorUnrealizedLevel}
		\Sigma^{\vec m} F(\vec d) \to \int\limits^{\vec c \in \Cin_*^{\wedge n}} \Cin(c_1,d_1)_{m_1} \otimes \cdots \otimes \Cin(c_n,d_n)_{m_n} \otimes \cdots \otimes F(\vec c)
	\end{equation}
	is a weak equivalence for every $n$-tuple $\vec m$ of natural numbers and every $n$-tuple $\vec d$ of objects in $\Cin$. If instead \ref{item:assumption2} holds and $F$ is strictly multilinear and projectively cofibrant, then this map is a weak equivalence in $\Fun^\mathrm{st}(\Cin_*^{\wedge n},\modN)$.
\end{lemma}

\begin{proof}
	This is analogous to the proofs of \Cref{lemma:CriterionStrictLinearity,lemma:MapUnitSpectralTensoredCaseUnrealizedEquivalence}.
\end{proof}

In the proof of \Cref{theorem:MainTheoremMultilinear1}, we will make use of \emph{$n$-fold spectrum objects} in a given $\Top_*$-category $\modM$; i.e., functors $\OIndex^{\wedge n} \to \modM$. We write $\Sp^n(\modM)$ for the category
\[\Sp^n(\modM) = \Fun(\OIndex^{\wedge n}, \modM)\]
of $n$-fold spectrum objects in $\modM$. If $\modM$ admits tensors by $S^1$, then we obtain an adjunction
\[ \Sigma^\infty : \modM \rightleftarrows \Sp^n(\modM) : \Omega^\infty \]
where $\Sigma^\infty M = \{\Sigma^{\vec m} M\}_{\vec m \in \bbN^n}$. If $\modM$ is a (co)complete spectral category, then we moreover obtain the adjunction
\begin{equation*}|-|^n_\SSS : \Sp^n(\modN) \rightleftarrows \modN : \Sing^n_\SSS,\end{equation*}
where \[\Sing_\SSS^n M = \{M^{\SSS^{-m_1} \otimes \cdots \otimes \SSS^{-m_n}}\}_{\vec m \in \bbN^n} \quad \text{and} \quad |\{N_{\vec m}\}_{\vec m \in \bbN^n}|^n_\SSS = \int^{\OIndex^{\wedge n}} \SSS^{-m_1} \otimes \cdots \otimes \SSS^{-m_n} \otimes N.\]
If $\modM$ is a spectral model category and $\Sp^n(\modM)$ is endowed with the projective model structure, then both of these adjunctions are Quillen pairs.

\begin{proof}[Proof of \Cref{theorem:MainTheoremMultilinear1}]
	One direction was proved in \cref{lemma:SpectralMultifunctorIsStrictlyLinear}. For the other direction, let a strictly multilinear functor $F \colon \Cin_*^{\wedge n} \to \modN$ be given. By the Quillen equivalence \Cref{equation:MultivariateRestrictionKanExtenstion}, there exists a projectively cofibrant $G \colon \wh \Cin_*^{\wedge n} \to \modN$ such that $p_!G \simeq F$ and $G \simeq p^*F$. In particular, it suffices to show that $G \to ULG$ is a pointwise weak equivalence, where $L \dashv U$ is the Quillen pair \Cref{equation:MultivariateSpectralization}. To this end, note that it follows as in \Cref{lemma:SpectralRealization1,lemma:SpectralRealization2} that $G \to ULG$ agrees with the map obtained by applying $|-|^n_\SSS$ to
	\[\Sigma^\infty G \cong \{\Sigma^{\vec m} G\}_{\vec{m} \in \bbN^n} \to \left\{\int\limits^{\vec c \in \wh \Cin_*^{\wedge n}} \wh \Cin(c_1,-)_{m_1} \otimes \cdots \otimes \wh \Cin(c_n,-)_{m_n} \otimes \cdots \otimes G(\vec c)\right\}_{\vec m \in \bbN^n}\]
	It now follows as in \Cref{proposition:StrictlyLinearImpliesUnitWE,lemma:DerivedUnitSpectralTensoredCase} that $G \to ULG$ is a weak equivalence.
\end{proof}

The analogue of \Cref{theorem:MainTheoremMultilinear1} for symmetric multifunctors is now deduced in the same way as \Cref{theorem:MainTheoremSymmetricReduced}.

\begin{theorem}\label{theorem:MainTheoremSymmetricMultilinear1}
	Let $\Cin$ and $\modN$ be as in \Cref{assumption:The-multilinear-case} and let a symmetric multifunctor $F \colon \Cin_*^{\wedge n} \to \modN$ be given. Then $F$ is strictly multilinear if and only if there exists a symmetric spectral multifunctor $G \colon \wh \Cin^{\otimes n} \to \modN$ such that $p^*F$ is naturally equivalent to $UG$ (as symmetric multifunctors).
\end{theorem}

\begin{remark}
	Both these results can be strengthened to zigzags of Quillen equivalences between suitable model categories, analogous to \Cref{theorem:MainTheoremLinearOrthogonal2,theorem:MainTheoremLinearOrthogonal2Tensored}. We leave this to the reader to work out.
\end{remark}